\def\editmode{0}

\def\bibfilenames{bibman_refs}

\def\singlenarrowcol{0}
\if\singlenarrowcol0
    \documentclass[conference]{IEEEtran}

    \IEEEoverridecommandlockouts

\else
    \documentclass[onecolumn]{IEEEtran}
    \usepackage[pass]{geometry}
    \newgeometry{textwidth=252pt, textheight=672pt}
    
\fi

\def\spsformat{0}

\usepackage{hyperref} 
  
\usepackage{include_v7}

\usepackage{algorithmic}   
\usepackage[ruled,vlined,resetcount,noend,linesnumbered]{algorithm2e}

\newenvironment{changes}{}{}

\begin{document}

\title{Aerial Base Station Placement via  \\Propagation Radio  Maps}

\if\spsformat1 \name{} \else

    \author{Daniel Romero, Pham Q. Viet, and Raju Shrestha\thanks{ The
            authors are with the Dept. of Information and Communication
            Technology, University of Agder, Jon Lilletunsvei 9, 4879 Grimstad,
            Norway. Email \{daniel.romero,viet.q.pham, raju.shrestha\}@uia.no.

            This
            work was supported by the Research Council of Norway through the
            IKTPLUSS Grant 311994.

            The present paper extends its conference precursor \cite{romero2022aerial} to
            accommodate constraints in the backhaul and more general propagation radio maps.
            It also includes a much more comprehensive simulation study where, among others,
            simulations using ray-tracing software are carried out.  }}

\fi

\maketitle

\newcommand{\bound}{\hc{b}}
\newcommand{\admmstop}{\hc{\epsilon}}

\newcommand{\dc}[1]{#1}

\newcommand{\includefig}[1]{\includegraphics[trim=25 0 50 30,clip,width=0.49\textwidth]{#1}}

\newcommand{\journal}[1]{\textcolor{purple}{#1}}
\newcommand{\journalc}[1]{#1} 
\newcommand{\conference}[1]{\textcolor{gray}{#1}}

\newcommand{\indicator}{\mathbb{I}}
\newcommand{\indicatorinf}{\mathcal{I}}

\newcommand{\groundregion}{{\hc{\mathcal{X}}}}
\newcommand{\height}{{\hc{h}}}
\newcommand{\flyregion}{{\hc{\mathcal{F}}}}

\newcommand{\usernum}{{\hc{M}}}
\newcommand{\userind}{{\hc{m}}}
\newcommand{\absnum}{{\hc{N}}}
\newcommand{\absind}{{\hc{n}}}

\newcommand{\slfgrid}{{\hc{\bar{\mathcal{X}}}}}
\newcommand{\slfgridptind}{{\hc{q}}}
\newcommand{\slfgridptnum}{{\hc{Q}}}
\newcommand{\slfgridside}{{\hc{\slfgridptnum}}_{0}}
\newcommand{\slfgridsidex}{{\hc{\slfgridptnum}}_x}
\newcommand{\slfgridsidey}{{\hc{\slfgridptnum}}_y}
\newcommand{\slfgridsidez}{{\hc{\slfgridptnum}}_z}
\newcommand{\slfgridpt}{{\hc{\bm x}}^\slfgrid}

\newcommand{\grid}{{\hc{\bar{\mathcal{F}}}}}
\newcommand{\gridptind}{{\hc{g}}}
\newcommand{\gridptnum}{{\hc{G}}}
\newcommand{\gridside}{{\hc{\gridptnum}}_{0}}
\newcommand{\gridsidex}{{\hc{\gridptnum}}_x}
\newcommand{\gridsidey}{{\hc{\gridptnum}}_y}
\newcommand{\gridsidez}{{\hc{\gridptnum}}_z}
\newcommand{\gridpt}{{\hc{\bm x}}^\grid}

\newcommand{\locs}{{\hc{x}}} 
\newcommand{\loc}{{\hc{\bm \locs}}}
\newcommand{\userloc}{{\hc{\bm x}}^\text{GT}}
\newcommand{\absloc}{{\hc{\bm x}}^\text{ABS}}
\newcommand{\wavelen}{\hc{\lambda}}
\newcommand{\gain}{\hc{\gamma}}
\newcommand{\shad}{\hc{\xi}}
\newcommand{\slf}{\hc{l}}
\newcommand{\weightfun}{\hc{w}}
\newcommand{\wfwidth}{\hc{\tilde\lambda}} 

\newcommand{\bandwidth}{\hc{W}}
\newcommand{\txpow}{\hc{P_\text{TX}}}
\newcommand{\noisepow}{\hc{\sigma^2_\text{N+I}}}

\newcommand{\capmat}{\hc{\bm C}}
\newcommand{\capsubmat}{\hc{\bbm C}}
\newcommand{\capfun}{\hc{C}}
\newcommand{\ratefun}{\hc{R}}
\newcommand{\maxratefun}{\hc{C}^\text{BH}}
\newcommand{\maxrate}[1]{\hc{c}^\text{BH}_{#1}}
\newcommand{\maxratevec}{\hc{\bm c}^\text{BH}}
\newcommand{\maxratesubvec}{\hc{\bbm c}^\text{BH}}

\newcommand{\caps}{\hc{c}}
\newcommand{\caprs}{\hc{\bar c}}
\newcommand{\capvec}{\hc{\bm \caps}}
\newcommand{\caprvec}{\hc{\bbm \caps}} 
\newcommand{\ratemat}{\hc{\bm R}}
\newcommand{\ratesubmat}{\hc{\bbm R}}
\newcommand{\rate}{\hc{ r}}
\newcommand{\rateuserabsinds}{\hc{ r}\entnot{\absind,\userind}}
\newcommand{\rateusergridptinds}{\hc{ r}_{\gridptind}\entnot{\userind}}
\newcommand{\raters}{\hc{ \bar r}} 
\newcommand{\ratevec}{\hc{\bm \rate}}
\newcommand{\ratervec}{\hc{\bbm \rate}} 
\newcommand{\raterowvec}{\hc{\bbm \rate}}
\newcommand{\minrate}{\hc{\ratefun}^\text{min}}
\newcommand{\acti}{\hc{\alpha}}
\newcommand{\actiold}{\hc{\tilde\alpha}}
\newcommand{\activec}{\hc{\bm \alpha}}
\newcommand{\invacti}{\hc{\beta}}

\newcommand{\sparsweight}{\hc{ w}}
\newcommand{\sparsweightvec}{\hc{\bm \sparsweight}}

\newcommand{\normvec}{\hc{\bm v}}
\newcommand{\norm}{\hc{ v}}
\newcommand{\vecent}[1]{[#1]}

\newcommand{\slack}{\hc{s}}
\newcommand{\slackvec}{\hc{\bm \slack}}
\newcommand{\slacklow}{\check \slack}
\newcommand{\slackhigh}{\hat \slack}

\newcommand{\dslackovec}{\hc{\bm \delta}_1}
\newcommand{\dslacktvec}{\hc{\bm \delta}_2}
\newcommand{\dslacktmat}{\hc{\bm \Delta}_2}
\newcommand{\dslackthvec}{\hc{\bm \delta}_3}
\newcommand{\dslackthmat}{\hc{\bm \Delta}_3}

\newcommand{\supervarvec}{\hc{\tbm x}}
\newcommand{\superweightvec}{\hc{\tbm w}}
\newcommand{\supervec}{\hc{\tbm b}}
\newcommand{\supermat}{\hc{\tbm A}}

\newcommand{\Xmat}{{\hc{\bm X}}}
\newcommand{\Ymat}{{\hc{\bm Y}}}
\newcommand{\Zmat}{{\hc{\bm Z}}}
\newcommand{\zrvec}{{\hc{\bbm z}}}
\newcommand{\zs}{{\hc{ z}}} 
\newcommand{\zrs}{{\hc{ \bar z}}} 
\newcommand{\zvec}{{\hc{\bm \zs}}}
\newcommand{\yvec}{{\hc{\bm y}}}
\newcommand{\ys}{{\hc{ y}}}
\newcommand{\tvec}{{\hc{\bm t}}}

\newcommand{\dvec}{{\hc{\bm d}}}

\newcommand{\Amat}{{\hc{\bm A}}}
\newcommand{\Bmat}{{\hc{\bm B}}}

\newcommand{\Emat}{{\hc{\bm E}}}
\newcommand{\Umat}{{\hc{\bm U}}}
\newcommand{\urvec}{{\hc{\bbm u}}}
\newcommand{\us}{{\hc{ u}}}
\newcommand{\urs}{{\hc{ \bar u}}}
\newcommand{\uvec}{{\hc{\bm \us}}}
\newcommand{\Ffun}{{\hc{F}}} 
\newcommand{\Gfun}{{\hc{G}}} 
\newcommand{\admmfunx}{\hc{f}}
\newcommand{\admmfunz}{\hc{h}}

\newcommand{\itnot}[1]{^{#1}}
\newcommand{\itind}{\hc{k}}
\newcommand{\entnot}[1]{[#1]} 

\newcommand{\gii}{_\gridptind\itnot{\itind}}
\newcommand{\admmstep}{\hc{\rho}}

\newcommand{\lambdas}{{\hc{\lambda}}}
\newcommand{\lambdaslow}{{\hc{\check\lambda}}}
\newcommand{\lambdashigh}{{\hc{\hat\lambda}}}
\newcommand{\nus}{{\hc{\nu}}}
\newcommand{\nuvec}{{\hc{\bm \nus}}}
\newcommand{\mus}{{\hc{ \mu}}}
\newcommand{\muvec}{{\hc{\bm \mus}}}
\newcommand{\mumult}{{\hc{\mu}}} 

\newcommand{\incs}{\hc{b}_\text{inc}}
\newcommand{\incvec}{\hc{\bm b}_\text{inc}}

\newcommand{\vis}{\hc{i}} 
\newcommand{\vivec}{\hc{\bm i}} 
\newcommand{\vicurrents}{\hc{i}_\text{current}}
\newcommand{\vicurrentvec}{\hc{\bm i}_\text{current}}
\newcommand{\viset}{\hc{\mathcal{I}}} 

\newcommand{\slftens}{\hc{L}}
\newcommand{\slften}{\hc{\bm \slftens}}

\newcommand{\deltaloc}{\hc{\bm \Delta}_{\loc}}
\newcommand{\deltagridvec}{\hc{\bm \delta}_{\slfgrid}}
\newcommand{\deltagrids}{\hc{ \delta}_{\slfgrid}}
\newcommand{\tcands}{\hc{ t}_\text{cand}}
\newcommand{\tcandvec}{\hc{\bm t}_\text{cand}}
\newcommand{\inext}{\hc{i}_\text{next}}
\newcommand{\tnext}{\hc{t}_\text{next}}
\newcommand{\integral}{\hc{{I}}}

\newcommand{\supmatsec}[1]{Sec.~\ref{#1} of the supplementary material}

\newcommand{\change}[1]{\textcolor{black}{#1}}

\begin{abstract}
    \begin{changes}
        The deployment of aerial base stations (ABSs) on unmanned aerial vehicles (UAVs) presents a promising solution for extending cellular connectivity to areas where terrestrial infrastructure is  overloaded, damaged, or absent. A pivotal challenge in this domain is to decide the locations of a set of ABSs to effectively serve ground-based users. Most existing approaches  oversimplify this problem by assuming that the channel gain between two points is  a function of solely distance and, sometimes, also the elevation angle. In turn, this paper leverages propagation radio maps to account for arbitrary air-to-ground channel gains. This methodology enables the identification of an approximately minimal set of locations where ABSs need to be deployed to  ensure that all ground terminals achieve a target service rate, while adhering to backhaul capacity limitations and avoiding designated no-fly zones. Relying on a convex relaxation technique and the alternating direction method of multipliers (ADMM), this paper puts forth a solver whose computational complexity scales linearly with the number of ground terminals. Convergence is established analytically and an extensive set of simulations corroborate  the merits of the proposed scheme relative to conventional methods.
    \end{changes}

\end{abstract}


\section{Introduction}
\label{sec:intro}

\begin{bullets}
    \blt[overview]
    \begin{bullets}
        \blt[motivation]Aerial base stations (ABSs), namely unmanned
        aerial vehicles (UAVs) equipped with on-board base stations, were
        conceived as a means to deliver cellular connectivity in areas
        where the terrestrial infrastructure is
        absent, overloaded, or damaged~\cite{zeng2019accessing}. This may occur for example in
        remote areas, in the vicinity of a crowded event, or after a
        natural disaster, such as a wildfire or a flood.
        \blt[ABSs] Users on the ground, here referred to as \emph{ground
            terminals} (GTs), are served by ABSs, which in turn connect to the
        terrestrial infrastructure through \emph{backhaul} links, possibly
        in multiple hops through other UAVs that act as relays.
        \blt[placement] Deploying ABSs involves addressing the problem of
        ABS placement, where one is given the locations of the GTs and
        must decide on a suitable set of spatial positions for the ABSs to
        effectively serve the GTs~\cite{viet2022introduction}.  This task
        is typically hindered by several challenges, remarkably:
        \begin{changes}
            \begin{itemize}
                \item[C1]                     the
                    uncertainty in the gain of the propagation channel between the GTs
                    and the potential ABS locations,
                \item[C2] the limited capacity of the
                    \emph{backhaul} links, and
                \item[C3] constraints on the positions that
                    the ABSs may adopt, often due to the presence of obstacles such as buildings or no-fly zones such as airports, natural parks,
                    embassies, or prisons.
            \end{itemize}
        \end{changes}

    \end{bullets}

    \blt[literature]
    \begin{bullets}
        \blt\cmt{one UAV} The problem of placing a single ABS has been
        extensively investigated in the literature; see
        e.g. \cite{han2009manet\journalc{,lee2011climbing},boryaliniz2016placement,chen2017map,wang2018adaptive}.
        \blt[multiple UAVs]However, since a single  ABS  may not suffice to cater for the needs of many situations,
        \begin{bullets}
            \blt[2D]a large number of works, including the
            present one, focus on placing multiple ABSs. A
            usual approach is to regard the height of the ABSs as given and
            address the problem of \emph{2D placement}, where the ABSs must be
            placed on a horizontal plane; see
            e.g.~\cite{lee2010decentralized,andryeyev2016selforganized,galkin2016deployment,lyu2017mounted,romero2019noncooperative,huang2020sparse,mach2021realistic,yin2021multiagent}.
            \blt[3D] Nevertheless, since the heights of the ABSs are useful
            degrees of freedom to optimize the target communication metric,
            the focus here is on \emph{3D placement}.

            Existing algorithms for  3D placement of multiple ABSs can be
            classified according to how they handle uncertainty in the
            air-to-ground channel.
            \begin{changes}
                \begin{enumerate}%
                    \item\cmt{channel agnostic} The first category comprises schemes
                    that do not explicitly model the channel. For example, in
                    \cite{park2018formation}, each GT associates with the ABS
                    from which it receives the strongest beacons, but nothing is
                    known or assumed about the channel gain from the GT
                    locations to a given location until an ABS is physically
                    there. This means that every iteration of the placement
                    algorithm requires arranging the ABSs at a particular set of
                    locations, which drastically limits convergence speed.
                    \item\cmt{free-space}The second class includes works that assume
                    free-space propagation and, therefore, the coverage area of
                    each ABS is a circle; see
                    e.g.~\cite{kim2018topology}. Unfortunately, this
                    assumption is too inaccurate in practice.
                    \item\cmt{empirical models\ra based on
                        \cite{alhourani2014urban,alhourani2014lap}}The third category
                    is made up of works that rely on the empirical model
                    from~\cite{alhourani2014urban\journalc{,alhourani2014lap}};
                    \begin{bullets}%
                        \blt[papers]see e.g.
                        \cite{kalantari2016number,hammouti2019mechanism,perabathini2019qos,liu2019deployment,shehzad2021backhaul}. 
                        \blt[limitations] These works use the mean provided by such
                        a model as the predictor of the channel gain, which is
                        equivalent to assuming that the gain of a link depends only
                        on its length and elevation. Again, the mismatch between this assumption and reality is significant  since two links with the same length and elevation
                        may exhibit highly different gains depending on whether
                        there are obstructions such as buildings between the
                        transmitter and the receiver.
                    \end{bullets}%
                    \item\cmt{Terrain map/3D}The fourth category, which can be termed
                    \emph{channel-aware}, is composed of works that rely on gain
                    predictions that do depend on the locations of the endpoints
                    of the link. To the best of our knowledge, only
                    \cite{qiu2020reinforcement,sabzehali2021orientation} fall in
                    this category. Unfortunately, these schemes incur prohibitive
                    complexity, assume an unlimited backhaul connection between
                    the ABSs and the terrestrial infrastructure, and cannot
                    guarantee a minimum service to the GTs.
                \end{enumerate}%
            \end{changes}
        \end{bullets}%
    \end{bullets}%

    \blt[contributions]\begin{changes} The main contribution  of this paper is to address the aforementioned limitations by proposing a scheme that relies on
        \emph{radio propagation maps}~\cite{romero2022cartography} to solve
        the problem of \emph{channel-aware 3D placement of multiple ABSs}. \end{changes}
    \begin{bullets}%
        \blt[channel map] Recall that a propagation map is a special kind
        of radio
        map~\cite{alayafeki2008cartography,romero2022cartography,shrestha2022surveying,teganya2019locationfree,teganya2020rme,zeng2021toward}
        that provides a channel metric of interest for every pair of
        transmitter and receiver locations. In this paper, this metric is
        the channel gain, which can be used to predict the capacity of the
        communication link between each candidate ABS location and every
        GT without deploying an ABS at that location to measure the
        channel. Two classes of propagation maps with complementary
        strengths will be considered, namely those obtained via
        ray-tracing and those that rely on the \emph{radio tomographic
            model}~\cite{patwari2008nesh,patwari2008correlated}.
        \begin{bullets}%
            \blt[ray tracing]The former are more suitable to frequency bands where the
            dominating propagation phenomena beyond free-space loss are
            reflection and diffraction.
            \blt[tomography]The latter are suitable to bands where
            the dominating propagation phenomenon is the absorption
            introduced by obstacles such as buildings. In this context, a
            secondary contribution of this paper is to adapt existing
            radio tomographic techniques to air-to-ground channels.
        \end{bullets}%

        \blt[Algorithm]
        \begin{bullets}
            \blt[novelty]To the best of our knowledge, the proposed scheme
            is the first for \emph{channel-aware ABS placement} that can
            guarantee a minimum rate for all GTs. Formally, the algorithm can
            find a feasible placement if it exists, where a feasible placement
            is an assignment of ABSs to spatial locations that ensures a
            target rate for all GTs according to the given propagation
            map.
            \blt[backhaul] Besides, unlike the vast majority of works in the literature,
            constraints in the backhaul link between the ABSs and the
            terrestrial infrastructure as well as
            \blt[no-fly zones]no-fly zones can be enforced.
            \blt[criterion]Subject to these constraints, the proposed algorithm
            approximately minimizes the number of ABSs required to serve all
            GTs. Note that this is of special interest in emergency scenarios,
            which is one of the main use cases of ABSs.
            \blt[Approach] The algorithm relies on a sparse optimization
            formulation that naturally arises from a discretization of the
            space of candidate ABS positions, as required to be able to
            utilize propagation maps in a tractable fashion. To counteract the
            high-dimensionality of the 3D placement problem, a
            linear-complexity and highly parallelizable algorithm is developed
            based on the alternating-direction method of multipliers
            (ADMM)~\cite{boyd2011distributed}.

        \end{bullets}
    \end{bullets}%

    \blt[simulations]Experiments with tomographic and ray-tracing models
    showcase a great reduction in the number of required ABSs as
    compared to existing algorithms. To complement this manuscript, an
    open-source simulator was released to allow developing and testing
    algorithms for ABS placement. This simulator and the code needed to
    reproduce all experiments is available at
    \url{https://github.com/uiano/ABS_placement_via_propagation_maps}.

    \begin{figure}[!t]
        \centering
        \includegraphics[width=0.48\textwidth]{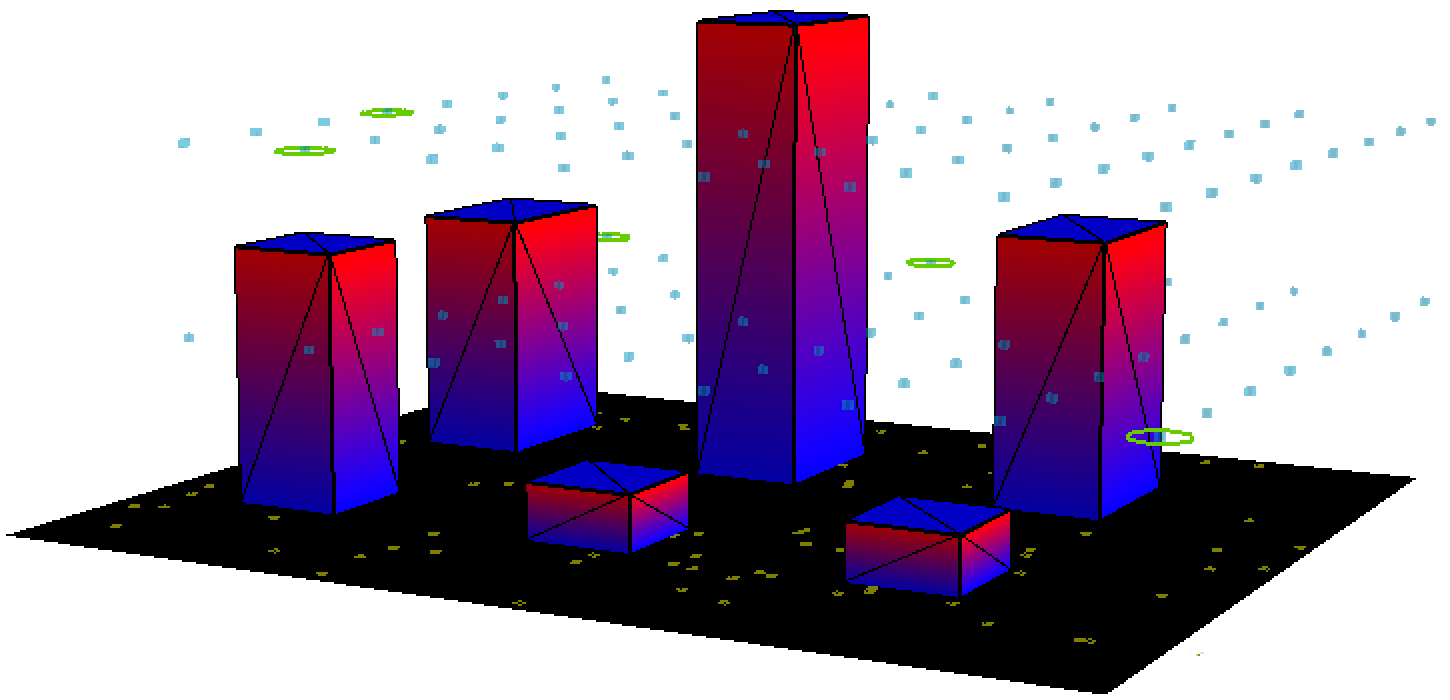}
        \caption{Example of ABS placement in an urban environment. GTs are
            represented by markers on the ground, flight grid points by blue dots,
            and ABS positions by green circles.}
        \label{fig:urban}
    \end{figure}

    \blt[paper structure]\emph{Paper structure.}  The rest of the paper
    is organized as follows. Sec.~\ref{sec:model} presents the model and
    formulates the problem. Two approaches for predicting the capacity
    of a link between arbitrary pairs of endpoints of the air-to-ground
    channel are then described in Sec.~\ref{sec:radiomaps}. The problem
    of ABS placement and rate allocation is then addressed in
    Sec.~\ref{sec:placement} and a solver with linear complexity is
    developed in Sec.~\ref{sec:solver}. Two extensions are presented in Sec.~\ref{sec:extensions}. Finally,
    Secs.~\ref{sec:experiments} and~\ref{sec:conclusions} respectively
    present numerical results and conclusions. \change{The supplementary material~\cite{romero2023placementarxiv} contains extended derivations, proofs, and additional simulations.}


    \blt[notation]\emph{Notation.}
    \begin{bullets}%
        \blt[] $\rfield_+$ is the set of non-negative real numbers and
        $\rfield_{++}$ is the set of positive real numbers.
        \blt[] Boldface uppercase (lowercase) letters denote matrices
        (column vectors).
        \blt[] $a\entnot{i}$ represents the $i$-th entry of vector $\bm a$.
        \blt[] Notation $\bm 0$ (respectively $\bm 1$) refers to the
        matrix of  the appropriate dimensions with all zeros (ones).
        \blt[] $\|\bm A\|_\frob$ denotes the Frobenius norm of matrix $\bm A$,
        whereas $\|\bm a\|_p$ denotes the $\ell_p$-norm of vector $\bm
            a$. With no subscript, $\|\bm a\|$ stands for the
        $\ell_2$-norm.
        \blt Inequalities between vectors or matrices must be
        understood entrywise.
        \blt The Kronecker product is denoted by $\otimes$.
        \blt If $\bm a$ and $\bm b$ are vectors of the same dimension,
        then $\bm a\odot \bm b$ is the entrywise product of $\bm a$ and
        $\bm b$, whereas $\bm a\div \bm b$ is the entrywise quotient of
        $\bm a$ and $\bm b$.

    \end{bullets}%
\end{bullets}%

\section{Model and Problem Formulation}
\label{sec:model}

\begin{bullets}%
    \blt[model]
    \begin{bullets}%
        \blt[users]Consider $\usernum$ GTs  at positions
        $\{\userloc_1, \ldots,\userloc_\usernum\} \subset \groundregion
            \subset\rfield^3$, where the region $\groundregion$ represents an
        arbitrary set of spatial locations, including for example points
        on the street, inside buildings, inside vehicles, and so on.
        \blt[also air]The proposed scheme carries over unaltered to the
        scenario where $\groundregion$ includes points in the airspace and
        some or all users are airborne, which can be of interest e.g. to
        deploy auxiliary ABSs as picocells. However, to simplify the
        exposition, this possibility is neglected and the users will be
        referred to throughout as GTs.

        \blt[abss]To provide connectivity to these GTs, $\absnum$ ABSs
        are deployed at locations
        $\{\absloc_1,\ldots,\absloc_\absnum \}$ $\subset\flyregion\subset
            \rfield^3$, where  $\flyregion$ comprises all
        spatial positions where a UAV is allowed to fly. This excludes no-fly
        zones, airspace occupied by buildings, and altitudes out
        of legal limits.

        \blt[channel]%
        \begin{bullets}
            \blt[abs-gt]
            \begin{bullets}
                \blt[uplink]For the sake of specificity, it will be assumed that
                data packets originated in a remote location are sent from the
                terrestrial infrastructure to the ABSs through a \textit{backhaul link}
                and the ABSs forward these packets to their intended GTs
                through the \emph{downlink} of the \emph{radio access
                    network}. However, the entire discussion applies also to the
                \emph{uplink}, i.e., when the data packets are originated at
                the GTs and sent to the terrestrial
                infrastructure through the ABSs.

                \blt[interference]\change{Following standard practice in the
                    literature (see
                    e.g.~\cite{sabzehali2021orientation,huang2020sparse,lyu2017mounted,kalantari2016number,galkin2016deployment,mach2021realistic,liu2019deployment,kim2018topology,lee2010decentralized,park2018formation,qiu2020reinforcement}),
                    the ABSs utilize orthogonal channel resources and therefore
                    do not introduce interference to GTs that are not associated
                    with them.  This is possible in most of the application
                    scenarios of interest, where operations take place in areas
                    without terrestrial base stations or where the existing ones
                    are damaged. Consequently, a large portion of spectral resources are available and, thus, the ABSs can use disjoint bands or
                    quasi-orthogonal spreading codes.}
                \blt[capacity]Disregarding also frequency-selective
                effects for simplicity, the capacity of the communication
                link between an ABS at position $\absloc\in\groundregion$
                and the $\userind$-th GT is given by
                \begin{align}
                    \label{eq:capfun}
                    \capfun_\userind(\absloc) = \bandwidth \log_{2}\left(
                    1 + \frac{\txpow 10^{\gain_\userind( \absloc)/10}}{\noisepow}
                    \right) ,
                \end{align}
                where $ \bandwidth$ denotes bandwidth, $\txpow$ represents
                the transmit power, $\noisepow$ \change{aggregates the
                    noise and interference power}, and $\gain_\userind(
                    \absloc)$ is the channel gain, which is described in
                Sec.~\ref{sec:radiomaps}.

            \end{bullets}
            \blt[backhaul]Unlike most schemes in the literature of ABS
            placement, the present work can accommodate constraints in the
            backhaul.
            \begin{bullets}
                \blt[max capacity]To formalize such constraints, let
                $\maxratefun(\absloc)$ denote the maximum rate of the link
                between the ABS at $\absloc$ and the terrestrial
                ground station(s) that serve(s) it. An equation like
                \eqref{eq:capfun} can also be established to express
                $\maxratefun(\absloc)$ in terms of the gain of the
                relevant channel(s). Note that, as the notation suggests,
                $\maxratefun(\absloc)$  generally depends on the ABS
                position $\absloc$ -- typically, the greater the distance
                from $\absloc$ to the terrestrial ground stations, the
                lower $\maxratefun(\absloc)$.
                \blt[bw allocation]With $ \ratefun_\userind(\absloc)$
                denoting the downlink rate that an ABS at $\absloc$
                allocates to the $\userind$-th GT, the backhaul rate
                constraint imposes that
                $\sum_\userind
                    \ratefun_\userind(\absloc)\leq \maxratefun(\absloc)$.

            \end{bullets}

            %
            %

        \end{bullets}%
    \end{bullets}%

    \blt[problem formulation]
    \begin{bullets}%
        \blt[criterion]The problem is to find a minimal number
        of ABS locations that guarantee that every user
        receives a rate of at least $\minrate$. This criterion
        arises naturally in some of the main use cases of
        UAV-assisted networks such as emergency response or
        disaster management. Motivated by this scenario and to
        enhance flexibility in the deployment, each GT may be
        served by multiple ABSs. This means that the rate that
        the $\userind$-th GT receives is
        $\sum_\absind\ratefun_\userind(\absloc_\absind)$, where
        $\absloc_\absind$ denotes the location of the
        $\absind$-th ABS. \change{The easiest way of
            implementing the communication between one GT and multiple ABSs is by time, code, or
            frequency multiplexing; see also Sec.~\ref{sec:minnumconnections}.}

        \blt[problem]To summarize,
        the problem can be formulated as follows:
        \begin{subequations}
            \label{eq:problemf}
            \begin{align}
                 & \hspace{-2cm} \minimize_{\absnum,
                \{\absloc_\absind\}_{\absind=1}^\absnum, \{\rateuserabsinds\}_{\userind=1,\absind=1}^{\usernum,\absnum}}~ \absnum \\
                \st~               \label{eq:problemfbkhaul}
                 & \sum_\userind \rateuserabsinds \leq
                \maxratefun(\absloc_\absind),
                \\
                 & \sum_\absind \rateuserabsinds
                \geq \minrate,                                                                                                    \\
                 & 0\leq \rateuserabsinds \leq \capfun_\userind(\absloc_\absind),
                \label{eq:problemfcap}
                \\
                \label{eq:problemfpos}
                 & \absloc_\absind \in \flyregion,
            \end{align}
        \end{subequations}
        where the constraints need to hold for all $\userind$ and
        $\absind$ and the earlier notation
        $\ratefun_\userind(\absloc_\absind)$ has been replaced with
        $\rateuserabsinds \define
            \ratefun_\userind(\absloc_\absind)$ to emphasize that it
        refers to optimization variables, not to functions.
        Observe that Problem \eqref{eq:problemf} constitutes a
        joint placement and rate-allocation problem. Note also that
        the same minimum rate $\minrate$ is imposed for all GTs,
        but different rates can be set up to straightforward
        modifications.

        \blt[challenges] From a practical perspective, solving
        \eqref{eq:problemf} involves two challenges. First,
        $\capfun_\userind(\absloc_\absind)$ and
        $\maxratefun(\absloc_\absind)$ depend on the channel gain
        of the corresponding downlink and backhaul links, which is
        generally unknown. This issue will be addressed in
        Sec.~\ref{sec:radiomaps}. Second, given
        $\capfun_\userind(\absloc_\absind)$ and
        $\maxratefun(\absloc_\absind)$, one needs to find the
        positions of the ABSs and the rate allocations that solve
        \eqref{eq:problemf}. This will be the subject of
        Sec.~\ref{sec:placement}.

    \end{bullets}%

\end{bullets}%

\section{Capacity Prediction via Propagation Radio Maps }
\label{sec:radiomaps}

\begin{figure}[!t]
    \centering
    \includegraphics[width=.49\textwidth]{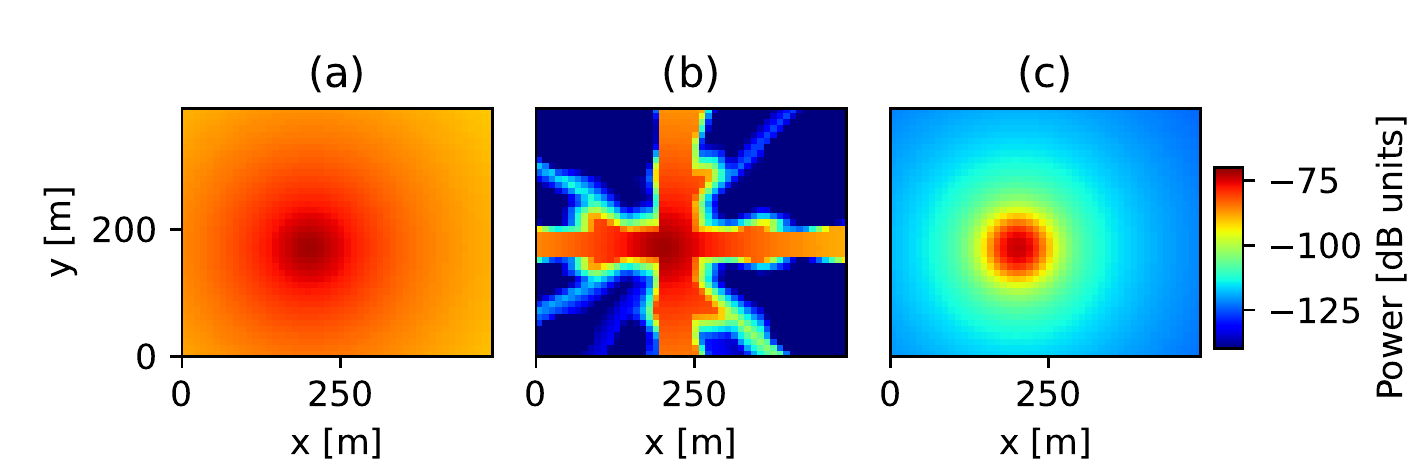}
    \caption{Radio maps that show $\gain_\userind(\absloc)$ vs. $\absloc\define(x,y,z)$ on a horizontal plane of height $z=40$ m. A GT is at a street location $\userloc_\userind$ in  a dense urban environment. Three channel models are used: (a) the free-space channel, (b) the tomographic channel, and (c) the model in~\cite{alhourani2014lap}. }
    \label{fig:maps}
\end{figure}

\begin{bullets}
    \blt[overview]
    \begin{bullets}
        \blt[Motivation]As described in Sec.~\ref{sec:intro}, nearly all\footnote{
            \change{In \cite{qiu2020reinforcement} and \cite{chen2021relay}, other channel models are used. However, they are restricted to the case of a single ABS and single GT.
            }
        }
        existing schemes  for ABS placement
        assume that $\capfun_\userind(\absloc_\absind)$ depends on
        $\userloc_\userind$ and $\absloc_\absind$ only through the length and
        (possibly) the elevation angle of the line segment between these two
        points. \change{More specifically, this applies e.g. to the free-space
            model used e.g. in
            \cite{romero2019noncooperative,huang2020sparse,galkin2016deployment,yin2021multiagent,kim2018topology,park2018formation},
            as well as the model~\cite{alhourani2014urban,alhourani2014lap} used e.g.
            in~ \cite{ mozaffari2016coverage, kalantari2016number,
                perabathini2019qos, shehzad2021backhaul, liu2019deployment,
                boryaliniz2016placement}.  Fig.~\ref{fig:maps} visually shows that
            the channel predicted by these two models is  radially symmetric on
            every horizontal plane.}
        This simplification may entail a
        significant performance degradation since channel gain in
        reality is heavily affected by the environment: two links
        with the same length and elevation may experience very
        different channel gain depending on the position, shape,
        and material of the surrounding obstacles and
        scatterers. The same observation applies to
        $\maxratefun(\absloc_\absind)$.

        \blt[Here]\begin{changes}~A
            propagation radio map is a  radio map that
            provides a certain channel metric
            between two given spatial
            locations. In this paper, the  metric is the channel gain, in which case a propagation radio map is often referred to as a \emph{channel-gain map}~\cite{patwari2008nesh,agrawal2009correlated,patwari2008correlated,romero2022cartography,
                romero2018blind,kim2011cooperative,wilson2009regularization,hamilton2014modeling}.\end{changes}
        This section describes how  propagation radio maps can be used to obtain
        $\capfun_\userind(\absloc_\absind)$ using two classes of existing channel models with complementary strengths.
        \blt[only downlink]\change{The same considerations apply to
            $\maxratefun(\absloc_\absind)$.}
    \end{bullets}

\end{bullets}

\subsection{Ray-tracing Models}

\begin{bullets}
    \blt[overview]Ray-tracing techniques~\cite{imai2017survey} can
    be used to predict $\gain_\userind( \absloc_\absind)$ for
    arbitrary pairs $(\userloc_\userind,\absloc_\absind)$ given a
    3D model of the propagation environment, which includes
    buildings, terrain elements, and other objects. The idea is to
    geometrically obtain the relevant propagation paths between
    the transmitter and the receiver and subsequently calculate
    the attenuation and delay or phase change of each path.

    \blt[benefits]
    \begin{bullets}
        \blt Ray-tracing models are especially attuned to
        high-frequency bands, where the dominating propagation
        phenomena are reflection and diffraction.
    \end{bullets}%
    \blt[limitations] The limitation is that, 
    to track changes in the channel, one  needs to track
    changes in the 3D model, which may not be viable.
    \blt[ABS placement]Thus, in the case of ABS placement, it
    makes sense to use ray-tracing to precompute the values of
    $\gain_\userind(\absloc_\absind)$ (cf.
    Sec.~\ref{sec:ray_exp}) for the given scenario and ignore the
    effects of changes in the channel. If the impact of the
    changes (e.g. moving vehicles) is not too large, using a
    ray-tracing propagation map would still yield better placement
    solutions than the approaches
    mentioned at the beginning of Sec.~\ref{sec:radiomaps}.

\end{bullets}

\subsection{Radio Tomographic Models}
\label{sec:radiotomographicmaps}


\begin{bullets}%

    \blt[overview]As opposed to ray-tracing models, which account for
    reflection and diffraction effects, the radio tomographic
    model\begin{changes}~\cite{patwari2008nesh,agrawal2009correlated,patwari2008correlated,romero2022cartography,
            romero2018blind,kim2011cooperative,wilson2009regularization,hamilton2014modeling}\end{changes} focuses on the absorption undergone by radio
    waves on the straight path between the transmitter and the receiver.
    \blt[benefits] This model is best suited to traditional cellular
    communication frequencies, where radio waves  readily penetrate
    structures such as buildings.
    \blt[only ground so far]However, the existing works in this context
    focus on ground-to-ground channels. To the best of our knowledge,
    this is the first work to apply  radio tomography to
    air-to-ground channels. This entails special challenges due to the
    high dimensionality of the underlying space, which render
    existing techniques unsuitable. After describing the radio
    tomographic model, this section proposes an algorithm to bypass
    these difficulties.

    \blt[model description]The radio tomographic model dictates that the
    channel gain between $\userloc_\userind$ and $\absloc$ can be
    decomposed into a free-space loss component and a shadowing component
    as\footnote{Expression \eqref{eq:gain} assumes that small-scale fading
        has been averaged out. Otherwise, a random term can be added to
        the right-hand side.  }
    \begin{bullets}%
        \blt[channel gain]
        \begin{align}
            \label{eq:gain}
            \gain_\userind(\absloc) = 20\log_{10}\left( \frac{\wavelen}{4\pi
                \|\userloc_\userind- \absloc\|} \right) -
            \shad(\userloc_\userind, \absloc),
        \end{align}
        where $\wavelen$ is the wavelength associated with the carrier
        frequency  and the shadowing function $\shad$
        is given by~\cite{patwari2008nesh}
        \blt[Shadowing]
        \begin{bullets}%
            \blt[integral]
            \begin{align}
                \label{eq:tomoint}
                \shad(\loc_1,\loc_2) = \frac{1}{ \| \loc_1 -\loc_2\|_2^{1/2}
                } \int_{\loc_1}^{\loc_2}\slf(\loc)d\loc.
            \end{align}
            The non-negative function $\slf$ inside the line integral is
            termed \emph{spatial loss field} (SLF) and quantifies the
            local attenuation (absorption) that a signal suffers at each
            position.

        \end{bullets}%
    \end{bullets}%
    \blt[tasks]Two tasks are of interest: (T1) evaluate
    $\gain_\userind(\absloc)$ for a pair of locations
    $(\userloc_\userind, \absloc)$, which involves evaluating the integral in
    \eqref{eq:tomoint} and substituting the result into
    \eqref{eq:gain}; (T2) estimate
    $\slf$ given a set of measurements of the form $(\absloc,
        \userloc_\userind, \gain_\userind(\absloc))$ collected beforehand.

    \blt[evaluation]In both tasks, $\slf$ needs to be
    discretized, which can be accomplished by storing its
    values
    $\slf(\slfgridpt_1),\ldots,\slf(\slfgridpt_\slfgridptnum)$
    on a regular grid of $\slfgridptnum$ points
    $\slfgrid\define\{\slfgridpt_1,\ldots,\slfgridpt_\slfgridptnum\}$. The
    rest of this section explains how
    $\shad(\loc_1,\loc_2)$ can be expressed in terms of
    these values. This addresses (T1) and enables (T2) in
    combination with standard estimators; see
    e.g.~\cite{wilson2009regularization,kanso2009compressed,romero2018blind}

    \begin{bullets}%
        \blt[conventional]
        \begin{bullets}%
            \blt[description]%
            \begin{bullets}%
                \blt[innner prod]The conventional approach
                approximates the right-hand side of \eqref{eq:tomoint}
                as the weighted sum~\cite{hamilton2014modeling}
                \begin{align}
                    \label{eq:tomointapprox}
                    \shad(\loc_1,\loc_2)
                    \approx \sum_\slfgridptind \weightfun(\loc_1, \loc_2, \slfgridpt_\slfgridptind)\slf(\slfgridpt_\slfgridptind),
                \end{align}
                \blt[weight fun]where the weight function
                $\weightfun(\loc_1, \loc_2, \slfgridpt)$ aims at
                assigning a non-zero weight only to those grid
                points $\slfgridpt$ lying close to the line segment between
                $\loc_1$ and $\loc_2$. Although there are some
                variations, the functions $\weightfun$ adopted in
                the literature are non-zero only when $\slfgridpt$
                lies inside an ellipsoid with foci at $\loc_1$ and
                $\loc_2$~\cite{hamilton2014modeling,romero2018blind,gutierrezestevez2021hybrid}; 
                see the ellipses in Fig.~\ref{fig:tomography} for a
                depiction in 2D.
            \end{bullets}

            \begin{figure}[!t]
                \centering
                \includegraphics[width=.49\textwidth]{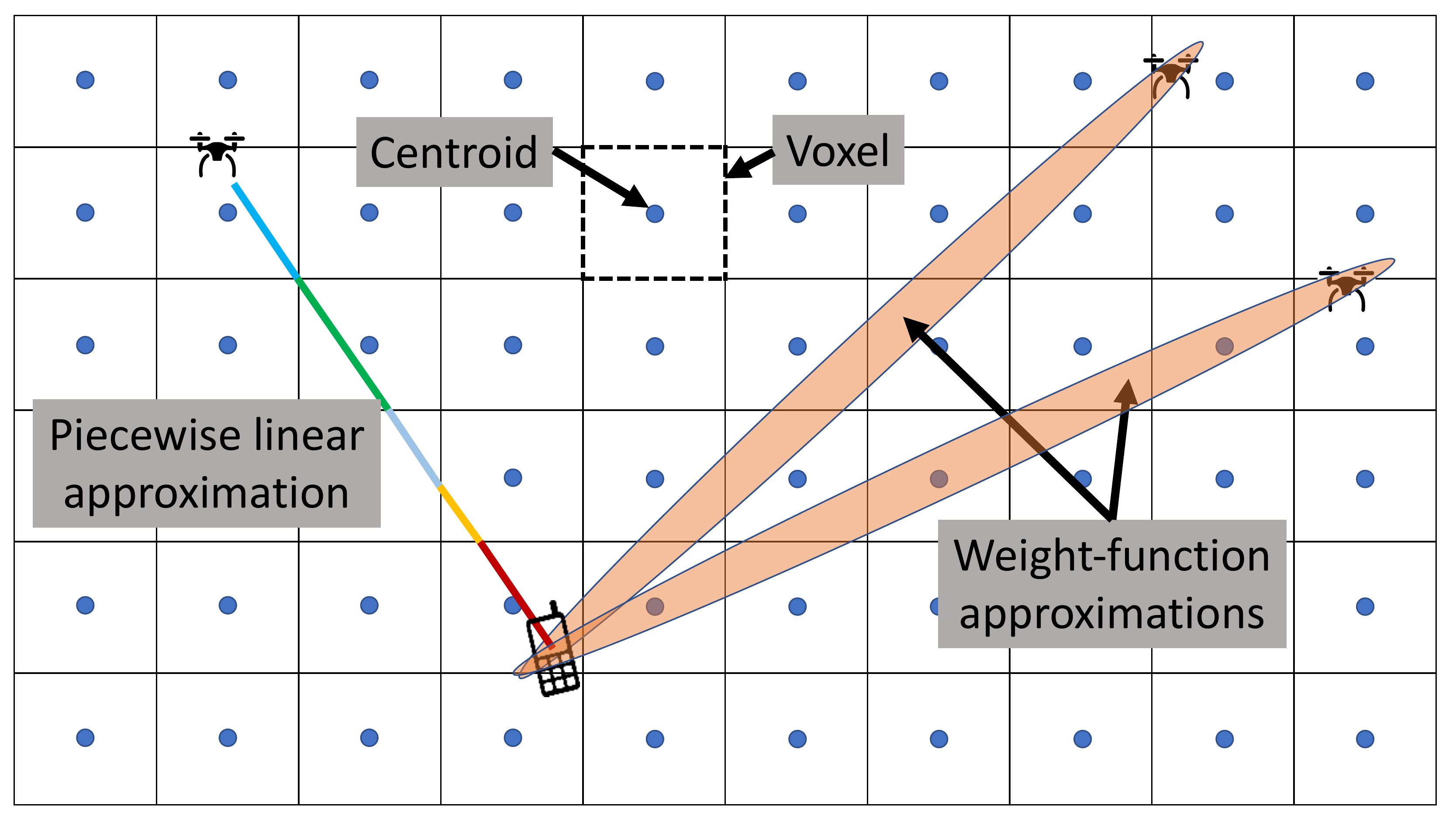}
                \caption{2D illustration of the conventional weight-function
                    approximation of the tomographic integral
                    \eqref{eq:tomoint} (orange ellipses) and the approximation
                    adopted here (colored line segment). Observe that the upper
                    ellipse contains no centroid and, therefore, the
                    approximation will yield zero attenuation regardless of the
                    values of the SLF. }
                \label{fig:tomography}
            \end{figure}

            \blt[Limitations]Such an approximation suffers from
            several limitations that render it impractical for the
            application at hand.
            \begin{bullets}%
                \blt[discontinuous]First, since existing choices of
                $\weightfun(\loc_1, \loc_2, \slfgridpt)$ are
                discontinuous functions, the resulting approximations
                of $\shad(\loc_1,\loc_2)$ are also discontinuous, which
                may lead to erratic behavior. For example, it may
                well happen that $\sum_\slfgridptind
                    \weightfun(\loc_1, \loc_2,
                    \slfgridpt_\slfgridptind)\slf(\slfgridpt_\slfgridptind)=0$
                even when $\loc_1 \neq \loc_2$ and
                $\slf(\slfgridpt_\slfgridptind)\neq 0~\forall
                    \slfgridptind$. This happens when $\loc_1$ and $
                    \loc_2$ are such that no grid point falls inside the
                elliptical support of $\weightfun(\loc_1, \loc_2,
                    \slfgridpt)$; see the upper ellipse in
                Fig.~\ref{fig:tomography}.
                \blt[large lambda] Since $\weightfun$ is typically
                chosen so that its support coincides with the first
                Fresnel ellipsoid, the length of its minor axis is
                roughly proportional to $\sqrt{\wavelen}$.  If the carrier frequency is low, the
                ellipsoid is large and, therefore, likely to contain
                a large number of grid points, which somehow
                mitigates the effects of the discontinuity.
                \blt[number of grid pts] Conversely, if the
                carrier frequency is high, one could think of
                creating a sufficiently dense grid so that the
                distance between grid points is small relative to
                $\wavelen$. However, it is easy to see that this
                may entail a prohibitively high $\slfgridptnum$. For example, if one wishes to set the
                grid points, say, 5 cm apart and the region of
                interest is 1 km $\times$ 1 km $\times$ 100 m,
                then the total number of grid points is
                $10^{11}$, which is obviously impractical.
                \blt[complexity] Another critical limitation is
                computational complexity. Observe that
                \eqref{eq:tomointapprox} generally requires
                evaluating $\weightfun(\loc_1, \loc_2,
                    \slfgridpt_\slfgridptind)$ and the product $
                    \weightfun(\loc_1, \loc_2,
                    \slfgridpt_\slfgridptind)\slf(\slfgridpt_\slfgridptind)$
                for each grid point. Thus, if the grid is
                $\slfgridside \times \slfgridside \times
                    \slfgridside$, the complexity of approximating
                $\shad(\loc_1,\loc_2)$ is
                $\mathcal{O}(\slfgridside^3)$. This is
                computationally problematic for the application at
                hand since $\shad(\loc_1,\loc_2)$ must be
                approximated at a large number of locations to solve
                \eqref{eq:problemf}.

            \end{bullets}%

        \end{bullets}%

        \blt[proposed]
        \begin{bullets}%
            \blt[overview] To remedy these issues, this paper
            advocates approximating the integral in
            \eqref{eq:tomoint} as a line integral of a piecewise
            constant approximation of $\slf$. 
            The resulting approximation is continuous, can be used
            with large grid point spacing, and can be computed with
            complexity only $\mathcal{O}(\slfgridside)$ for a
            $\slfgridside\times\slfgridside\times\slfgridside$
            grid. 
            %
            \blt[piecewise approx] This technique, commonly used in
            other disciplines (see references in
            \cite{mitchell1990comparison}), involves splitting the
            3D space into voxels centered at the grid points
            $\slfgrid\define\{\slfgridpt_1,\ldots,\slfgridpt_\slfgridptnum\}$
            and approximating $\slf$ by a function that takes the
            value $\slf(\slfgridpt_\slfgridptind)$ at all points of
            the $\slfgridptind$-th voxel. The resulting piecewise
            constant approximation of $\slf$ can be integrated by
            determining the positions of the crossings between the
            voxel boundaries and the line segment between $\loc_1$
            and $\loc_2$; see~Fig.~\ref{fig:tomography}.

            \begin{algorithm}[t!]
                \caption{Tomographic Integral Approximation}
                \label{algo:tomo}
                \begin{minipage}{200cm}
                    \begin{algorithmic}[1]
                        \STATE Input: $\loc_1$, $\loc_2$, grid spacing vector
                        $\deltagridvec\in\rfield^3$, \dc{\\\quad\quad~} SLF tensor
                        $\slften\in\rfield^{\slfgridsidex\times \slfgridsidey\times \slfgridsidez}$.
                        \STATE Initialize: \dc{
                        }
                        {\\\quad\quad~} $\deltaloc=\loc_2-\loc_1$,  $\incvec =
                            \sign(\deltaloc)$, $\integral=0$, $t=0$.
                        \STATE Set zero entries of $\deltaloc$ to 1 \# To avoid dividing by 0
                        \STATE Set $\vicurrentvec = \text{round}( \loc_1\div\deltagridvec)$ \#
                        Voxel indices
                        \WHILE{$t<1$}

                        \STATE Set $\tcandvec =( \deltagridvec \odot (\vicurrentvec + \incvec/2) -
                            \loc_1)\div \deltaloc$
                        \label{step:crossings}

                        \STATE Set $\inext = \argmin_i
                            \tcands\entnot{i}~\st~\incs\entnot{i}\neq 0$

                        \STATE Set $\tnext = \min(1,\tcands\entnot{\inext})$
                        \label{step:nextcrossing}

                        \STATE Set $\integral =\integral + (\tnext - t) \slftens\entnot{\vicurrentvec}$

                        \STATE Set $t = \tnext$

                        \STATE Set $\vicurrents\entnot{\inext} = \vicurrents\entnot{\inext} + \incs\entnot{\inext}$

                        \ENDWHILE
                        \STATE \textbf{return} $\|\loc_2-\loc_1\|^{1/2}\integral$
                    \end{algorithmic}
                \end{minipage}
            \end{algorithm}

            \blt[derivation] \change{Algorithm~\ref{algo:tomo} is our implementation of
                this approach. Details on its derivation can be found in
                \supmatsec{sec:algotomoderivation}.}
            \blt[algo benefits]This algorithm solves the
            limitations of the conventional approximation outlined earlier.
            \begin{bullets}%
                \blt[continuous]First, Algorithm~\ref{algo:tomo} yields an
                approximation of $\shad(\loc_1,\loc_2)$ that is a continuous
                function of $\loc_1$ and $\loc_2$ since the line integral of a
                piecewise constant function is a continuous function of the
                endpoints. Besides, the algorithm does not suffer from the
                issue of the approximation becoming zero when the elliptical
                support of the weight function in \eqref{eq:tomointapprox}
                misses all grid points.
                \blt[num. voxels]For this reason, the voxels can now be kept
                large regardless of the wavelength and, therefore, the total
                number of voxels can be kept low enough to be handled given the
                available computational resources.
                \blt[complexity]Finally, as indicated earlier in
                Sec.~\ref{sec:radiomaps}, the computational complexity of
                Algorithm~\ref{algo:tomo} is much smaller than the one of the
                conventional approximation. Specifically, one can observe in
                Algorithm~\ref{algo:tomo} that a constant number of products and
                additions is required for each crossing. The total number of
                crossings is at most
                $\slfgridsidez+\slfgridsidey+\slfgridsidez$, which means that,
                if $\slfgridsidex=\slfgridsidey=\slfgridsidez=\slfgridside$,
                then the total complexity of Algorithm~\ref{algo:tomo} is
                $\mathcal{O}(\slfgridside)$, whereas the complexity of the            standard approximation is, recall,~$\mathcal{O}(\slfgridside^3)$.

            \end{bullets}%

        \end{bullets}%

    \end{bullets}%
\end{bullets}%

\begin{changes}

    \begin{myremark}
        \label{remark:choice}
        The choice between a radio tomographic and a ray-tracing model depends on several aspects. First, tomographic models are more suitable to lower frequencies, where the dominating propagation phenomenon is absorption. In contrast, ray-tracing models are more suitable to higher frequencies, where reflection and diffraction have a greater impact than absorption.
        Second, the nature of the information that will be used to construct the radio map is also important. Tomographic models require radio measurements whereas  ray-tracing models require the geometry and electromagnetic properties of the objects in the environment. This is discussed further in Remark~\ref{remark:requiredinfo}.
        Third, ray-tracing models are more computationally intensive than tomographic models. More details are provided in Remark~\ref{remark:mobile}.
    \end{myremark}

    \begin{myremark}
        \label{remark:requiredinfo}
        Ray-tracing methods require the geometry and electromagnetic properties of the objects in the environment. Technologies that can be utilized to acquire this kind of information include satellite imaging, lidar on mobile robots (e.g. UAVs), simultaneous localization and mapping (SLAM)~\cite{thrun2008simultaneous}, and integrated sensing and communication (ISAC)~\cite{liu2022integrated} to name a few. However, needless to say, the accuracy of this information may be variable.
        In the case of tomographic maps, a set of measurements of the channel gain between pairs of locations must be collected. The precise number of such pairs that is required to attain a target estimation accuracy will generally depend on the complexity of the environment as well as the estimation method.
        In any case, deviations between the radio map estimates and the true radio maps are unavoidable. The impact of these deviations on the proposed scheme will be addressed in Sec.~\ref{sec:masnumservedusers} and at the end of Sec.~\ref{sec:tomographic_exp}.
    \end{myremark}

    \begin{myremark}
        \label{remark:mobile}
        Whether the GTs are mobile or not constitutes an important aspect when dealing with propagation radio maps in practice. With tomographic models, the spatial loss field is stored instead of the radio map itself. Whenever a prediction of the channel gain between two points is required, a line integral is computed using an algorithm such as Algorithm~\ref{algo:tomo}. The computational complexity of such an operation is low and can therefore be performed in real time. When using a ray-tracing map, a similar approach may be attempted: the 3D model of the environment can be stored and the ray-tracing algorithm can be executed every time a channel gain prediction is required. Unfortunately, the computational complexity of such an operation may be prohibitive for real-time implementation. If the GTs are static, one may  precompute the channel gain between the GT locations and
        a sufficiently large number of potential locations for the ABSs, e.g. the grid  $\grid$ in Sec.~\ref{sec:placement}. In contrast, if the GTs are mobile, precomputing the channel gain between a sufficiently large number of pairs of potential locations of the ABSs and GTs may not be feasible in terms of memory. For this reason, when the GTs are mobile users, it may be preferable to evaluate the radio map in real time by adjusting the parameters of the ray-tracing algorithm to reduce its computational complexity at the expense of decreasing accuracy. As discussed in Sec.~\ref{sec:masnumservedusers}, this may degrade the performance  of the placement algorithm to some extent.
    \end{myremark}

\end{changes}

\section{Discretization and Convexification}
\label{sec:placement}



\begin{bullets}%

    \blt[overview]

    \begin{bullets}%
        \blt[channel known]The approaches in Sec.~\ref{sec:radiomaps} make
        it possible to predict the channel gain
        $\gain_\userind(\absloc_\absind)$ for all required pairs of GT
        location $\loc_\userind$ and ABS location $\absloc_\absind$. These
        predictions can then be substituted into \eqref{eq:capfun} to
        obtain the capacity values $\capfun_\userind(\absloc_\absind)$ that
        appear in the constraint \eqref{eq:problemfcap}. A similar
        observation applies to $\maxratefun(\absloc_\absind)$ in
        \eqref{eq:problemfbkhaul}.
        \blt[nonconvex]Unfortunately, solving \eqref{eq:problemf} is
        challenging: even if $\absnum$ were known and one just needed to find
        feasible $\{\absloc_\absind\}_{\absind=1}^\absnum$, the problem would
        still be non-convex since the right-hand sides of
        \eqref{eq:problemfbkhaul} and \eqref{eq:problemfcap} are, in general,  non-concave
        functions of $\absloc_\absind$.
        \blt[discretization]To bypass this difficulty, the flight
        region $\flyregion$ will be discretized into
        a \emph{flight grid} $\grid\define
            \{\gridpt_1,\ldots,\gridpt_\gridptnum \}\subset\flyregion\subset
            \rfield^3$; see~Fig.~\ref{fig:urban}. 
    \end{bullets}%

    \blt[reformulation]
    \begin{bullets}%
        \blt[constraint]Replacing $\absloc_\absind \in \flyregion$
        in \eqref{eq:problemfpos} with $\absloc_\absind \in \grid$
        means that optimizing with respect to $\absnum$ and  the ABS positions
        $\{\absloc_\absind\}_{\absind=1}^\absnum$ is equivalent to choosing
        the smallest subset of points in $\grid$ for which there
        exists a feasible rate allocation, i.e., for which there
        exist
        $\{\rateuserabsinds\}_{\userind=1,\absind=1}^{\usernum,\absnum}$
        satisfying the constraints in \eqref{eq:problemf}.
        \blt[activations] Equivalently, one can consider all grid
        points by setting each $\absloc_\absind$ to be one of the
        grid points and then ``disable'' those $\absloc_\absind$
        where no ABS is going to be present. To this end, let
        $\acti_\gridptind$ be 1 if there is an ABS at
        $\gridpt_\gridptind$ and 0 otherwise, in which case the
        rates from that location need to be zero. By following this approach, Problem
        \eqref{eq:problemf} becomes
        \begin{subequations}
            \label{eq:problemfalpha}
            \begin{align}
                \minimize_{
                \{\acti_\gridptind\}_{\gridptind=1}^\gridptnum, \{\rateusergridptinds\}_{\userind=1,\gridptind=1}^{\usernum,\gridptnum}}~ & \sum_{\gridptind=1}^{\gridptnum}\acti_\gridptind \\
                \st~
                                                                                                                                          & \sum_\userind \rateusergridptinds \leq
                \maxratefun(\gridpt_\gridptind),
                \\                                                                                                                       & \sum_\gridptind \rateusergridptinds
                \geq \minrate,                                                                                                                                                               \\
                                                                                                                                          &
                0\leq \rateusergridptinds \leq \acti_\gridptind\capfun_\userind(\gridpt_\gridptind), \label{eq:problemfalphacap}
                \\                                                                                                                       & \acti_\gridptind \in \{0,1\},
            \end{align}
        \end{subequations}
        where $\rateusergridptinds$ is the rate that an ABS at $\gridpt_\gridptind$ would assign to the $\userind$-th user
        and now the constraints are for all $\userind$ and
        $\gridptind$. Note the presence of $\acti_\gridptind$ in the right-hand
        side of \eqref{eq:problemfalphacap}, which forces the rate to be 0 at
        grid points without an actual ABS.

        To simplify notation, it is convenient to express Problem
        \eqref{eq:problemfalpha} in terms of
        $\ratevec_\gridptind\define[\rate_\gridptind\entnot{1},\ldots,
                \rate_\gridptind\entnot{\usernum}]\transpose$,
        $\ratemat\define[\ratevec_1,\ldots,\ratevec_\gridptnum]$,
        $\maxratevec\define[ \maxratefun(\gridpt_1),\ldots,
                \maxratefun(\gridpt_\gridptnum)]\transpose$, and
        $\capvec_\gridptind\define[\capfun_1(\gridpt_\gridptind),\ldots,
                \capfun_\usernum(\gridpt_\gridptind)]\transpose$ as
        \begin{subequations}
            \label{eq:problemfalphamat}
            \begin{align}
                \minimize_{
                \{\acti_\gridptind\}_{\gridptind=1}^\gridptnum,\ratemat}~ & \sum_{\gridptind=1}^{\gridptnum}\acti_\gridptind                            \\
                \st~                                                      & \ratemat\transpose \bm 1 \leq \maxratevec\label{eq:problemfalphamatmaxrate} \\
                                                                          & \ratemat \bm 1 \geq \minrate \bm 1\label{eq:problemfalphamatminrate}
                \\
                                                                          & \bm
                0
                \leq
                \ratevec_\gridptind
                \leq\acti_\gridptind
                \capvec_\gridptind \label{eq:problemfalphamatrange}
                \\                                                                                                                       & \acti_\gridptind \in \{0,1\}.
            \end{align}
        \end{subequations}

        \blt[complexity]
        \begin{bullets}
            \blt[fixed activ] Observe that, if
            $\{\acti_\gridptind\}_{\gridptind=1}^\gridptnum$ were
            given, then Problem \eqref{eq:problemfalphamat} would be a
            convex feasibility problem in $\ratemat$.
            \blt[not fixed activ]However, since one needs to optimize
            over $\{\acti_\gridptind\}_{\gridptind=1}^\gridptnum$ as
            well, Problem \eqref{eq:problemfalphamat} becomes
            combinatorial. In fact, the following result establishes
            the NP-hardness of \eqref{eq:problemfalphamat}.

            \begin{mytheorem}
                \label{prop:nphard}
                Problem \eqref{eq:problemfalphamat} is NP-hard unless
                P$=$NP.
            \end{mytheorem}

            \begin{myproofsketch}
                The idea is to establish that a special case of
                \eqref{eq:problemfalphamat} is a multidimensional
                knapsack problem.

                Let the $\gridptind$-th entry of
                $\maxratevec$ be at least as large as $\bm 1\transpose
                    \capvec_\gridptind$ and note that, due to
                \eqref{eq:problemfalphamatrange}, constraint
                \eqref{eq:problemfalphamatmaxrate} holds regardless of
                the choice of
                $\{\acti_\gridptind\}_{\gridptind=1}^\gridptnum$ and
                $\ratemat$, meaning that
                \eqref{eq:problemfalphamatmaxrate} can be removed.  Then
                it can be easily shown that $\ratevec_\gridptind$ can be
                replaced with $\acti_\gridptind\capvec_\gridptind$.  The
                left-hand side of \eqref{eq:problemfalphamatminrate} can
                then be written as $\ratemat \bm 1=\sum_\gridptind
                    \ratevec_\gridptind = \sum_\gridptind
                    \acti_\gridptind\capvec_\gridptind$.  Finally, applying
                the change of variables $\invacti_\gridptind\leftarrow
                    1-\acti_\gridptind$, the objective becomes
                $\gridptnum-\sum_\gridptind \invacti_\gridptind$ and the
                left-hand side of \eqref{eq:problemfalphamatminrate}
                becomes $\sum_\gridptind
                    (1-\invacti_\gridptind)\capvec_\gridptind
                    =\sum_\gridptind \capvec_\gridptind - \sum_\gridptind
                    \invacti_\gridptind\capvec_\gridptind $. Applying these
                changes onto \eqref{eq:problemfalphamat} yields an
                instance of the so-called multidimensional knapsack
                problem, which has been shown to be NP-hard unless
                P$=$NP~\cite{gens1980complexity}. For an extended version of this proof,
                see Sec.~\ref{proof:nphard} of the supplementary
                material~\cite{romero2023placementarxiv}.

            \end{myproofsketch}


        \end{bullets}

        \blt[remove actis]
        It is useful to reduce Problem \eqref{eq:problemfalphamat} to
        the optimization with respect to $\ratemat$ alone. To this
        end, one can note that, for an arbitrary $\ratemat$, the
        corresponding optimal
        $\{\acti_\gridptind\}_{\gridptind=1}^\gridptnum$ satisfy that
        $\acti_\gridptind=0$ if $\ratevec_\gridptind=\bm 0$ and
        $\acti_\gridptind=1$ otherwise. This means that
        \eqref{eq:problemfalphamat} can be written as
        \begin{subequations}
            \begin{align}
                \minimize_{
                    \ratemat}~
                     & \sum_{\gridptind=1}^{\gridptnum}\indicator[\ratevec_\gridptind\neq
                \bm 0]                                                                    \\
                \st~ & \ratemat\transpose \bm 1 \leq \maxratevec                          \\
                     & \ratemat \bm 1 \geq \minrate \bm 1
                \\
                     & \bm
                0
                \leq
                \ratevec_\gridptind
                \leq \capvec_\gridptind,
            \end{align}
        \end{subequations}
        where $\indicator[\cdot]$ equals 1 if the condition inside
        brackets is true and 0 otherwise.

        \blt[objective]The objective
        $\sum_{\gridptind=1}^{\gridptnum}
            \indicator[\ratevec_\gridptind\neq \bm 0]$ can be equivalently
        expressed as
        $\sum_{\gridptind=1}^{\gridptnum}
            \indicator[\|\ratevec_\gridptind\|_\infty\neq  0]$, where the
        $\ell_\infty$-norm $\|\bm v\|_\infty$ equals the largest absolute
        value of the entries of vector $\bm v$. Clearly,
        $\sum_{\gridptind=1}^{\gridptnum}
            \indicator[\|\ratevec_\gridptind\|_\infty\neq  0] = \|
            [\|\ratevec_1\|_\infty,\ldots,
            \|\ratevec_\gridptnum\|_\infty]\transpose\|_0$, which suggests the
        relaxation
        $ \| [\|\ratevec_1\|_\infty,\ldots,
            \|\ratevec_\gridptnum\|_\infty]\transpose\|_1=\sum_\gridptind
            \|\ratevec_\gridptind\|_\infty$, or its reweighted version
        $\sum_\gridptind \sparsweight_\gridptind
            \|\ratevec_\gridptind\|_\infty$, where
        $\{\sparsweight_\gridptind\}_\gridptind$ are non-negative
        constants set as in~\cite{candes2008reweighted}.
        The problem thereby becomes
        \begin{subequations}
            \label{eq:mainproblem}
            \begin{align}
                \minimize_{
                    \ratemat}~
                     & \sum_\gridptind \sparsweight_\gridptind
                \|\ratevec_\gridptind\|_\infty                                         \\
                \st~ & \ratemat\transpose \bm 1 \leq \maxratevec                       \\
                     & \ratemat \bm 1 \geq \minrate \bm 1\label{eq:mainproblemminrate}
                \\
                     & \bm
                0
                \leq
                \ratemat
                \leq \capmat,
                \label{eq:mainproblemcap}
            \end{align}
        \end{subequations}
        \blt[interpretation]where the $(\userind,\gridptind)$-th
        entry of
        $\capmat\define[\capvec_1,\ldots,\capvec_\gridptnum]\in
            \rfield_+^{\usernum\times\gridptnum}$ is given by
        $\caps_{\userind,\gridptind}\define\capfun_\userind(\gridpt_\gridptind)$,
        i.e., the capacity of the link between the $\userind$-th
        GT and the $\gridptind$-th grid point. The
        $(\userind,\gridptind)$-th entry of $\ratemat$ therefore
        satisfies $0\leq \rate_{\userind,\gridptind}\leq
            \caps_{\userind,\gridptind}$, which means that it can be
        interpreted as the rate at which a \emph{virtual
            ABS}~\cite{huang2020sparse} placed at grid point
        $\gridpt_\gridptind$ communicates with the $\userind$-th
        GT. In case that $\rate_{\userind,\gridptind}=0$ for
        all $\userind$, then no \emph{actual} ABS needs to be
        deployed at $\gridpt_\gridptind$. In other words, the
        virtual ABS at $\gridpt_\gridptind$ corresponds to an
        actual ABS if and only if $\rate_{\userind,\gridptind}\neq 0$
        for some $\userind$.

    \end{bullets}
\end{bullets}

\section{Numerical Solver}
\label{sec:solver}

\begin{bullets}
    \blt[Overview] Observe that \eqref{eq:mainproblem} is a convex
    optimization problem and therefore it can be numerically solved in
    polynomial time.
    \blt[solvers]%
    \begin{bullets}%
        \blt[interior]In view of the inequality constraints, the first
        possibility one may consider is to apply an interior-point solver,
        as described in \supmatsec{app:interior}.
        \begin{bullets}%
            \blt[why not interior]Unfortunately, such an approach is only
            suitable for relatively small values of $\usernum$ and $\gridptnum$
            given the poor scalability of interior-point methods with the
            number of variables and constraints~\cite{lin2021admm}. Indeed, in
            this application, $\gridptnum$ can be in the order of tens of
            thousands, which would render the (at least
            cubic; cf. \supmatsec{app:interior}) complexity
            of interior-point methods prohibitive.
        \end{bullets}
        \blt[admm]In contrast, the
        rest of this section presents a solver whose complexity is linear
        in~$\gridptnum$ and $\usernum$ by building upon the so-called
        \emph{alternating-direction method of multipliers}
        (ADMM)~\cite{boyd2011distributed} and exploiting the special
        structure of the problem.

    \end{bullets}

\end{bullets}

\subsection{ADMM Decomposition}

\begin{bullets}%

    \blt[overview]As outlined below, ADMM alternately solves two
    optimization subproblems in the primal variables and performs a
    gradient step along the dual
    variables~\cite{boyd2011distributed}. But before decomposing
    \eqref{eq:mainproblem} into such subproblems, a few
    manipulations are in order.

    \blt[Simplification]
    \begin{bullets}
        \blt[equality]The first is to replace the inequality constraint
        \eqref{eq:mainproblemminrate} with an equality constraint. To this
        end, let $\ratervec_\userind$ denote the $\userind$-th column of
        $\ratemat\transpose$ and suppose that $\ratemat$ is feasible. Then,
        it follows from \eqref{eq:mainproblemminrate} that
        $\ratervec_\userind\transpose \bm 1\geq \minrate$. If one replaces
        $\ratervec_\userind$ with $\ratervec'_\userind \define
            (\minrate/(\ratervec_\userind\transpose \bm 1))\ratervec_\userind$,
        the entries of $\ratervec'_\userind$ are non-negative and not
        greater than the entries of $\ratervec_\userind$ since $0\leq
            (\minrate/(\ratervec_\userind\transpose \bm 1))\leq 1$. Hence, the resulting $\ratemat$ still satisfies all other constraints
        and the $\userind$-th constraint in \eqref{eq:mainproblemminrate}
        now holds with equality. Besides, the resulting objective  will
        not be greater. Therefore, if $\ratemat$ is optimal, scaling any of
        its rows in this way yields another optimal $\ratemat$. Applying
        this reasoning to all rows (i.e., for all $\userind$) shows that
        \eqref{eq:mainproblemminrate} can be replaced with an equality
        constraint without loss of optimality.

        \blt[objective]Second, due to \eqref{eq:mainproblemcap},
        the entries of $\ratevec_\gridptind$ are non-negative and,
        thus, $\|\ratevec_\gridptind\|_\infty$ equals the largest entry
        of $\ratevec_\gridptind$. This means that the objective can be
        replaced with $\sum_\gridptind \sparsweight_\gridptind
            \slack_\gridptind$ upon introducing  the auxiliary variables $
            \slack_\gridptind$ and  constraints $\ratevec_\gridptind\leq
            \slack_\gridptind\bm 1$ for each $\gridptind$.
        This gives rise to
        \begin{subequations}
            \label{eq:mainproblemeq}
            \begin{align}
                \minimize_{
                    \ratemat, \slackvec}~
                     & \sparsweightvec\transpose\slackvec        \\
                \st~ & \ratemat\transpose \bm 1 \leq \maxratevec \\
                     & \ratemat \bm 1 = \minrate \bm 1
                \\
                     & \bm
                0
                \leq
                \ratemat
                \leq \capmat                                     \\
                     & \ratemat \leq  \bm 1\slackvec\transpose,
            \end{align}
        \end{subequations}
        where
        $\sparsweightvec\define[\sparsweight_1,\ldots,\sparsweight_\gridptnum]\transpose$
        and $\slackvec\define [\slack_1,\ldots,\slack_\gridptnum]\transpose$.

    \end{bullets}

    \blt[formulation in ADMM form]The next step is to
    express \eqref{eq:mainproblemeq} in a form amenable to application of ADMM.
    \begin{bullets}%
        \blt[Theory]Specifically,  \eqref{eq:mainproblemeq} will be
        expressed in the homogeneous form
        \begin{subequations}
            \label{eq:admmstd}
            \begin{align}
                \minimize_{\Xmat, \Zmat}\quad & \admmfunx(\Xmat)+\admmfunz(\Zmat)                                               \\
                \st \quad                     & \Amat_1\Xmat \Amat_2+ \Bmat_1\Zmat\Bmat_2 = \bm 0,      \label{eq:admmstdconst}
            \end{align}
        \end{subequations}
        for which  the ADMM iteration becomes
        \cite[Sec.~3.1.1]{boyd2011distributed}
        \begin{subequations}
            \label{eq:admmstdit}
            \begin{align}
                \label{eq:admmstditx}
                \Xmat\itnot{\itind+1} & = \argmin_\Xmat \admmfunx(\Xmat) + \frac{\admmstep}{2}\| \Amat_1\Xmat \Amat_2+ \Bmat_1\Zmat\itnot{\itind}\Bmat_2 + \Umat\itnot{\itind}\|_\frob^2   \\
                \label{eq:admmstditz}
                \Zmat\itnot{\itind+1} & = \argmin_\Zmat \admmfunz(\Zmat) + \frac{\admmstep}{2}\| \Amat_1\Xmat\itnot{\itind+1} \Amat_2+ \Bmat_1\Zmat\Bmat_2 + \Umat\itnot{\itind}\|_\frob^2 \\
                \label{eq:admmstditu}
                \Umat\itnot{\itind+1} & =     \Umat\itnot{\itind} + \Amat_1\Xmat\itnot{\itind+1} \Amat_2+ \Bmat_1\Zmat\itnot{\itind+1}\Bmat_2
            \end{align}
        \end{subequations}
        for $\itind=1,2\ldots$.
        Here, $\Xmat\itnot{\itind}$ and $\Zmat\itnot{\itind}$ collect the
        primal variables,
        $\Umat\itnot{\itind}$ is a matrix of scaled dual variables, and
        $\admmstep>0$ is the step-size parameter.

        \blt[assignments]Each possible correspondence that one may establish
        between the variables, constants, and functions of \eqref{eq:admmstd}
        and those of \eqref{eq:mainproblemeq}  results in a 
        different ADMM algorithm.  Finding a good correspondence is typically
        the most critical step and takes multiple attempts since, unless properly
        accomplished,  the complexity of the  subproblems
        \eqref{eq:admmstditx} and \eqref{eq:admmstditz} will be comparable to
        the complexity of the original problem. For the problem at hand, the
        following assignment was found to yield subproblems
        that separate along the rows and columns of $\ratemat$:
        \begin{subequations}
            \label{eq:ubassignments}
            \begin{align}
                \Xmat            & \rightarrow [\ratemat\transpose , \slackvec]\transpose                                                                                                     \\
                \Zmat            & \rightarrow \ratemat                                                                                                                                       \\
                \admmfunx(\Xmat) & \rightarrow     \sparsweightvec\transpose \slackvec + \sum_\gridptind
                \indicatorinf[\ratevec_\gridptind\leq \slack_\gridptind
                    \bm 1] + \sum_\gridptind\indicatorinf[
                    \bm 1\transpose\ratevec_\gridptind \leq \maxrate{\gridptind}
                ]                                                                                                                                                                             \\
                \admmfunz(\Zmat) & \rightarrow \indicatorinf[\ratemat \bm 1 = \minrate \bm 1] + \indicatorinf[\bm 0 \leq \ratemat \leq \capmat]                                               \\
                \Amat_1          & \rightarrow [\bm I_\usernum, \bm 0], ~\Amat_2 \rightarrow \bm I_{\gridptnum }, ~\Bmat_1\rightarrow - \bm I_\usernum, ~\Bmat_2\rightarrow \bm I_\gridptnum.
            \end{align}
        \end{subequations}
        Here,
        $\maxrate{\gridptind}$ is the $\gridptind$-th entry of $\maxratevec$
        and  $\indicatorinf[\cdot]$ is a function that takes the value 0 when
        the condition inside brackets holds and $\infty$ otherwise.
        \blt[overall structure] Note that, given \eqref{eq:ubassignments}, it
        follows that $\Amat_1\Xmat \Amat_2+
            \Bmat_1\Zmat\Bmat_2=\ratemat-\Zmat$ and, therefore, \eqref{eq:admmstdconst} imposes that $\ratemat=\Zmat$.  Thus, intuitively
        speaking, each subproblem \eqref{eq:admmstditx} and
        \eqref{eq:admmstditz} tries to find values for their respective
        variables that satisfy the structure promoted by the first terms in
        \eqref{eq:admmstditx} and \eqref{eq:admmstditz} while, at the same time,
        the second terms in these expressions as well as \eqref{eq:admmstditu}
        push towards an agreement between the solutions of both subproblems.

        \blt[Overview] The next two subsections will be respectively concerned
        with finding the solutions of \eqref{eq:admmstditx} and
        \eqref{eq:admmstditz}. Afterwards, both solutions are put together
        in Sec.~\ref{sec:proposedsolver} to
        obtain the desired algorithm.

    \end{bullets}%

\end{bullets}%

\subsection{The $X$-subproblem}
\label{sec:xsubproblem}

\begin{bullets}%
    \blt[Overview]This section decomposes the $\Xmat$-update
    \eqref{eq:admmstditx} into $\gridptnum$ smaller problems. The latter
    can be efficiently solved by finding a root of a scalar equation
    through the bisection algorithm.

    \blt[Special expression] In view of \eqref{eq:ubassignments},
    \eqref{eq:admmstditx} can be expressed as
    \begin{subequations}
        \label{eq:admmstditxexp}
        \begin{align}
            (\ratemat\itnot{\itind+1},\slackvec\itnot{\itind+1}) = & \arg\min_{\ratemat, \slackvec}~
            \sparsweightvec\transpose \slackvec + \sum_\gridptind
            \indicatorinf[\ratevec_\gridptind\leq \slack_\gridptind \bm 1]
            \nonumber                                                                                \\ &+\sum_\gridptind\indicatorinf[
                \bm 1\transpose\ratevec_\gridptind \leq \maxrate{\gridptind}
            ]+ \frac{\admmstep}{2}\|\ratemat - \Zmat\itnot{\itind
            }+\Umat\itnot{\itind}\|_\frob^2
            \\
            =                                                      & \arg\min_{\ratemat, \slackvec}~
            \sum_\gridptind \big[ \sparsweight_\gridptind \slack_\gridptind +
                \indicatorinf[\ratevec_\gridptind\leq \slack_\gridptind \bm 1]
            \nonumber                                                                                \\&+ \indicatorinf[
                    \bm 1\transpose\ratevec_\gridptind \leq \maxrate{\gridptind}
                ]+ \frac{\admmstep}{2}\|\ratevec_\gridptind - \zvec_\gridptind\itnot{\itind
                }+\uvec_\gridptind\itnot{\itind}\|_2^2                                          \big],
        \end{align}
    \end{subequations}
    where $\zvec_\gridptind\itnot{\itind }$ and $\uvec_\gridptind\itnot{\itind}$
    respectively denote the $\gridptind$-th column of $\Zmat\itnot{\itind}$
    and $\Umat\itnot{\itind}$.
    \blt[separation] This problem clearly separates into
    $\gridptnum$ problems of the form
    \begin{subequations}
        \label{eq:ratevecindividualconst}
        \begin{align}
            (\ratevec_\gridptind\itnot{\itind+1},\slack_\gridptind\itnot{\itind+1})
            = \arg\min_{\ratevec_\gridptind, \slack_\gridptind}~ & \sparsweight_\gridptind \slack_\gridptind
            + \frac{\admmstep}{2}\|\ratevec_\gridptind - \zvec_\gridptind\itnot{\itind
            }+\uvec_\gridptind\itnot{\itind}\|_2^2                                                                                                                       \\
            \st\quad                                             & \ratevec_\gridptind\leq \slack_\gridptind \bm 1                                                       \\
                                                                 & \bm 1\transpose\ratevec_\gridptind \leq \maxrate{\gridptind}.\label{eq:ratevecindividualconstmaxrate}
        \end{align}
    \end{subequations}

    \blt[Two cases]There are two cases: C1) constraint
    \eqref{eq:ratevecindividualconstmaxrate} is active (i.e. it holds
    with equality) at the optimal solution; C2)
    \eqref{eq:ratevecindividualconstmaxrate} is inactive (i.e.  it
    holds with strict inequality) at the optimal solution. Thus, to
    solve \eqref{eq:ratevecindividualconst}, one can apply the
    following strategy. First, solve the problem that results from
    removing \eqref{eq:ratevecindividualconstmaxrate}:
    \begin{subequations}
        \label{eq:ratevecindividual}
        \begin{align}
            (\ratevec_\gridptind\itnot{\itind+1},\slack_\gridptind\itnot{\itind+1})
            = \arg\min_{\ratevec_\gridptind, \slack_\gridptind}~ & \sparsweight_\gridptind \slack_\gridptind
            + \frac{\admmstep}{2}\|\ratevec_\gridptind - \zvec_\gridptind\itnot{\itind
            }+\uvec_\gridptind\itnot{\itind}\|_2^2                                                                  \\
            \st\quad                                             & \ratevec_\gridptind\leq \slack_\gridptind \bm 1.
        \end{align}
    \end{subequations}
    If the solution to \eqref{eq:ratevecindividual} satisfies
    \eqref{eq:ratevecindividualconstmaxrate}, then it is also the optimal
    solution of \eqref{eq:ratevecindividualconst}. Else, due to the
    convexity of the problem, \eqref{eq:ratevecindividualconstmaxrate}
    must necessarily be active at the optimum. In this case, the optimal
    solution of \eqref{eq:ratevecindividualconst} can be found by
    replacing \eqref{eq:ratevecindividualconstmaxrate} with an equality
    constraint.
    \begin{bullets}
        \blt[without max rate constraint]
        Thus, let us start by solving \eqref{eq:ratevecindividual}.

        \begin{bullets}

            \blt[nec. and suff. conds] 

            \begin{myproposition}
                \thlabel{prop:eqxstep1}
                Let
                $\ratevec_\gridptind\itnot{\itind+1}$ and
                $\slack_\gridptind\itnot{\itind+1} $ be  given by \eqref{eq:ratevecindividual}. It holds that
                \begin{subequations}
                    \label{eq:xkkt}
                    \begin{align}
                        \label{eq:xkktr}
                         & \ratevec_\gridptind\itnot{\itind+1} = \min(\zvec_\gridptind\itnot{\itind
                        }-\uvec_\gridptind\itnot{\itind}, \slack_\gridptind\itnot{\itind+1}\bm 1)   \\
                        \label{eq:xkks}
                         & \bm 1\transpose \max(\zvec_\gridptind\itnot{\itind
                        }-\uvec_\gridptind\itnot{\itind} -\slack_\gridptind\itnot{\itind+1}\bm 1, \bm 0) = \frac{\sparsweight_\gridptind}{\admmstep},
                    \end{align}
                \end{subequations}
                where  $\min$ and $\max$ operate  entrywise.
            \end{myproposition}

            \begin{myproofsketch}
                Since Problem \eqref{eq:ratevecindividual} is convex differentiable
                and Slater's conditions are satisfied, it follows that the
                Karush-Kuhn-Tucker (KKT) conditions are sufficient and
                necessary~\cite[Sec.~5.5.3]{boyd}. Thus, one can start from the
                Lagrangian of \eqref{eq:ratevecindividual}, which is given by
                \begin{align}
                    \lagrangian(\ratevec_\gridptind, \slack_\gridptind;\nuvec)
                    = \sparsweight_\gridptind \slack_\gridptind
                    + \frac{\admmstep}{2}\|\ratevec_\gridptind - \zvec_\gridptind\itnot{\itind
                    }+\uvec_\gridptind\itnot{\itind}\|_2^2 + \nuvec\transpose(   \ratevec_\gridptind- \slack_\gridptind \bm 1),
                \end{align}
                and write the KKT conditions $
                    \nabla_{\ratevec_\gridptind}\lagrangian(\ratevec_\gridptind,
                    \slack_\gridptind;\nuvec) = \admmstep(\ratevec_\gridptind -
                    \zvec_\gridptind\itnot{\itind }+\uvec_\gridptind\itnot{\itind}) +
                    \nuvec = \bm 0$, $
                    \nabla_{\slack_\gridptind}\lagrangian(\ratevec_\gridptind,
                    \slack_\gridptind;\nuvec) = \sparsweight_\gridptind - \bm 1\transpose
                    \nuvec=0$, $ \ratevec_\gridptind\leq \slack_\gridptind \bm 1$, $
                    \nuvec\geq \bm 0$, and
                $\nus\entnot{\userind}(\rate_\gridptind\entnot{\userind} -
                    \slack_\gridptind) = 0~\forall \userind$. The latter equations
                constitute a non-linear system of equations that reduces, after some
                algebraic manipulations, to \eqref{eq:xkkt}.  See
                Sec.~\ref{proof:eqxstep1} of the supplementary material~\cite{romero2023placementarxiv} for details.
            \end{myproofsketch}


            \blt[solving the conditions]
            Observe that \eqref{eq:xkktr} can be used to obtain
            $\ratevec_\gridptind\itnot{\itind+1}$ if
            $\slack_\gridptind\itnot{\itind+1}$ is given, whereas \eqref{eq:xkks}
            does not depend on $\ratevec_\gridptind\itnot{\itind+1}$. Therefore, a
            solution to \eqref{eq:xkkt} can be found by first solving
            \eqref{eq:xkks} for $\slack_\gridptind\itnot{\itind+1}$ and then
            substituting the result into \eqref{eq:xkktr}
            to recover $\ratevec_\gridptind\itnot{\itind+1}$. To this end, we have
            the following:

            \begin{myproposition}
                \thlabel{prop:rootssg}
                Equation \eqref{eq:xkks} has a unique root. This root lies in the interval $[\slacklow\gii, \slackhigh\gii]$, where
                \begin{subequations}
                    \begin{align}
                        \label{eq:slacklow}
                        \slacklow\gii  & \define \min_\userind\left( \zs\gii\entnot{\userind} - \us\gii\entnot{\userind}\right) - \frac{\sparsweight_\gridptind}{\usernum \admmstep}  \\
                        \label{eq:slackhigh}
                        \slackhigh\gii & \define \max_\userind\left( \zs\gii\entnot{\userind} - \us\gii\entnot{\userind}\right) - \frac{\sparsweight_\gridptind}{\usernum \admmstep}.
                    \end{align}
                \end{subequations}

            \end{myproposition}


            \begin{myproofsketch}
                Consider the function $\Ffun(\slack)\define \bm 1\transpose
                    \max(\zvec_\gridptind\itnot{\itind }-\uvec_\gridptind\itnot{\itind}
                    -\slack\bm 1, \bm 0) = \sum_\userind
                    \max(\zs_\gridptind\itnot{\itind
                    }\entnot{\userind}-\us_\gridptind\itnot{\itind}\entnot{\userind}
                    -\slack, 0)$. Since $\Ffun$ is the sum of non-increasing
                piecewise linear functions, so is $\Ffun$. Since
                $\Ffun(\slack)\rightarrow \infty$ as $\slack\rightarrow -\infty$ and
                $\Ffun(\slack)=0$ for a sufficiently large $\slack$, it follows that
                \eqref{eq:xkks} has at least one root. Uniqueness of the root
                follows readily by noting that $\Ffun$ is strictly decreasing
                whenever $\Ffun(\slack)>0$.  The proof is easily completed by
                establishing that
                $\Ffun(\slacklow\gii)\geq{\sparsweight_\gridptind}/
                            \admmstep$ whereas $\Ffun(\slackhigh\gii)\leq
                    {\sparsweight_\gridptind}/ \admmstep$.  Further details
                are provided in Sec.~\ref{proof:rootssg} of the supplementary
                material~\cite{romero2023placementarxiv}.
            \end{myproofsketch}

            \blt[interpretation]Observe that \thref{prop:rootssg} essentially
            provides the bounds that are required to find the unique root of
            \eqref{eq:xkks} via the well-known bisection algorithm. Recall
            that the latter is a very efficient algorithm as it geometrically
            reduces the uncertainty at every iteration.

        \end{bullets}%

        \blt[max rate constraint is active]In accordance with the strategy
        mentioned earlier, one also necessitates a means to solve
        \eqref{eq:ratevecindividualconst} when
        \eqref{eq:ratevecindividualconstmaxrate} is replaced with an equality
        constraint:
        \begin{subequations}
            \label{eq:ratevecindividualconsteq}
            \begin{align}
                (\ratevec_\gridptind\itnot{\itind+1},\slack_\gridptind\itnot{\itind+1})
                = \arg\min_{\ratevec_\gridptind, \slack_\gridptind}~ & \sparsweight_\gridptind \slack_\gridptind
                + \frac{\admmstep}{2}\|\ratevec_\gridptind - \zvec_\gridptind\itnot{\itind
                }+\uvec_\gridptind\itnot{\itind}\|_2^2                                                                                                                      \\
                \st\quad                                             & \ratevec_\gridptind\leq \slack_\gridptind \bm 1                                                      \\
                                                                     & \bm 1\transpose\ratevec_\gridptind = \maxrate{\gridptind}.\label{eq:ratevecindividualconstmaxrateeq}
            \end{align}
        \end{subequations}
        \begin{bullets}
            \blt[nec. and suff. conds]
            \begin{myproposition}
                \thlabel{prop:eqxstep2}
                Let      $\ratevec_\gridptind\itnot{\itind+1}$ and
                $\slack_\gridptind\itnot{\itind+1} $ be given by \eqref{eq:ratevecindividualconsteq}. Then, it holds that
                \begin{subequations}
                    \label{eq:eqxkkt}
                    \begin{align}
                        \label{eq:eqxkktr}
                         & \ratevec_\gridptind\itnot{\itind+1} = \min(\zvec_\gridptind\itnot{\itind
                        }-\uvec_\gridptind\itnot{\itind} - (\mumult/\admmstep)\bm 1, \slack_\gridptind\itnot{\itind+1}\bm 1) \\
                        \label{eq:eqxkks}
                         & \bm 1\transpose\max(  \mumult\bm 1, \admmstep(\zvec_\gridptind\itnot{\itind
                            }-\uvec_\gridptind\itnot{\itind}-\slack_\gridptind\itnot{\itind+1}\bm
                            1))  = \sparsweight_\gridptind + \mumult\usernum,
                    \end{align}
                \end{subequations}
                where
                \begin{align}
                    \label{eq:eqmumult}
                    \mumult
                    \define\frac{ -\admmstep\maxrate{\gridptind}
                        +\admmstep\bm 1\transpose ( \zvec_\gridptind\itnot{\itind
                        }-\uvec_\gridptind\itnot{\itind})-\sparsweight_\gridptind}{\usernum}.
                \end{align}

            \end{myproposition}


            \begin{myproofsketch}
                The proof follows along the same lines as the proof of
                \thref{prop:eqxstep1}.  See Sec.~\ref{proof:eqxstep2} of the
                supplementary material~\cite{romero2023placementarxiv} for more details.
            \end{myproofsketch}

            \blt[solving the conditions] Observe that \eqref{eq:eqxkktr} and
            \eqref{eq:eqmumult} can be used to obtain
            $\ratevec_\gridptind\itnot{\itind+1}$ if
            $\slack_\gridptind\itnot{\itind+1}$ is given, whereas \eqref{eq:eqxkks}
            does not depend on $\ratevec_\gridptind\itnot{\itind+1}$. Therefore, a
            solution to \eqref{eq:ratevecindividualconsteq} can be found by first solving
            \eqref{eq:eqxkks} for $\slack_\gridptind\itnot{\itind+1}$ and then
            substituting the result into  \eqref{eq:eqxkktr}
            to recover $\ratevec_\gridptind\itnot{\itind+1}$. To this end, we have
            the following:

            \begin{myproposition}
                \thlabel{prop:eqrootssg}
                Equation \eqref{eq:eqxkks} has a unique root. This root lies in the interval $[\slacklow\gii, \slackhigh\gii]$, where
                \begin{subequations}
                    \begin{align}
                        \label{eq:eqslacklow}
                        \slacklow\gii  & \define \min_\userind\left( \zs\gii\entnot{\userind} - \us\gii\entnot{\userind}\right) - \frac{\sparsweight_\gridptind}{\usernum \admmstep}-\frac{\mumult}{\admmstep}  \\
                        \label{eq:eqslackhigh}
                        \slackhigh\gii & \define \max_\userind\left( \zs\gii\entnot{\userind} - \us\gii\entnot{\userind}\right) - \frac{\sparsweight_\gridptind}{\usernum \admmstep}-\frac{\mumult}{\admmstep}.
                    \end{align}
                \end{subequations}

            \end{myproposition}



            \begin{myproofsketch}
                The proof follows along the same lines as the proof of
                \thref{prop:rootssg}. Details can be found in
                Sec.~\ref{proof:eqrootssg} of the supplementary material~\cite{romero2023placementarxiv}.
            \end{myproofsketch}

        \end{bullets}
    \end{bullets}
    \blt[summary]To sum up, \eqref{eq:admmstditxexp} can be solved separately
    for each column of $\ratemat$ and entry of $\slackvec$. Each of these
    $\gridptnum$ problems can be solved by first solving \eqref{eq:ratevecindividual}
    and checking whether \eqref{eq:ratevecindividualconstmaxrate} holds
    for the obtained solution. If it does not hold,
    one must solve
    \eqref{eq:ratevecindividualconsteq}. \thref{prop:eqxstep1} and
    \thref{prop:eqxstep2} respectively establish that a solution can be
    found for each of these problems just by solving a scalar
    equation.  \thref{prop:rootssg} and  \thref{prop:eqrootssg}
    prove uniqueness of the solutions of these two equations and provide
    an interval where they can be sought using the bisection algorithm.


\end{bullets}

\subsection{The $Z$-subproblem}

\begin{bullets}%
    \blt[Overview] This section describes how a solution to
    \eqref{eq:admmstditz} can be found by solving a bisection problem per
    row of $\Zmat$.
    \blt[Specialization]To this end, start by noting that  it follows from \eqref{eq:admmstditz} and \eqref{eq:ubassignments} that
    \begin{subequations}
        \begin{align}
            \Zmat\itnot{\itind+1} = \argmin_\Zmat ~~ & \bigg[ \indicatorinf[\Zmat \bm 1 = \minrate \bm 1] + \indicatorinf[\bm 0 \leq \Zmat \leq \capmat]                                                \\&+ \frac{\admmstep}{2}\| \ratemat\itnot{\itind+1} -\Zmat + \Umat\itnot{\itind}\|_\frob^2\bigg]\\
            = \argmin_\Zmat ~~                       & \sum_\userind\bigg[ \indicatorinf[\zrvec_\userind\transpose \bm 1 = \minrate] + \indicatorinf[\bm 0 \leq \zrvec_\userind \leq \caprvec_\userind] \\&+ \frac{\admmstep}{2}\| \ratervec_\userind\itnot{\itind+1} -\zrvec_\userind + \urvec_\userind\itnot{\itind}\|_\frob^2\bigg],
        \end{align}
    \end{subequations}
    where $\zrvec_\userind$, $\caprvec_\userind$, and $\urvec_\userind\itnot{\itind}$ respectively denote the $\userind$-th column of $\Zmat\transpose$, $\capmat\transpose$, and $(\Umat\itnot{\itind})\transpose$.
    \blt[separation] Clearly, this separates into $\usernum$ problems of the form
    \begin{subequations}
        \label{eq:ubzrvecs}
        \begin{align}
            \zrvec_\userind\itnot{\itind+1} = \argmin_{\zrvec_\userind} ~~ & \frac{1}{2}\| \ratervec_\userind\itnot{\itind+1} -\zrvec_\userind + \urvec_\userind\itnot{\itind}\|_\frob^2 \\
            \st \quad                                                      & \bm 1\transpose \zrvec_\userind = \minrate,~~ \bm 0 \leq \zrvec_\userind \leq \caprvec_\userind.
        \end{align}
    \end{subequations}

    \blt[Nec. and suff. conds]
    \begin{myproposition}
        \label{prop:zsol}
        If
        $\bm 1\transpose \caprvec_\userind< \minrate$, then
        \eqref{eq:ubzrvecs} is infeasible.  If
        $\bm 1\transpose \caprvec_\userind\geq \minrate$, the solution to
        \eqref{eq:ubzrvecs} is given by
        \begin{align}
            \label{eq:ubzrvecsol}
            \zrvec_\userind\itnot{\itind+1} =
            \max(\bm 0,    \min(\caprvec_\userind, \ratervec_\userind\itnot{\itind+1} +
                \urvec_\userind\itnot{\itind} - \lambdas\bm 1)),
        \end{align}
        where $\lambdas$ satisfies
        \begin{align}
            \label{eq:ubzrveclambda}
            \bm 1\transpose\max(\bm 0,    \min(\caprvec_\userind, \ratervec_\userind\itnot{\itind+1} +
                \urvec_\userind\itnot{\itind} - \lambdas\bm 1)) = \minrate.
        \end{align}

    \end{myproposition}


    \begin{myproofsketch}
        The proof follows along the same lines as the proof of
        \thref{prop:eqxstep1}.  See Sec.~\ref{proof:zsol} of the
        supplementary material~\cite{romero2023placementarxiv} for more details.
    \end{myproofsketch}

    \blt[solving the conditions] Thus, as in Sec.~\ref{sec:xsubproblem}, one needs
    to solve the scalar equation \eqref{eq:ubzrveclambda}. The following
    result is the counterpart of \thref{prop:rootssg} for the
    $\Zmat$-subproblem.

    \begin{myproposition}
        \thlabel{prop:ubrootlambda} If
        $\bm 1\transpose \caprvec_\userind< \minrate$, then equation
        \eqref{eq:ubzrveclambda} has no roots. If
        $\bm 1\transpose \caprvec_\userind\geq \minrate$, then
        \eqref{eq:ubzrveclambda} has a unique root. This root lies in the
        interval
        $[\lambdaslow_\userind\itnot{\itind},
                    \lambdashigh_\userind\itnot{\itind}]$, where
        \begin{subequations}
            \begin{align}
                \label{eq:ublambdaslow}
                \lambdaslow_\userind\itnot{\itind}=  &
                \min_\gridptind[\raters_\userind\itnot{\itind+1}\entnot{\gridptind} +
                    \urs_\userind\itnot{\itind}\entnot{\gridptind} -\caprs_\userind\entnot{\gridptind}]
                \\
                \label{eq:ublambdashigh}
                \lambdashigh_\userind\itnot{\itind}= &
                \max\{ \raters_\userind\itnot{\itind+1}\entnot{\gridptind} +
                \urs_\userind\itnot{\itind}\entnot{\gridptind}:\gridptind\in                 \{
                \gridptind: \caprs_\userind\entnot{\gridptind}>\frac{\minrate}{\gridptnum}
                \}
                \}\dc{\nonumber                        \\&}-\frac{\minrate}{\gridptnum}
            \end{align}
        \end{subequations}
    \end{myproposition}



    \begin{myproofsketch}
        The proof follows along the same lines as the proof of
        \thref{prop:rootssg}. Details can be found in
        Sec.~\ref{proof:ubrootlambda} of the supplementary material~\cite{romero2023placementarxiv}.
    \end{myproofsketch}

    \blt[interpretation]\thref{prop:ubrootlambda} establishes uniqueness
    and provides the bounds needed to solve \eqref{eq:ubzrveclambda},
    and therefore \eqref{eq:ubzrvecs}, via the bisection algorithm.

\end{bullets}%


\subsection{The Proposed Solver}
\label{sec:proposedsolver}


\newcommand{\widemargin}[1]{\parbox{\dimexpr\textwidth-2\algomargin\relax}{#1}}

\begin{algorithm}[t]
  \SetAlgoLined
  \textbf{Input:}~{$\capmat\in \rfield_+^{\usernum \times \gridptnum}$,
    $\minrate\in \rfield_+$, $\{\sparsweight_\gridptind\}_\gridptind\subset
    \rfield_+$,$\{\maxrate{\gridptind}\}_\gridptind\subset
    \rfield_+$,~$\admmstep>0$}\\
  \textbf{Initialize}~ $\Umat\itnot{1}\in \rfield^{\usernum \times \gridptnum}$ and $\Zmat\itnot{1}\in \rfield_+^{\usernum \times \gridptnum}$ \\
  \For{$\itind=1,2,\ldots$}{
    \For{$\gridptind=1,2,\ldots,\gridptnum$}{
      Bisection: find $\slack_\gridptind\itnot{\itind+1}$ s.t. \widemargin{$\bm 1\transpose \max(\bm 0, \admmstep(\zvec_\gridptind\itnot{\itind
      }-\uvec_\gridptind\itnot{\itind} -\slack_\gridptind\itnot{\itind+1}\bm 1)) = {\sparsweight_\gridptind}$}\\
      Set $\ratevec_\gridptind\itnot{\itind+1} = \min(\zvec_\gridptind\itnot{\itind
      }-\uvec_\gridptind\itnot{\itind}, \slack_\gridptind\itnot{\itind+1}\bm 1)$\\
      \If{$\bm 1\transpose\ratevec_\gridptind\itnot{\itind+1}>\maxrate{\gridptind}$}{
        Set $\mumult
    =({ -\admmstep\maxrate{\gridptind}
      +\admmstep\bm 1\transpose ( \zvec_\gridptind\itnot{\itind
  }-\uvec_\gridptind\itnot{\itind})-\sparsweight_\gridptind})/{\usernum}$\\
Bisection: find $\slack_\gridptind\itnot{\itind+1}$ s.t. \widemargin{
$\bm 1\transpose\max(  \mumult\bm 1, \admmstep(\zvec_\gridptind\itnot{\itind
  }-\uvec_\gridptind\itnot{\itind}-\slack_\gridptind\itnot{\itind+1}\bm
            1))  = \sparsweight_\gridptind + \mumult\usernum$
}  \\
        Set $\ratevec_\gridptind\itnot{\itind+1} = \min(\zvec_\gridptind\itnot{\itind
      }-\uvec_\gridptind\itnot{\itind} - (\mumult/\admmstep)\bm 1, \slack_\gridptind\itnot{\itind+1}\bm 1)$
}    
    }
    \For{$\userind=1,2,\ldots,\usernum$}{
      Bisection: find $\lambdas$~s.t.
      \widemargin{$\bm 1\transpose\max(\bm 0,    \min(\caprvec_\userind, \ratervec_\userind\itnot{\itind+1} +
        \urvec_\userind\itnot{\itind} - \lambdas\bm 1)) =\minrate $}\\[1ex]
      \widemargin{Set $    \zrvec_\userind\itnot{\itind+1} =
\max(\bm 0,    \min(\caprvec_\userind, \ratervec_\userind\itnot{\itind+1} +
    \urvec_\userind\itnot{\itind} - \lambdas\bm 1))
$}\\
}
Set $\Umat\itnot{\itind+1} = \Umat\itnot{\itind} +
    \ratemat\itnot{\itind+1}  -
    \Zmat\itnot{\itind+1}$\\
    \textbf{If} convergence(~) \textbf{then}  return $\ratemat\itnot{\itind+1}$ \label{step:convergence}
  }
  \caption{Group-sparse Placement Algorithm (GSPA).}
\label{algo:placement}    
\end{algorithm}

\begin{bullets}
    \blt[subproblems] Having addressed both the $\Xmat$- and $\Zmat$-subproblems, it remains only
    \blt[U-step]to obtain the $\Umat$-update in
    \eqref{eq:admmstditu}, which for the assignments in
    \eqref{eq:ubassignments} becomes
    \begin{align}
        \Umat\itnot{\itind+1} = \Umat\itnot{\itind} +
        \ratemat\itnot{\itind+1}  -
        \Zmat\itnot{\itind+1}.
    \end{align}

    \blt[algo]This completes the derivation of the proposed scheme,
    summarized as Algorithm~\ref{algo:placement} and referred to as
    the \textit{group-sparse placement algorithm} (GSPA) since it
    promotes group sparsity~\cite{yuan2006grouplasso} in the columns of
    $\ratemat$, that is, only a few columns of the matrix
    $\ratemat\itnot{\itind+1}$ returned by the algorithm are expected to
    be non-zero. Recall that the non-zero columns  indicate which grid
    points will be occupied by ABSs.
    \begin{bullets}%
        \blt[notation]             In the notation used in Algorithm~\ref{algo:placement}, if $\bm A$
        is a matrix, then $\bm a_m$ is its $m$-th column and $\bbm
            a_n\transpose$ its $n$-th row. Furthermore, superscripts indicate
        the iteration index, $\admmstep>0$ is the step size, and the $\min$
        and $\max$ operators act entrywise.
        \blt[stopping] The criterion on line \ref{step:convergence},
        which determines whether the algorithm has converged, is detailed
        in \supmatsec{sec:stopping}.
    \end{bullets}

    \blt[complexity] Observe that the main strategy in the previous two
    subsections was to exploit the structure of
    Problem~\eqref{eq:mainproblem} to decompose it into one subproblem
    per row and column of $\ratemat$.  Each of these subproblems involves
    solving a bisection task of a 1D monotonically decreasing function
    and therefore, given a target accuracy, can be solved with $\mathcal{O}(1)$ evaluations. The
    total complexity is $\mathcal{O}(\gridptnum\usernum)$, much smaller
    than the $\mathcal{O}(\gridptnum^3\usernum^3)$ complexity per inner
    iteration of  interior-point methods;
    cf.~\supmatsec{app:interior}.

\end{bullets}

\subsection{Convergence Analysis}
\label{sec:convergence}
\begin{changes}
    To prove the convergence of Algorithm~\ref{algo:placement}, it will be shown that the optimization problem in \eqref{eq:admmstd} with the assignments in \eqref{eq:ubassignments} satisfies the sufficient conditions in \cite[Sec.~3.2]{boyd2011distributed}, in which case ADMM converges.

    \begin{bullets}
        \blt[epigraphs]The first condition is that the epigraph of $\admmfunx(\Xmat) + \admmfunz(\Zmat)$ must be a closed nonempty convex set. This readily follows by noting that this function is the sum of a linear function and the indicators of closed nonempty convex sets.

        \blt[saddle point]The second condition is that the unaugmented Lagrangian
        \begin{equation*}
            \lagrangian_0(\Xmat,\Zmat,\Ymat) = \admmfunx(\Xmat) + \admmfunz(\Zmat) + \Ymat\transpose(\Amat_1\Xmat \Amat_2+ \Bmat_1\Zmat\Bmat_2)
        \end{equation*}
        has a saddle point, meaning that $\exists\Xmat^*\in \rfield^{(\usernum+1)\times\gridptnum}$, $\Zmat^*\in \rfield^{\usernum\times\gridptnum}$, and  $\Ymat^*\in \rfield^{\usernum\times\gridptnum}$ such that
        \begin{equation}
            \lagrangian_0(\Xmat^*, \Zmat^*, \Ymat) \leq \lagrangian_0(\Xmat^*, \Zmat^*, \Ymat^*) \leq \lagrangian_0(\Xmat, \Zmat, \Ymat^*)
        \end{equation}
        for all $\Xmat,\Zmat,\Ymat.$ This will be shown by establishing that the
        primal and dual optima of \eqref{eq:admmstd}  are attained~\cite[Ch. 5]{bertsekas1999}.
        \begin{bullets}
            \blt[Primal problem] To establish that the primal optimum is attained, just note that  \eqref{eq:admmstd} is equivalent to \eqref{eq:mainproblemeq}. The last constraint in \eqref{eq:mainproblemeq} can be replaced with the $\gridptnum$ equality constraints $\slack_\gridptind = \|\ratevec_\gridptind\|_\infty$,  $\gridptind =1,\ldots, \gridptnum$. The feasible set is therefore bounded and closed, which implies that it is compact. From Weierstrass' Theorem~\cite[Prop.~A8]{bertsekas1999}, the primal optimum is attained.

            \newcommand{\varvec}{\hc{\bm x}}
            \newcommand{\linvec}{\hc{\bm a}}
            \newcommand{\linmat}{\hc{\bm A}}
            \newcommand{\fset}{\hc{\mathcal S}}
            \newcommand{\dualfun}{\hc{g}}
            \newcommand{\dualvec}{\hc{\bm y}}
            \newcommand{\dualvecdir}{\hc{\bbm y}}
            \newcommand{\dualvecamp}{\hc{\alpha}}
            \newcommand{\dualvecset}{\hc{\mathcal{Y}}}
            \renewcommand{\indicator}{\hc{\mathcal{I}}}

            \blt[dual problem] To establish that the dual optimum is attained, it is worth simplifying the notation. To this end, observe that,
            upon letting $\varvec \define [\vect(\Xmat)\transpose, \vect(\Zmat)\transpose]\transpose$,
            \eqref{eq:admmstd} is of the general form
            \begin{salign}
                \minimize_{\varvec}\quad & \linvec\transpose\varvec + \indicator[\varvec\in\fset]\\
                \st \quad                                   & \linmat\varvec=\bm 0
            \end{salign}
            for some column vector $\linvec$, matrix $\linmat$, and set $\fset$. The dual function can therefore be expressed as
            \begin{salign}
                \dualfun(\dualvec) &= \inf_{\varvec}\left\{\linvec\transpose\varvec + \indicator[\varvec\in\fset] + \dualvec\transpose\linmat\varvec\right\}\\
                \label{eq:convproofdualfun}
                &= \inf_{\varvec \in \fset}\left\{\linvec\transpose\varvec + \dualvec\transpose\linmat\varvec\right\},
            \end{salign}
            where $\dualvec$ is the dual variable. Now let $
                \dualvecset \define \{ \dualvec~:~\dualfun(\dualvec) > -\infty\}$.
            Since strong duality holds for \eqref{eq:admmstd}, it follows that $\dualvecset$ is nonempty. Since $\dualfun$ is concave, it follows that $\dualvecset$ is convex.

            If the optimal $\dualvec$ is not attained, there must exist a vector
            $\dualvecdir$ and a sequence $\dualvecamp_n\in \rfield_+$  such that
            $\lim_{n\rightarrow \infty} \dualvecamp_n = \infty$ and the sequence
            $\dualfun(\dualvecamp_n \dualvecdir)$ is strictly increasing (it
            follows from concavity of $\dualfun$ and   Weierstrass'
            theorem~\cite[Prop.~A8]{bertsekas1999}).  The rest of the proof will show that the existence of such a sequence would result in a contradiction.

            Observe that
            there are two
            classes of directions $\dualvecdir$: Those for which any sequence
            $\dualvec_n=\dualvecamp_n \dualvecdir$
            with $\dualvecamp_n\rightarrow \infty $ satisfies that  $\exists
                n_0$ such that $\dualvec_n \in \dualvecset$ for all $n\geq n_0$ and
            those for which this condition does not hold. Directions of the latter type do not result in an increasing $\dualfun(\dualvecamp_n \dualvecdir)$, so they can be ruled out.

            Considering a direction of the first type, observe that $\dualfun(\dualvec_n)>-\infty $
            for $n\geq n_0$
            and, therefore, there exists a finite $\varvec^*$ such that
            $\dualfun(\dualvec_n) = \linvec\transpose\varvec^* + \dualvec_n\transpose\linmat\varvec^*$. Since the argument of  the infimum in \eqref{eq:convproofdualfun}    is linear in $\varvec$, such an $\varvec^*$ can be found on the boundary of $\fset$. Since $\fset$ is defined by affine equality and inequality constraints, following the same arguments as in the fundamental theorem of linear programming \cite{tardella2011linearprogramming}, $\varvec^*$ can be set to be one of the finitely many vertices $\{\varvec_i\}$ of $\fset$, which implies that
            $\dualfun(\dualvec_n) =\min_i [\linvec\transpose\varvec_i + \dualvec_n\transpose\linmat\varvec_i]
                =\min_i [\linvec\transpose\varvec_i + \dualvecamp_n\dualvecdir\transpose\linmat\varvec_i]
            $. For $\dualfun(\dualvec_n)$ to be increasing, it is necessary that $\dualvecdir\transpose\linmat\varvec_i>0$ for all $i$. However, in this case, it would follow that $\dualfun(\dualvec_n)\rightarrow \infty$, which contradicts the fact that $\dualfun(\dualvec)$ is upper bounded due to duality.

        \end{bullets}

    \end{bullets}

\end{changes}

\section{Extensions}
\label{sec:extensions}
\begin{changes}
    The present section proposes two auxiliary algorithms that complement  GSPA. First, Sec.~\ref{sec:minnumconnections} describes a scheme that minimizes the number of connections between the ABSs and the GTs. Second, Sec.~\ref{sec:masnumservedusers} outlines an optimization approach that can be used to   maximize the number of served users in the presence of a mismatch between the true and estimated radio maps.
\end{changes}

\subsection{Minimization of the Number of Connections}
\label{sec:minnumconnections}
\begin{changes}
    \begin{bullets}%
        \blt[introduction] The problem that Algorithm~\ref{algo:placement} aims to solve focuses on minimizing the number of ABSs. To attain a low number of ABSs, the scheme allows the GTs to be simultaneously served by multiple ABSs. However, after a low number of ABSs is obtained, it is desirable to minimize the number of connections between the ABSs and the GTs given the resulting placement.

        \blt[minimize num. connections]Thus, once the ABSs have been placed according to GSPA, the problem can be formulated as
        \begin{salign}[eq:minnumconnectionsprobx]                      \minimize_{\ratesubmat\in\rfield^{\usernum\times\absnum}}\quad & \sum_{\userind=1}^\usernum\sum_{\absind=1}^{\absnum}\indicator[\bar r_{\userind,\absind}\neq 0] \\
            \st \quad                                                      & \bm 1\transpose\ratesubmat \geq \minrate\bm 1,                                                  \\
            & \ratesubmat\transpose\bm 1 \leq  \maxratesubvec,                                                \\
            & \bm 0 \leq \ratesubmat \leq \bar{\capmat},
        \end{salign}
        where  $\ratesubmat\in\rfield^{\usernum\times\absnum}$ is a matrix whose $(\userind,\absind)$-th entry  $\bar r_{\userind,\absind}$  is the rate that the $\absind$-th ABS allocates to the $\userind$-th GT,
        $\ratervec_\userind\in \rfield^{\absnum}$ is the $\userind$-th column of $\ratesubmat\transpose$,
        $\maxratesubvec\in\rfield^{\absnum}$ collects the backhaul capacity for each ABS, and matrix $\bar{\capmat}\in\rfield^{\usernum\times\absnum}$ provides the capacity between each GT and each ABS. The objective function, which equals the number of non-zero entries of $\ratesubmat$, is the number of connections between the ABSs and the GTs.
        Since \eqref{eq:minnumconnectionsprobx} is again a combinatorial problem, it is useful to relax it as
        \begin{align}
            \minimize_{\ratesubmat\in\rfield^{\usernum\times\absnum}}\quad & \bm 1\transpose\ratesubmat\bm 1                                                                                                                                                                      \\
            \st \quad                                                      & \bm 1\transpose\ratesubmat \geq \minrate\bm 1,          ~                                      \ratesubmat\transpose\bm 1 \leq  \maxratesubvec, ~ \bm 0 \leq \ratesubmat \leq \bar{\capmat},\nonumber
        \end{align}
        which is a linear program and can be solved via standard methods. Note that, due to the nature of the above convex relaxation,  a greater reduction in the number of connections can be attained via reweighting, along the same lines of \eqref{eq:mainproblem}.
    \end{bullets}%
\end{changes}

\subsection{Coping with Radio Map Mismatches}
\label{sec:masnumservedusers}

\begin{changes}
    \begin{bullets}
        \blt[motivation \ra mismatch]Recall that GSPA uses a radio map for predicting the capacity of the ABS-GT channels and a (possibly separate) radio map for predicting the capacity of the backhaul links. Although the ABSs may collect measurements during their flights to update the radio maps, a certain mismatch between the radio maps used by GSPA and the true channel gain must be expected in practice. In the case of tomographic radio maps, this is caused by (i) the finiteness of the number of measurements used to estimate the radio map and (ii) the possible changes in the environment between the time when the radio map was constructed and when GSPA is used. In the case of ray-tracing maps, the mismatch is caused by the differences between the 3D model of the environment adopted by the ray-tracing algorithm and the actual environment. As a result, once the ABSs are placed at the locations indicated by GSPA, the capacity of the backhaul links as well as the channels between the ABSs and each GT will generally differ from what was expected. If the map estimates are too optimistic, the true capacity will be lower than expected and, therefore, some GTs cannot be served. Conversely, if the map estimates are too pessimistic, a larger number of ABSs than necessary will be deployed.
        \blt[overview]In either case, one must decide the rate that each ABS  allocates to each GT. Following the motivation behind GSPA, the goal will be to maximize the number of served users. To this end, this section proposes a simple rate allocation scheme based on a convex optimization problem.

        \blt[setup]Suppose that the $\absnum$ ABSs already adopted the positions provided by Algorithm~\ref{algo:placement} and measured the actual backhaul capacity and the actual capacity between each GT and each ABS. Consider the notation in Sec.~\ref{sec:minnumconnections}, with the difference that  $\maxratesubvec\in\rfield^{\absnum}$ and  $\bar{\capmat}\in\rfield^{\usernum\times\absnum}$ now contain these measured values rather than the values provided by the radio map.
        \blt[problem formulation] Given $\maxratesubvec$, $\bar{\capmat}$, and $\minrate$, which is the minimum rate that a GT must receive to be regarded as served,  the problem is to  minimize the number of  users that are not served:
        \begin{subequations}
            \begin{align}
                \minimize_{\ratesubmat\in\rfield^{\usernum\times\absnum}}\quad & \sum_{\userind}\indicator[\ratervec_\userind\transpose\bm 1 < \minrate] \\
                \st \quad                                                      & \ratesubmat\transpose\bm 1 \leq  \maxratesubvec,                           \\
                                                                               & \bm 0 \leq \ratesubmat \leq \bar{\capmat}.
            \end{align}
        \end{subequations}

        \blt[Approximate solution]
        \begin{bullets}
            \blt[relaxed problem]Due to the  combinatorial nature of the problem, it is convenient to relax the objective function as follows:
            \begin{subequations}
                \begin{align}
                    \minimize_{\ratesubmat\in\rfield^{\usernum\times\absnum}}\quad & \sum_{\userind}\max(\minrate - \ratervec_\userind\transpose\bm 1, 0)                                             \\
                    \st \quad                                                      & \ratesubmat\transpose\bm 1 \leq  \maxratesubvec,                    ~ \bm 0 \leq \ratesubmat \leq \bar{\capmat}.
                \end{align}
            \end{subequations}
            \blt[equivalently]This problem can be clearly rewritten as
            \begin{subequations}
                \label{eq:problemwithsumofwhys}
                \begin{align}
                    \minimize_{\ratesubmat\in\rfield^{\usernum\times\absnum}, \yvec\in\rfield^{\usernum}}\quad & \sum_{\userind}\ys_\userind                                                                                     \\
                    \st \quad                                                                                  & \ratesubmat\transpose\bm 1 \leq  \maxratesubvec,                    ~ \bm 0 \leq \ratesubmat \leq \bar{\capmat}, \\
                                                                                                               & \ys_\userind=\max(\minrate - \ratervec_\userind\transpose\bm 1, 0),
                \end{align}
            \end{subequations}
            where $\ys_\userind$ is the $\userind$-th entry of the auxiliary variable $\yvec$. Imposing non-negativity of $\yvec$, one can simplify \eqref{eq:problemwithsumofwhys}     as
            \begin{subequations}
                \begin{align}
                    \minimize_{\ratesubmat\in\rfield^{\usernum\times\absnum}, \yvec\in\rfield^{\usernum}}\quad & \bm 1\transpose\yvec                                                                                                                     \\
                    \st \quad                                                                                  & \ratesubmat\transpose\bm 1 \leq  \maxratesubvec,~\bm 0 \leq \ratesubmat \leq \bar{\capmat},                                               \\
                                                                                                               & \bm 0 \leq \yvec,                                                                           ~ \minrate\bm 1 - \ratesubmat\bm 1\leq \yvec,
                \end{align}
            \end{subequations}
            which is a linear program and, therefore, can be solved by standard methods. Again, applying the reweighting technique  to this problem would generally result in a greater number of served users; cf.~\eqref{eq:mainproblem}.
        \end{bullets}
    \end{bullets}
\end{changes}

\section{Numerical Experiments}
\label{sec:experiments}



\begin{figure}
    \centering
    \begin{subfigure}[b]{.45\linewidth}
        \centering
        \includegraphics[width=\linewidth,trim={2cm 0 2 2}]{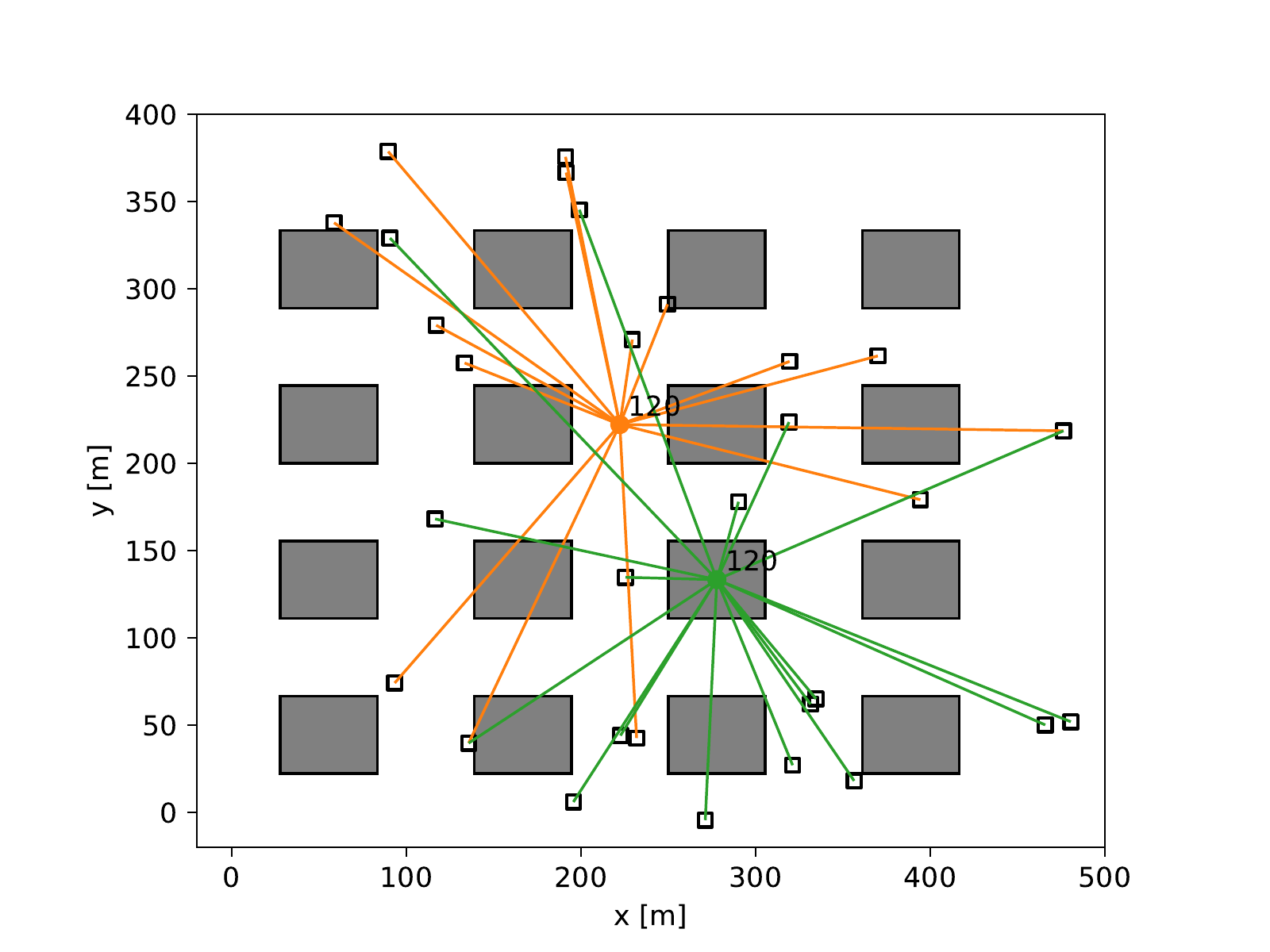}
        \caption{Building absorption: 0.25 dB/m.}
    \end{subfigure}
    \hfill
    \begin{subfigure}[b]{.45\linewidth}
        \centering
        \includegraphics[width=\linewidth,trim={2cm 0 2 2}]{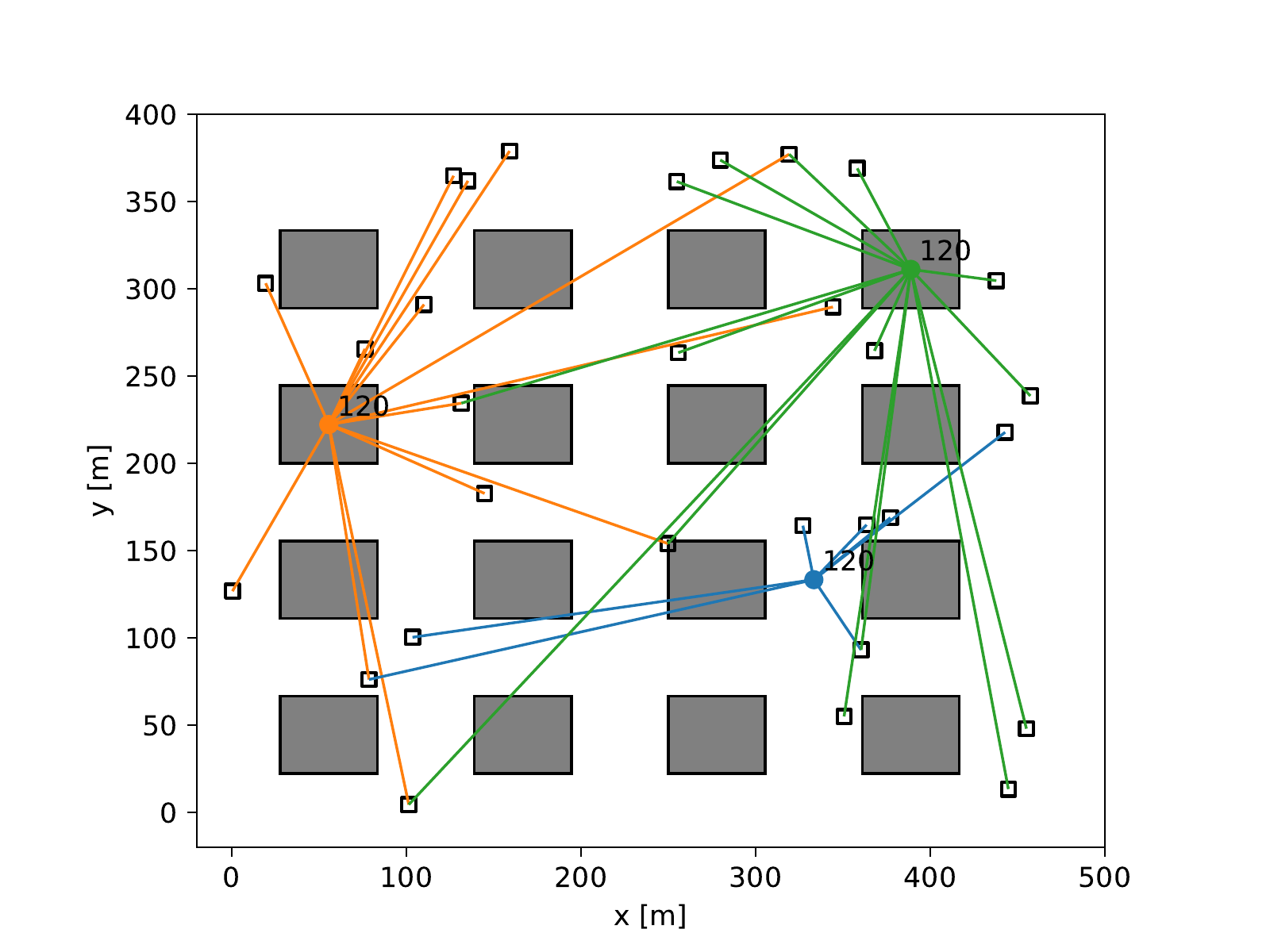}
        \caption{Building absorption: 1.25 dB/m.}
    \end{subfigure}
    \hfill
    \begin{subfigure}[b]{.45\linewidth}
        \centering
        \includegraphics[width=\linewidth,trim={2cm 0 2 2}]{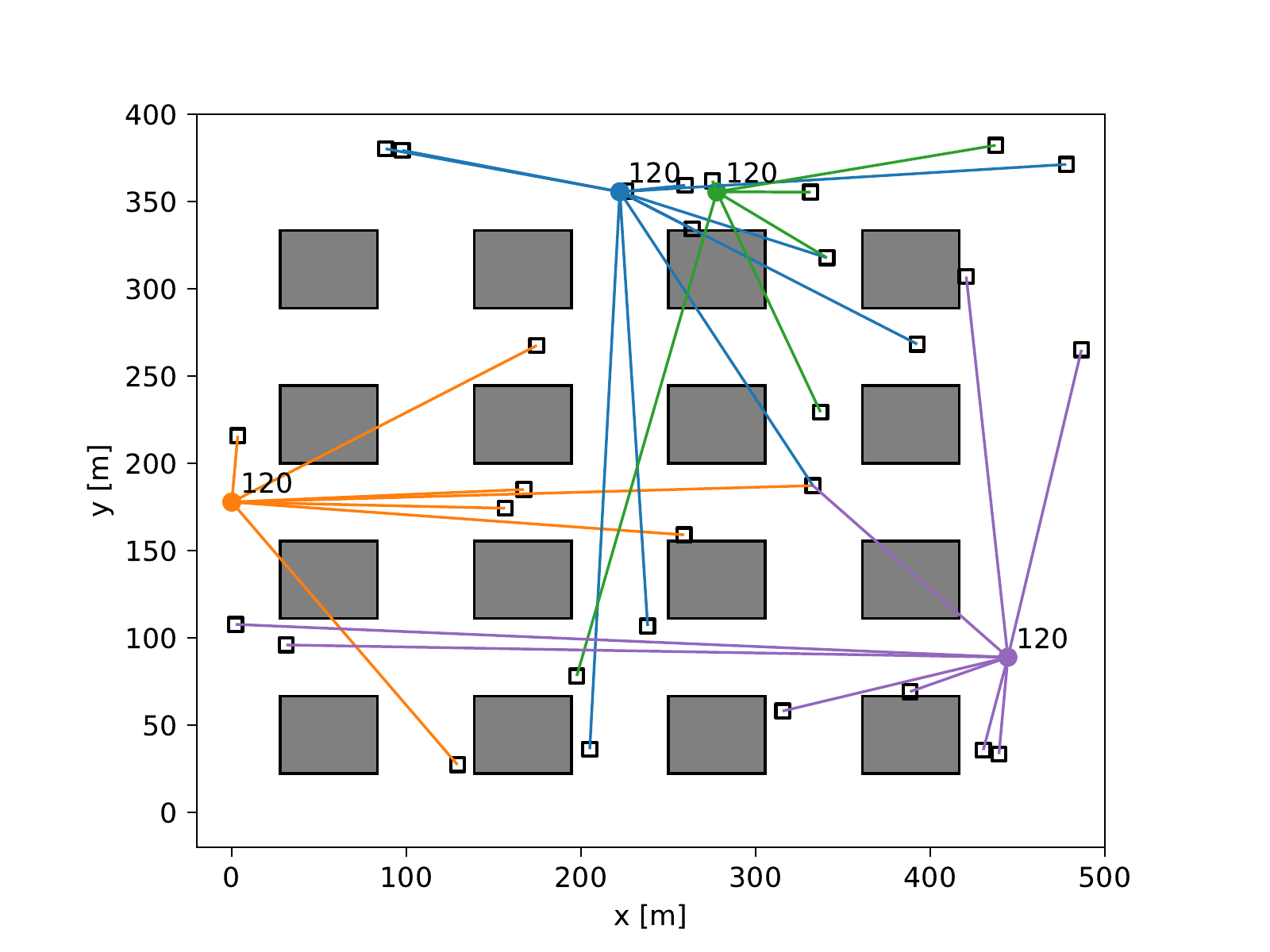}
        \caption{Building absorption: 1.75 dB/m.}
    \end{subfigure}
    \hfill
    \begin{subfigure}[b]{.45\linewidth}
        \centering
        \includegraphics[width=\linewidth,trim={2cm 0 2 2}]{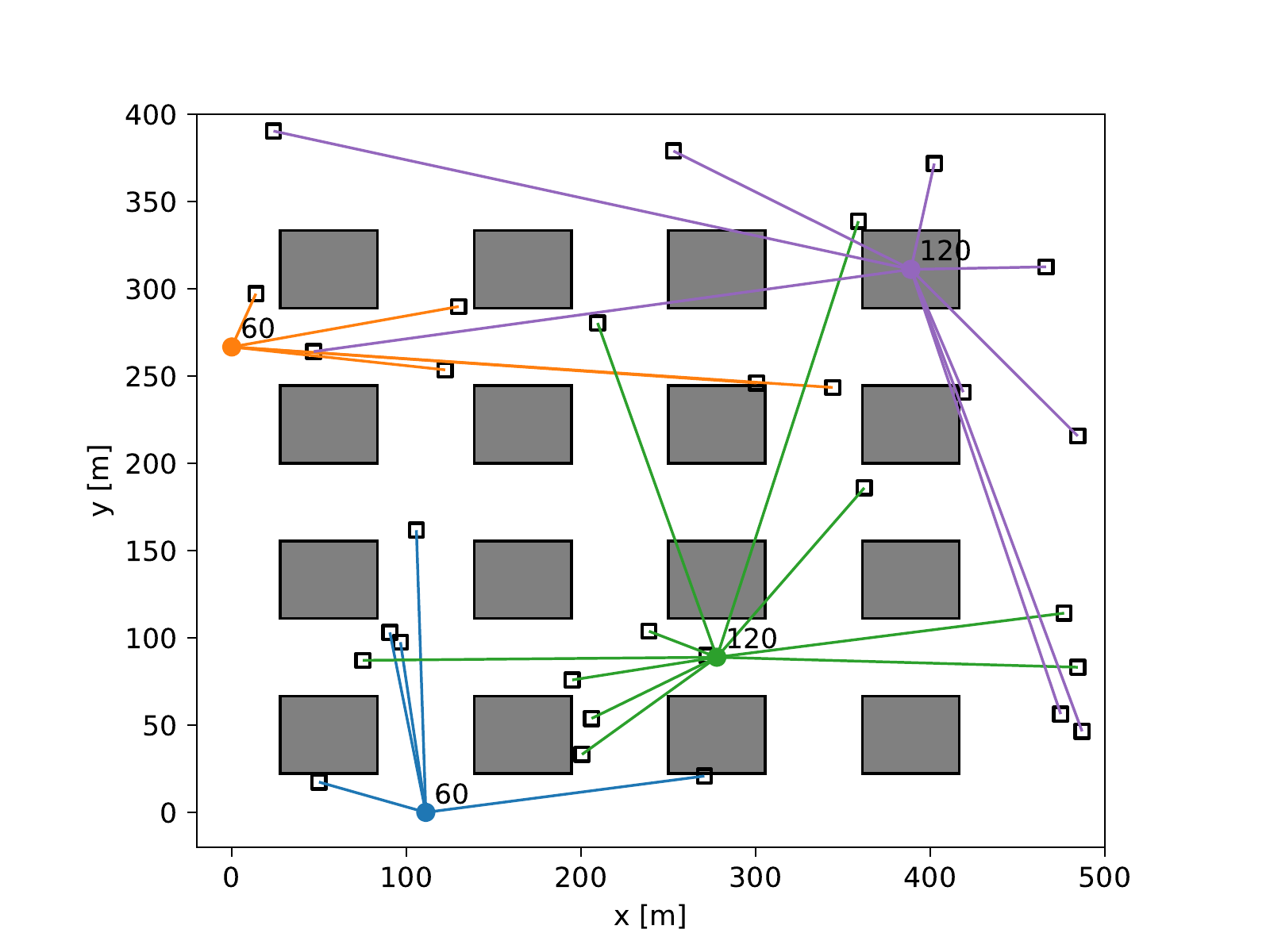}
        \caption{Building absorption: 2.5 dB/m.}
    \end{subfigure}
    \caption{
        \begin{changes}
            Plan view of four placements provided by GSPA. Gray rectangles represent buildings, which are 63 m high. There are 30 GTs on the ground, denoted by black squares. Circles of different colors represent the locations of the ABSs. Their  connections with the GTs are indicated by line segments of the same color. The heights of the ABSs are annotated right next to their markers ($\minrate$~=~1~Mbps, $\maxrate{}~=~100~\text{Mbps}$).
        \end{changes}
    }
    \label{fig:connections}
\end{figure}

\begin{figure}[t]
    \centering
    \centering
    \captionsetup{width=.9\linewidth}
    \includegraphics[width=1\linewidth]{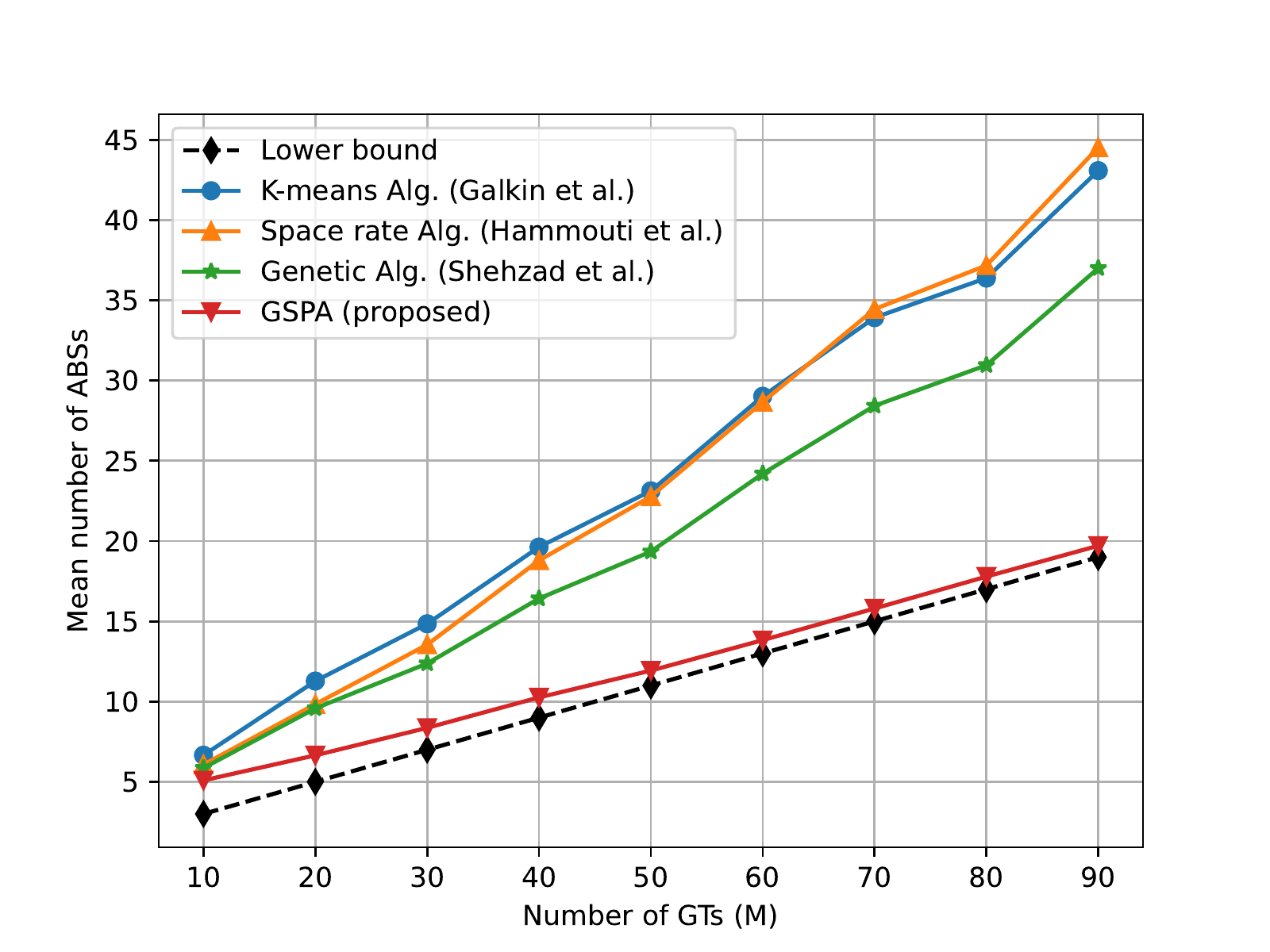}
    \captionof{figure}{Mean number of ABSs vs. number of GTs
        ($\minrate$ = 20 Mbps, $\maxrate{}$ = 100
        Mbps).}
    \label{fig:vs_num_users}
\end{figure}

\begin{figure}[t]
    \centering
    \captionsetup{width=.9\linewidth}
    \includegraphics[width=1\linewidth]{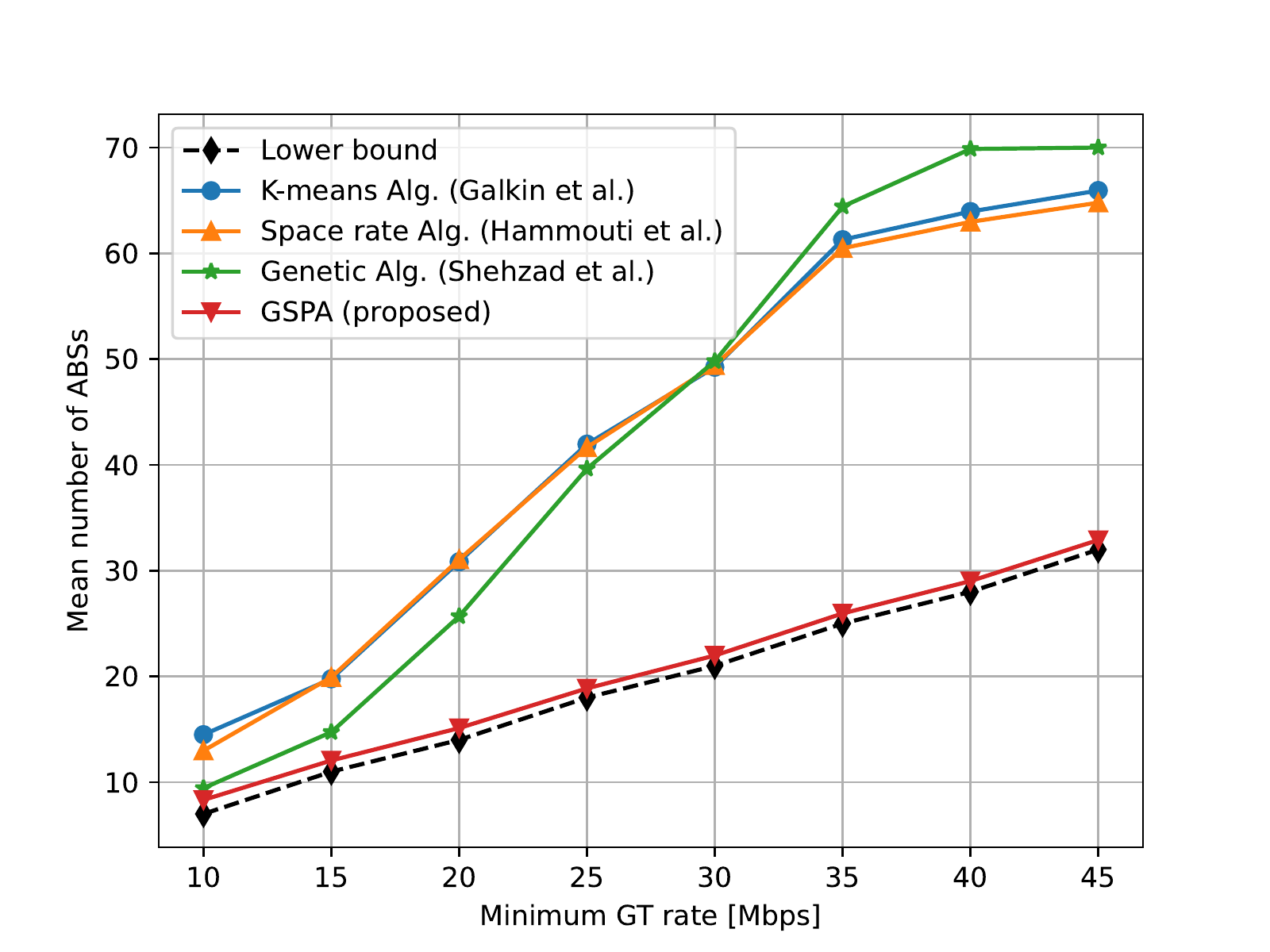}
    \captionof{figure}{Mean number of ABSs vs. $\minrate$ ($\maxrate{} = 100$ Mbps).}
    \label{fig:vs_minRate}
\end{figure}

\begin{bullets}
    \blt[Overview]This section empirically validates the performance of
    the proposed algorithm by means of numerical experiments with
    channel data generated using the tomographic model
    (Sec.~\ref{sec:tomographic_exp}), ray-tracing software
    (Sec.~\ref{sec:ray_exp}), and the model in \cite{alhourani2014lap} (Sec.~\ref{sec:alhourani}). All details of the simulations, such as
    the locations and size of the buildings, cannot be presented here due to
    space limitations, but the code and data necessary to reproduce the
    experiments is available at
    \url{https://github.com/uiano/ABS_placement_via_propagation_maps}.

    \blt[Data generation]
    \begin{bullets}
        \blt[Channel generation \ra following subsections]The channel gain,
        obtained from the aforementioned models when the carrier frequency is $2.4$
        GHz, is
        substituted  into  \eqref{eq:capfun} with
        \blt[Capacity matrix generation]
        \begin{bullets}%
            \blt[] $\bandwidth=20$ MHz,
            \blt[] $\txpow=20$ dBm, and
            \blt[]  $\noisepow = -96$ dBm
        \end{bullets}%
        to form the capacity matrix $\capmat$.
        \blt[BH capacity]For simplicity, the backhaul capacity is set to
        a common value $\maxrate{1}=\ldots=\maxrate{\gridptnum}=\maxrate{}$.
    \end{bullets}

    \blt[Compared algorithms \ra list and say parameters] The proposed
    placement algorithm is compared with three benchmarks:
    \begin{bullets}%
        \blt[Kmeans] i) the \textit{K-means placement algorithm} by Galkin et al. \cite{galkin2016deployment},
        \blt[SpaceRateKmeans]ii) the iterative \textit{space rate K-means placement algorithm} by Hammouti et al. \cite{hammouti2019mechanism}, and
        \blt[Genetic] iii) the \textit{genetic placement algorithm} by Shehzad et al. \cite{shehzad2021backhaul} with
        \begin{bullets}%
            \blt[num-sol-per-gen]50 solutions per generation.
        \end{bullets}%
        \blt[self-implemented]Since we did not manage to obtain the code
        used by the authors of these works, we implemented their
        algorithms ourselves. The resulting implementations are
        available in the repository mentioned earlier.
        \blt[fly grid] The first two algorithms are unable to enforce
        no-fly zones, which results in ABSs placed inside buildings. To
        avoid this behavior, the ABS locations provided by these
        algorithms are projected onto the same flight grid as the rest of
        algorithms.
        \blt[Parameters of the proposed algorithm] The proposed GSPA
        algorithm utilizes, unless otherwise stated,
        \begin{bullets}%
            \blt[]  a step size of $\admmstep=  10^{-7}$ and
            \blt[stopping criterion] stopping criterion parameters
            $\admmstop_{\rm{abs}}=\admmstop_{\rm{rel}}=10^{-4}$; cf.  Sec. \ref{sec:stopping} of the supplementary material.
        \end{bullets}%
    \end{bullets}%

    \blt[Performance metrics]To quantify performance, the number of ABSs
    required by each algorithm to guarantee a rate $\minrate$ for every
    GT is considered as a performance metric. This metric is averaged
    using Monte Carlo simulation across realizations of the GT
    locations. As a reference, figures will also include a lower bound
    on this metric; cf. Appendix~\ref{app:bound}.

\end{bullets}

\subsection{Experiments with the Tomographic Model}
\label{sec:tomographic_exp}

\begin{bullets}%
    %
    \blt[Channel generation]In the experiments of this section, the
    channel gain is generated using Algorithm~\ref{algo:tomo} in an
    environment like the one in Fig.~\ref{fig:urban}.  The SLF takes a
    constant value, termed \emph{building absorption}, inside the
    buildings and 0 outside.  Unless otherwise stated, the simulation
    parameters in this section are as follows.  Region size: 500 m
    $\times$ 400 m $\times$ 150 m; dimensions of the SLF grid: 50
    $\times$ 40 $\times$ 15; dimensions of the flight grid: 9 $\times$ 9
    $\times$ 5; minimum flight height: 50 m; height of the buildings: 63
    m; building absorption: 1
    dB/m; number of GTs: $\usernum= 70$.


    \blt[Experiments]
    \begin{bullets}%
        \begin{changes}
            \blt[Exp - realizations] Before characterizing performance via Monte Carlo simulations, the ABS placements corresponding to four realizations of the GT locations are shown in Fig.~\ref{fig:connections} for four values of the building absorption. As expected, the number of ABSs required to serve the GTs increases with the building absorption, as the propagation environment becomes less favorable. For large building absorptions, most links have line of sight, which also agrees with intuition. In all cases, most GTs are served by just one ABS, which
            validates the  algorithm from Sec.~\ref{sec:minnumconnections}. Some of the GTs are however served by multiple ABSs, especially in the case of GTs with poor channels to the ABSs. The number of GTs served by multiple ABSs can be made even smaller if the reweighting technique mentioned at the end of Sec.~\ref{sec:minnumconnections} is applied.

        \end{changes}

        \blt[Exp. - vs num-users]Fig. \ref{fig:vs_num_users}
        investigates the influence of the number of GTs ($\usernum$) on the
        performance of the compared algorithms.
        \begin{bullets}%
            %
            %
            \blt[Interpretation]Several observations are in
            order.
            \begin{bullets}
                \blt[proportional]First, the mean number of ABSs
                is seen to increase roughly proportionally to
                $\usernum$ for all the algorithms.
                \blt[outperform]Second, the proposed GSPA algorithm not
                only yields a lower mean number of ABSs than the
                competing algorithms, but its slope is smaller,
                which means that the margin by which GSPA
                outperforms the benchmarks increases with
                $\usernum$.
                \blt[lowerbound]Third, GSPA asymptotically
                approaches the lower bound
                (cf. Appendix~\ref{app:bound}), which means that
                its efficiency, quantified as the number of GTs per
                ABS, increases with $\usernum$. In contrast, the
                opposite is true for the other algorithms.
                %
            \end{bullets}
        \end{bullets}

        \blt[Exp. - vs min-user-rate]To investigate the influence
        of the GT requirements on performance,
        Fig.~\ref{fig:vs_minRate} depicts the mean number of ABSs
        vs.   $\minrate$   when
        \begin{bullets}
            \blt[How it is run] $\maxrate{} = 100$ Mbps.
            \blt[figure]%
            %
            \blt[Interpretation]
            \begin{bullets}%
                \blt[proportional]As in Fig.~\ref{fig:vs_num_users},
                the mean number of ABSs seems to increase roughly
                linearly. However, in Fig.~\ref{fig:vs_minRate} a
                saturation phenomenon arises: the mean number of
                ABSs cannot be greater than the number of GTs
                $\usernum=70$ since one ABS placed approximately
                above each GT suffices to serve all GTs.
                \blt[genetic] It is also observed that the
                performance of the genetic placement algorithm degrades faster
                than the rest of algorithms for high $\minrate$. The
                reason may be that this algorithm essentially tests
                multiple placements and the number of placements
                increases drastically with the number of ABSs, which
                is larger when $\minrate$ is larger.

            \end{bullets}
        \end{bullets}

        \begin{figure}[t]
            \centering \centering \captionsetup{width=.9\linewidth}
            \includegraphics[width=1.05\linewidth]{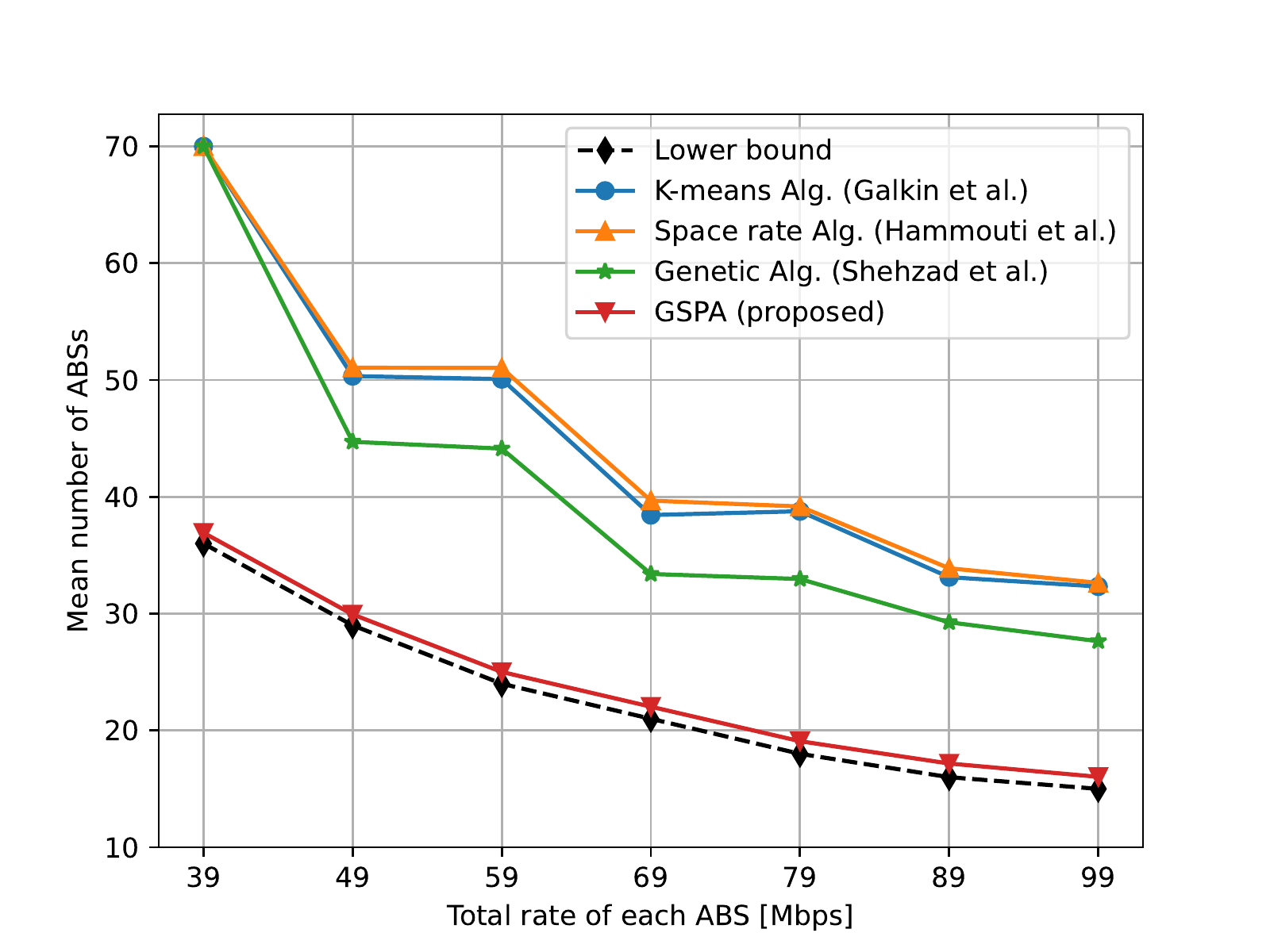}
            \captionof{figure}{Mean number of ABSs vs. the capacity
                $\maxrate{}$ of the backhaul link ($\minrate$ = 20
                Mbps).  }
            \label{fig:vs_maxRate}
        \end{figure}%

        \blt[Exp. - vs max-uav-total-rate]Fig.~\ref{fig:vs_maxRate}
        plots the mean number of ABSs vs. the backhaul capacity $\maxrate{}$
        \begin{bullets}
            \blt[How it is run]when $\minrate$ = 20 Mbps.
            \blt[Figure]
            %
            \blt[Interpretation]
            \begin{bullets}
                \blt[multiples]The values of $\maxrate{}$ on the
                horizontal axis are selected so that each $\maxrate{}$
                is not an integer multiple of $\minrate$. This gives
                rise to a ``staircase'' behavior for the benchmarks,
                which do not exploit the fact that a GT can be
                served by multiple ABSs. In contrast, GSPA exploits
                this fact, as corroborated by the smoothness of its
                curve.
                \blt[twice better]It is also seen that the mean number
                of ABSs required by GSPA is roughly half the one of the best
                competing alternative.
            \end{bullets}
        \end{bullets}

        \begin{changes}
            \blt[Exp. impact of radio map mismatch]%
            \begin{bullets}%
                \blt[motivation \ra mismatch] The last experiment of this section investigates the impact of radio map mismatches, as described in Sec.~\ref{sec:masnumservedusers}. To this end, mismatched radio map estimates are generated by assuming a wrong value of the building absorption. This provides a clean comparison, as introducing other deviations, such as modifying the building locations or dimensions, would change the flight grid $\grid$ and therefore introduce complicated effects that would obscure the phenomenon under study. For the same reason, the backhaul capacity is assumed infinite.

                \blt[figure] Fig.~\ref{fig:fractvsabsorption}  shows the mean number of ABSs and the mean fraction of served users vs. the building absorption used to create the radio map estimate used by GSPA. The true building absorption is 1.25 dB/m, which means that the radio map estimates are optimistic to the left of 1.25 and  pessimistic to the right of 1.25. The algorithm from Sec.~\ref{sec:masnumservedusers} is used to maximize the number of served users. It is observed that the fraction of served users degrades as the building absorption becomes smaller. However, this degradation is more pronounced when the $\minrate$ is small. The reason is that a larger $\minrate$ gives rise to a larger number of ABSs, which in turn means that GTs are more likely to have a good channel to at least one ABS. In the limit case where $\minrate$ is very small and the building absorption is 0, GSPA places a single ABS to serve all users. To sum up, GSPA is affected by radio map mismatches, but their impact, which is more pronounced when $\absnum$ is small, is mitigated by the algorithm in Sec.~\ref{sec:masnumservedusers}.

                \blt[interpretation]

            \end{bullets}%

        \end{changes}

    \end{bullets}

\end{bullets}%

\begin{figure}[t]
    \centering
    \includegraphics[width=\linewidth]{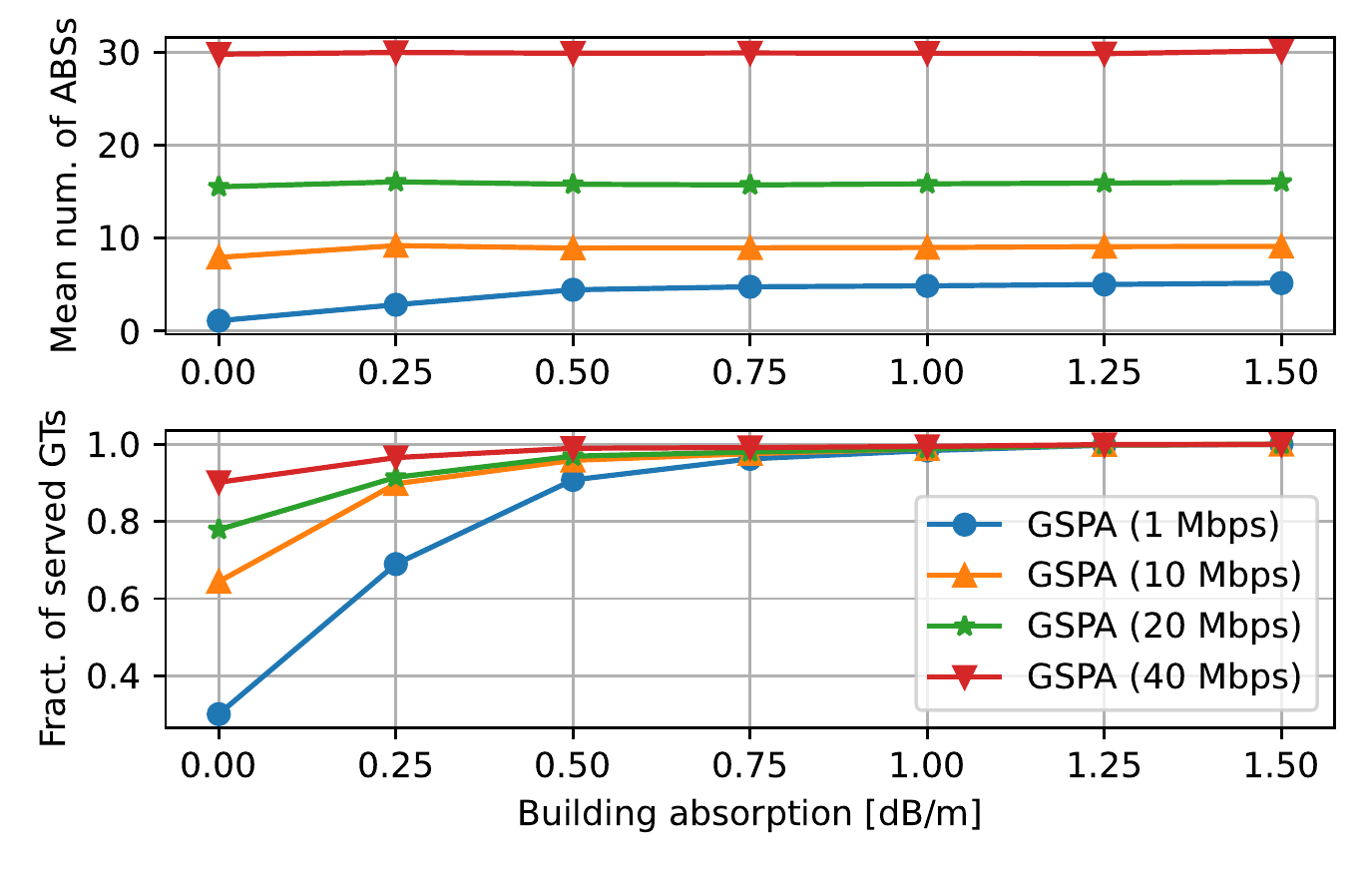}
    \caption{\begin{changes}
            Impact of the mismatch between the true radio map and the radio map estimate on the performance of the proposed algorithm. The x-axis indicates the value of the building
            absorption assumed by the algorithm, whereas the true building absorption is 1.25 dB/m. The rate in the legend indicates $\minrate$.
        \end{changes}
    }
    \label{fig:fractvsabsorption}
\end{figure}

Additional experiments with the tomographic channel model are
presented in Sec.~\ref{sec:additionaltomo} of the supplementary material.

\subsection{Experiments with the Ray-Tracing Model}
\label{sec:ray_exp}

\begin{bullets}
    \blt[Overview] This section corroborates the main findings of
    Sec.~\ref{sec:tomographic_exp} when the channel is obtained via the
    ray-tracing model, which is more accurate than the tomographic
    channel model for higher carrier frequencies.
    \blt[Channel generation]
    \begin{bullets}%
        \blt[Channel model] To this end, a data set was generated using
        the X3D ray-tracing model of the software Wireless Insite with six reflections
        and one diffraction in  a 400 $\times$ 600~m area of the city
        of Ottawa. The data set is also published in our
        repository.
        \blt[Fly grids and users]The channel was computed between all
        points of the flight grid and all points of a GT grid. The former
        comprises 35 points at each of the heights of 40, 60, and 80
        m. The latter is a 2D regular grid of 2501 GT locations at
        a height of 2 m spaced uniformly with a distance of 10 m on each
        axis. Points inside the buildings are removed.
        \blt[GT locations] At each Monte Carlo realization, the GT
        locations are generated by drawing $\usernum$ points uniformly at
        random without replacement from the GT grid.
    \end{bullets}

    \blt[Experiments]
    \begin{bullets}

        \begin{figure}[t]
            \centering
            \centering
            \captionsetup{width=.9\linewidth}
            \includegraphics[width=1\linewidth]{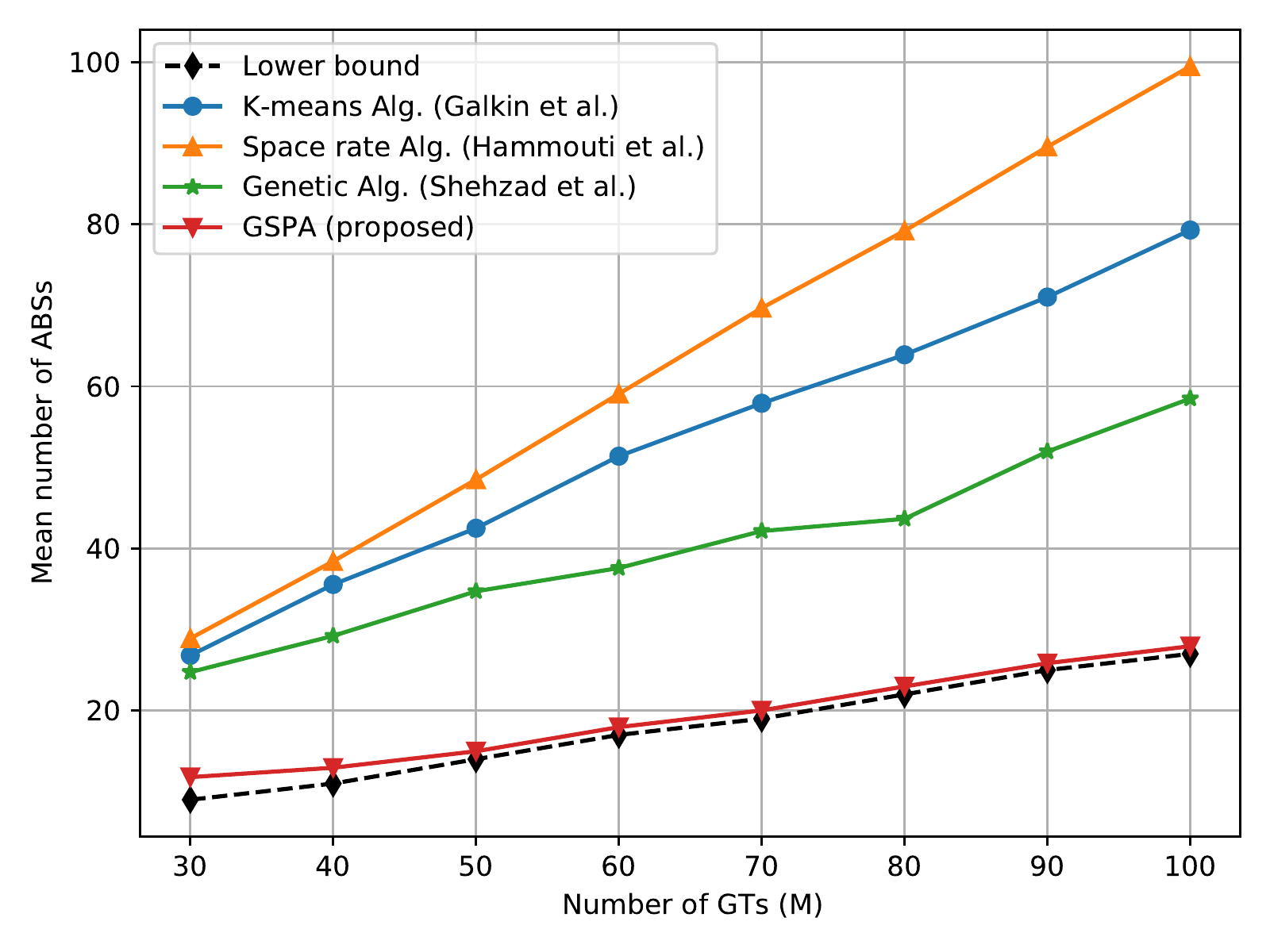}
            \captionof{figure}{Mean number of ABSs vs. number of
                GTs  ($\minrate = 20$ Mbps, $\maxrate{}=
                    74$ Mbps).}
            \label{fig:vs_num_users_ray_tracing}
        \end{figure}%

        \blt[no. uav vs. no. users]
        Figs.~\ref{fig:vs_num_users_ray_tracing}
        and~\ref{fig:vs_minRate_maxRate_100_ray_tracing} are the
        counterparts of Figs.~\ref{fig:vs_num_users} and
        \ref{fig:vs_minRate} for ray-tracing channel data.
        \begin{bullets}%
            \blt[Interpretation]%
            \begin{bullets}%
                \blt[differences] It is observed that the mean number of ABSs
                is also an approximately linear function of $\usernum$ for all
                algorithms and that GSPA roughly attains the lower
                bound. However, here the differences among benchmarks are
                greater. The fact that the K-means algorithm outperforms the
                space rate algorithm suggests that the channel changes rapidly
                with respect to the GT location, since the latter algorithm
                relies on clustering vectors of channel gains. Despite this
                fact, GSPA performs almost optimally.
                %
                %
            \end{bullets}%
        \end{bullets}%

        \blt[no. uav vs. min. GT rate]
        Fig.~\ref{fig:vs_minRate_GroupSparsePlacer_ray_tracing} 
        investigates the impact of the rate constraints. Observe that the tightness of the
        bounds in Fig.~\ref{fig:vs_minRate_GroupSparsePlacer_ray_tracing}
        increases (i) for larger $\minrate$ and (ii) for smaller
        $\maxrate{}$. Those are precisely the situations where the
        backhaul limitations become stricter. Indeed, it can be easily
        seen that the bound in Appendix~\ref{app:bound} can be attained
        with equality when the entries of $\capmat$ approach infinity.


        %

        \begin{figure}[h!]
            \centering
            \captionsetup{width=.9\linewidth}
            \includegraphics[width=1\linewidth]{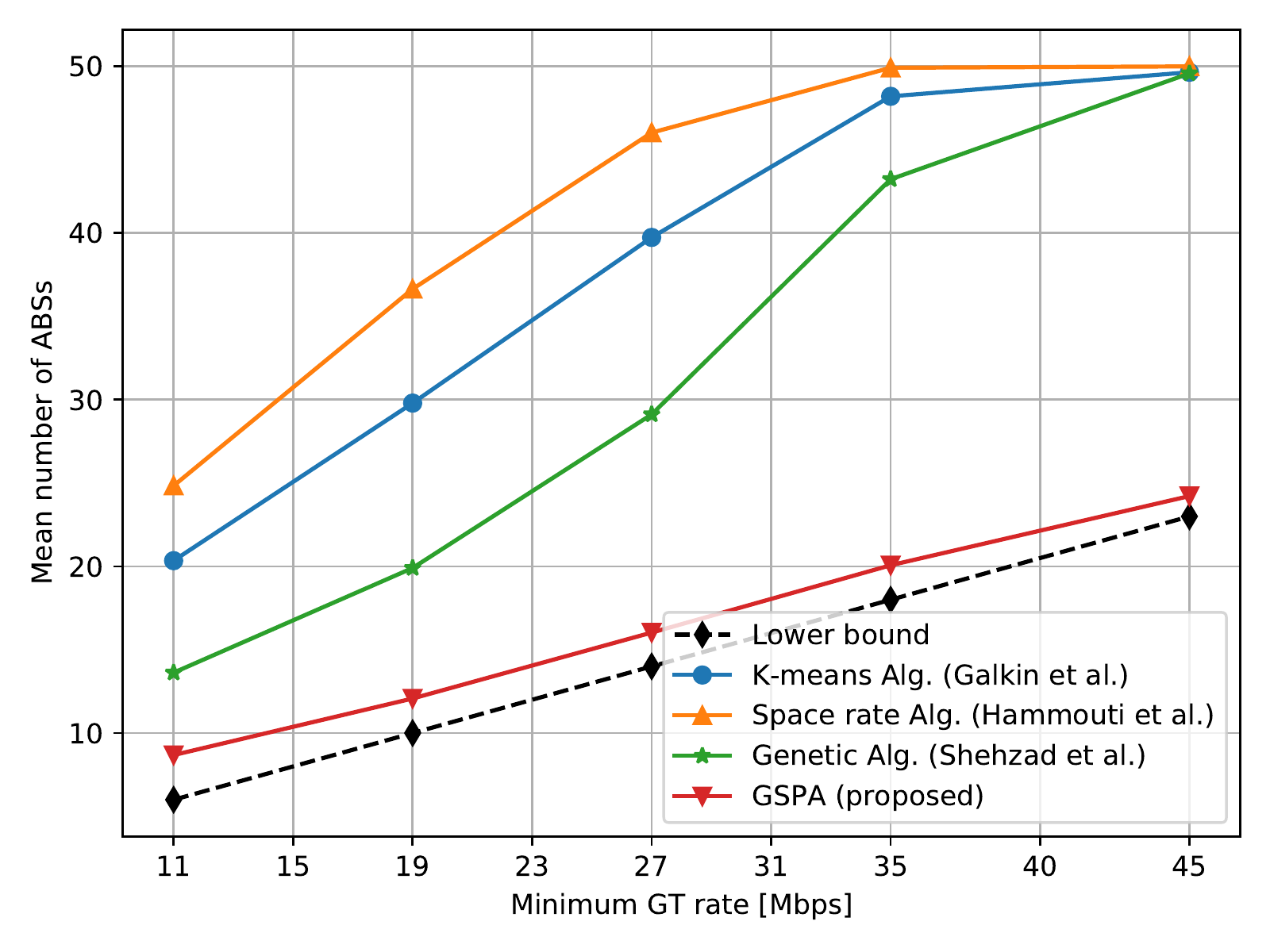}
            \captionof{figure}{Mean number of ABSs vs. minimum
                GT rate $\minrate$ ($\usernum=50$ GTs, $\maxrate{}$ = 100 Mbps).}
            \label{fig:vs_minRate_maxRate_100_ray_tracing}
        \end{figure}

        \begin{figure}[h!]
            \centering
            \captionsetup{width=.9\linewidth}
            \includegraphics[width=1\linewidth]{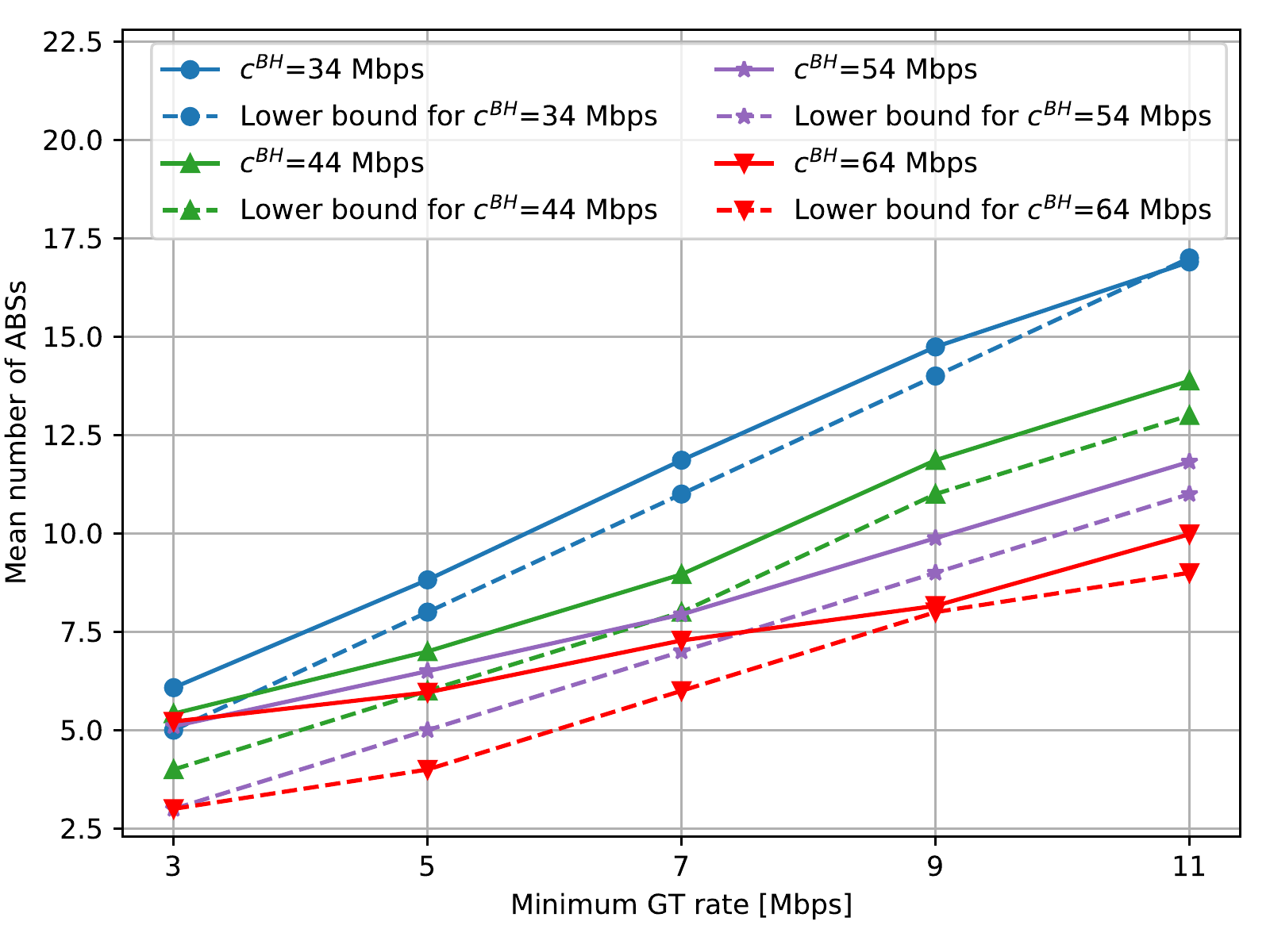}
            \captionof{figure}{Mean number of ABSs vs. minimum GT
                rate $\minrate$ for GSPA ($\usernum=50$ GTs). }
            \label{fig:vs_minRate_GroupSparsePlacer_ray_tracing}
        \end{figure}%

        %
        %



        \blt[no. uav vs. max. uav rate]Finally, the experiment in
        Fig.~\ref{fig:vs_maxRate} is performed with ray-tracing
        channel data in Fig.~\ref{fig:vs_maxRate_ray_tracing}.
        \begin{bullets}
            \blt[How it is run]
            \blt[Interpretation]
            \begin{bullets}%
                \blt[Outperform] The proposed algorithm still outperforms the
                other three benchmarks by a wide margin, requiring roughly $50\%$ of
                ABSs.
                \blt[saturation] However, relative to
                Fig.~\ref{fig:vs_maxRate}, one can observe that the benchmark
                algorithms saturate for large $\maxrate{}$, which indicates
                that $\capmat$ imposes a more stringent constraint than
                $\maxrate{}$ in the scenario of
                Fig.~\ref{fig:vs_maxRate_ray_tracing}. When it comes to GSPA,
                the greater tightness of the bound for small $\maxrate{}$ is a
                manifestation of the same effect, as discussed earlier.

            \end{bullets}
        \end{bullets}


    \end{bullets}
\end{bullets}

\subsection{Experiment with the model by Al-Hourani et al.}
\label{sec:alhourani}

\begin{bullets}%

    \blt[Exp. - vs Al-Hourani] Fig.~\ref{fig:vs_alhourani} compares the performance metric obtained using the tomographic model and the model in \cite{alhourani2014lap} whose parameters $(a,b,\eta_\text{LOS},\eta_\text{NLOS})=(12.08, 0.11, 2.3, 34)$ are taken from the  dense urban and high-rise urban environments. Overall, the number of ABSs also increases roughly linearly with the number of GTs when the model in \cite{alhourani2014lap} is used. Note that a placement may be valid according to one channel model and not according to the other. In this specific environment, fewer ABSs are required if the tomographic model is adopted, but the results must be compared with care since none of the channels constitutes the ground truth.

\end{bullets}%

\section{Conclusions}
\label{sec:conclusions}

Whereas virtually all existing algorithms for ABS placement assume that the
channel gain depends only on the length and (possibly)
elevation of each link, this paper presents a scheme that can
accommodate an arbitrary dependence of the gain on the position
of the ABSs and GTs. This enables the utilization of radio maps
for multiple-ABS placement. The proposed algorithm determines a set of
ABS locations that approximately minimizes the number of ABSs
required to guarantee a minimum rate to all GTs. Relative to
existing schemes, the proposed algorithm has low
complexity, accounts for limited backhaul capacity, and can
accommodate flight restrictions such as no-fly zones or
airspace occupied by buildings. A solver whose complexity is
linear in the number of GTs was derived based on the
alternating-direction method of multipliers and the problem of
evaluating tomographic integrals was revisited and extended to
air-to-ground channels. An extensive set of simulations
demonstrate that the proposed GSPA algorithm outperforms
competing algorithms by a wide margin with three types of channels. Remarkably, it was observed in the
numerical experiments that the proposed algorithm is the only
one among the compared schemes whose efficiency, measured in
terms of the number of GTs served per ABS, increases with the number
of GTs. This fundamental distinction renders GSPA especially
suitable for scenarios with a large number of GTs.

Future directions include approaches for tracking air-to-ground
propagation maps, possibly based on online kernel
methods~\cite{romero2015onlinesemiparametric,romero2017spectrummaps},
and algorithms that can adapt to changes in the GT locations.

\section{Acknowledgements}

The authors would like to thank Prof. Geert Leus for insightful
discussions.

\section*{Appendices}

\renewcommand{\thesubsection}{\Alph{subsection}}

\subsection{Lower  Bound for  the Number of ABSs}
\label{app:bound}

\begin{bullets}%
    \blt[Overview]This appendix presents a lower bound for the number of
    ABSs and, therefore, also for the \emph{mean} number of ABSs, which is
    the performance metric adopted in Sec.~\ref{sec:experiments}. This
    bound constitutes a fundamental limit for the problem of ABS placement
    and, hence, it applies regardless of the adopted algorithm.

    \blt[lower bound] Let $\absnum$ denote the smallest number of
    ABSs required to serve all $\usernum$ GTs with rate at least
    $\minrate$.  This means that the total backhaul rate available to
    all ABSs together, which is not greater than $\absnum
        \max_\gridptind \maxrate{\gridptind}$, cannot be less than the
    total rate $\usernum\minrate$ demanded by the GTs. It follows
    that $\absnum \max_\gridptind \maxrate{\gridptind}\geq
        \usernum\minrate$ and, therefore, $\absnum \geq\lceil
        \usernum\minrate/\max_\gridptind \maxrate{\gridptind}\rceil$,
    where $\lceil z\rceil$ denotes the smallest integer greater than
    or equal to $z$.

    \begin{figure}[t]
        \centering
        \captionsetup{width=.9\linewidth}
        \includegraphics[width=1\linewidth]{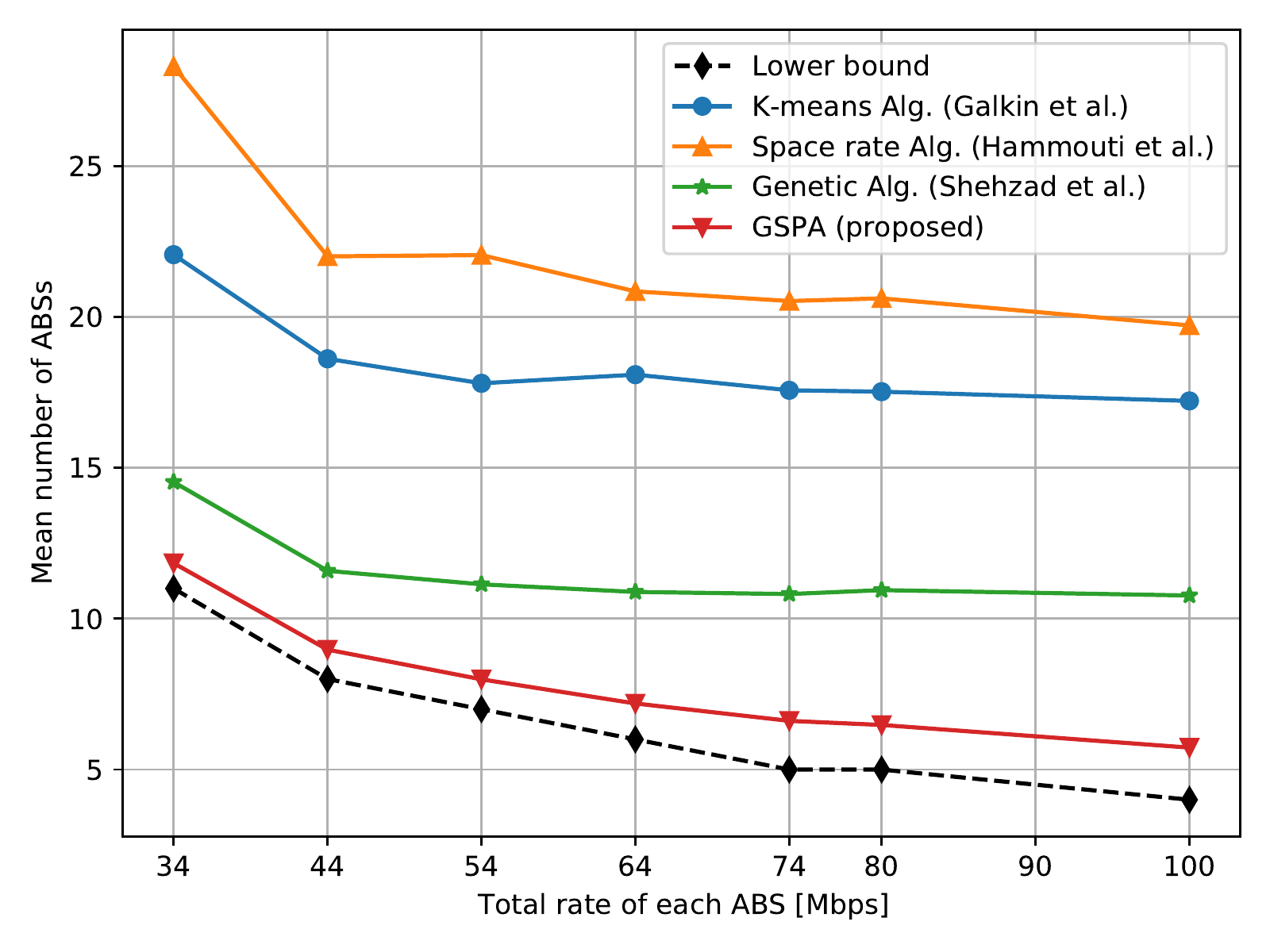}
        \captionof{figure}{Mean number of ABSs vs. backhaul
            link capacity $\maxrate{}$ ($\minrate = 7$ Mbps, $\usernum=50$ GTs).}
        \label{fig:vs_maxRate_ray_tracing}
    \end{figure}

    \begin{figure}[t]
        \centering
        \includegraphics[width=.5\textwidth]{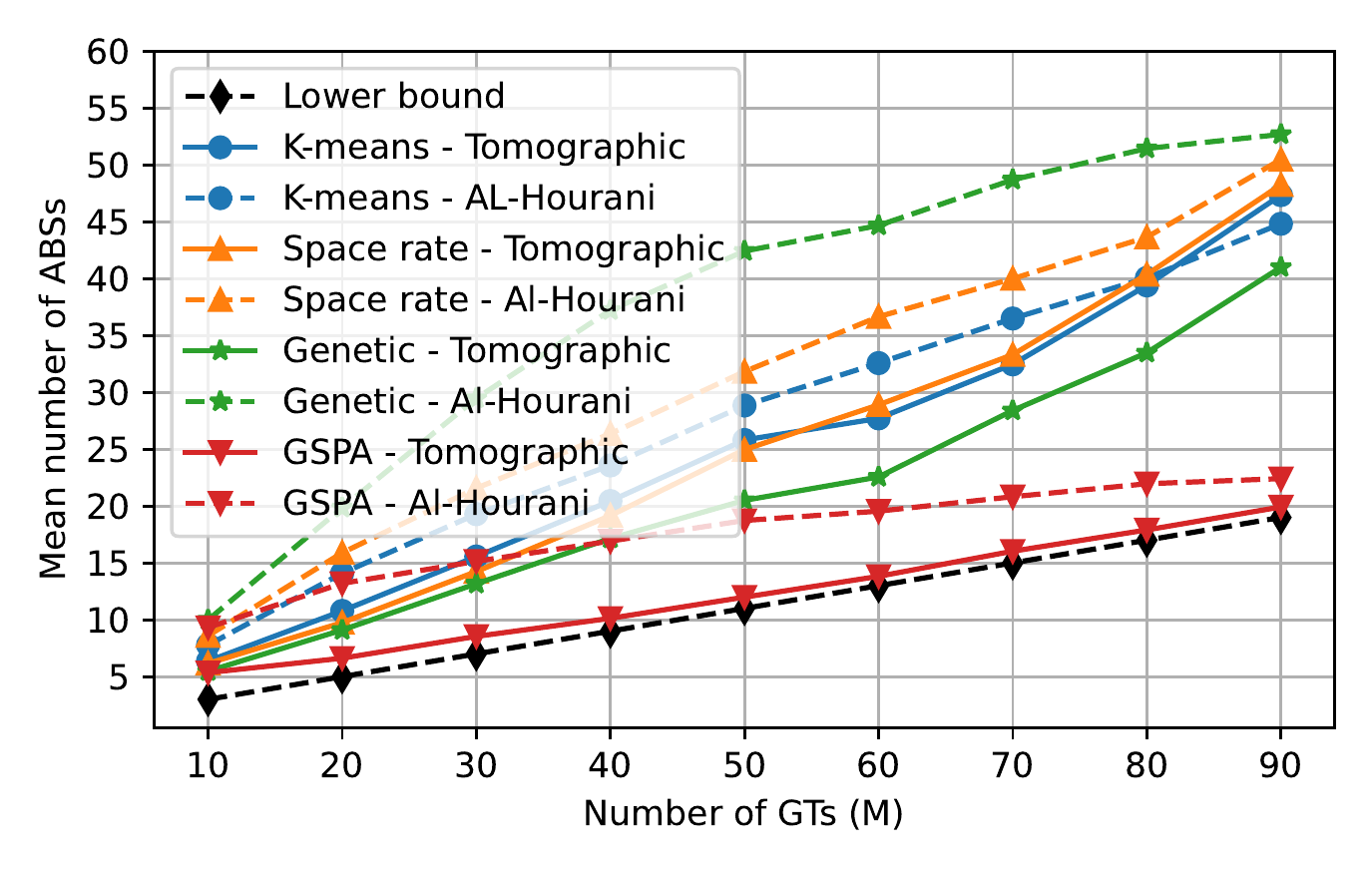}
        \caption{Mean number of ABSs vs. number of GTs ($\minrate= 20$ Mbps, $\maxrate{}= 99$
            Mbps).
        }
        \label{fig:vs_alhourani}
    \end{figure}







\end{bullets}


\printmybibliography






~\\[.2cm]
\setcounter{section}{0}
\renewcommand*{\thesection}{S\arabic{section}}

\textbf{Supplementary Material}

\section{Derivation of Algorithm~\ref{algo:tomo}}
\label{sec:algotomoderivation}

Algorithm~\ref{algo:tomo}, which can be classified as a
parametric, floating point, and zeroth-order algorithm
according to the terminology of
\cite[Sec. I-B-1]{mitchell1990comparison}, is our approach (yet
others are possible) to approximate the tomographic integral by
computing exactly the integral of a piecewise constant
approximation of the SLF.
\begin{bullets}%
    \blt[parameterization] The idea is to parameterize the line
    segment between $\loc_1$ and $\loc_2$ as
    $\loc(t)=\loc_1 + t(\loc_2-\loc_1)$, where $t\in[0,1]$, and
    identify the values $t_1<t_2<\ldots<t_T$ for which the
    boundary between two adjacent voxels is crossed.  \blt[sum]
    Since
    $\|\loc(t_i)-\loc(t_{i-1})\| =
        (t_i-t_{i-1})\|\loc_2-\loc_1\|$ whenever  $t_i>t_{i-1}$, the approximation is then
    \begin{salign}
        \shad(\loc_1,\loc_2)&\approx
        \frac{\sum_{i=2}^T(t_i-t_{i-1})\|\loc_2-\loc_1\|
        \slf(\slfgridpt_{\slfgridptind_i})}{\|\loc_2-\loc_1\|^{1/2}}\\
        &= \|\loc_2-\loc_1\|^{1/2}\sum_{i=2}^T(t_i-t_{i-1})\slf(\slfgridpt_{\slfgridptind_i}),
    \end{salign}
    where ${\slfgridptind_i}$ is the index of the voxel that
    contains the $i$-th segment $\{\loc(t)~:~t\in (t_{i-1},t_i)\}$.
    \blt[crossings]Since $\slfgrid$ is a 3D grid, each point in
    $\{\slfgridpt_1,\ldots,\slfgridpt_\slfgridptnum\}$ can also be
    indexed by a vector $\vivec$ of 3 indices that lies in the set
    $\viset\define\{1,\ldots,\slfgridsidex\}\times\{1,\ldots,\slfgridsidey\}\times\{1,\ldots,\slfgridsidez\}$. The
    values of the SLF can also be collected in a tensor
    $\slften\in\rfield^{\slfgridsidex\times \slfgridsidey\times
            \slfgridsidez}$, whose entry $\slftens[\vivec]$ is the value of
    $\slf$ at the $\vivec$-th grid point.  If
    $\deltagridvec\in \rfield_{++}^3$ denotes a vector whose $j$-th
    entry $\deltagrids\entnot{j}$ represents the spacing between grid
    points along the $j$-th axis, the coordinates of the $\vivec$-th
    grid point are clearly $\vivec\odot\deltagridvec$, where $\odot$
    denotes entrywise product. Similarly, the boundaries between
    adjacent voxels along the $j$-th axis occur at values of the
    $j$-th coordinate given by $\deltagrids\entnot{j}(\vis \pm 1/2)$,
    where $\vis$ is an integer.  It is then clear that
    steps~\ref{step:crossings}-\ref{step:nextcrossing} in
    Algorithm~\ref{algo:tomo} simply find the next value of $t$ for
    which the segment crosses a voxel boundary along one of the axes
    by solving the equation
    \begin{align}
        \locs_1\entnot{j} + t  (\locs_2\entnot{j}-\locs_1\entnot{j})=
        \deltagrids\entnot{j}(\vicurrents\entnot{j} \pm {1}/{2})
    \end{align}
    for $t$ along each axis $j$ and taking the minimum across
    axes. The $\pm$ becomes a plus sign for the $j$-th axis if the
    segment is increasing along this axis (i.e.
    $\locs_1\entnot{j}\leq \locs_2\entnot{j}$)  and
    a minus sign otherwise.

    \blt[alternative]An alternative implementation of the same
    integral approximation with smaller computational complexity but greater
    memory complexity could be obtained by creating 3 lists
    corresponding to the values of $t$ for which the line segment
    between $\loc_1$ and $\loc_2$ intersects each axis and then
    merging those lists into a list with non-decreasing values of $t$.

\end{bullets}%

\newcommand{\stopMatDual}{\hc{\bm P}}
\newcommand{\stopMatPrimal}{\hc{\bm Q}}

\section{Interior-Point Solver}
\label{app:interior}
\begin{bullets}

    \blt[overview]The present section illustrates how
    \eqref{eq:mainproblemeq} can be solved using an interior-point
    algorithm. Although such a solver is not utilized in this paper, the
    ensuing derivation provides its computational complexity, which
    motivates the ADMM algorithm from Sec.~\ref{sec:solver}.

    \blt[slack]It is convenient to start by expressing \eqref{eq:mainproblemeq} in a
    canonical form with only non-negativity constraints and linear
    equality constraints. To this end, introduce the slack variables
    $\dslackovec$, $\dslacktmat$, and $\dslackthmat$ to write
    \eqref{eq:mainproblemeq} as
    \begin{subequations}
        \label{eq:mainproblemslack}
        \begin{align}
            \minimize_{
                \ratemat, \slackvec,
                \dslackovec, \dslacktmat,
                \dslackthmat
            }~
                 & \sparsweightvec\transpose\slackvec                   \\
            \st~ & \ratemat\transpose \bm 1 + \dslackovec = \maxratevec \\
                 & \ratemat \bm 1 = \minrate \bm 1
            \\
                 &
            \ratemat + \dslacktmat = \capmat                            \\
                 & \ratemat + \dslackthmat =  \bm 1\slackvec\transpose  \\
                 &
            \ratemat\geq \bm 0, \dslackovec\geq \bm 0,~\dslacktmat\geq \bm 0,~\dslackthmat\geq \bm 0.
        \end{align}
    \end{subequations}

    \blt[kron]With this formulation, it is easy to see that
    \eqref{eq:mainproblemslack} is equivalent to
    \begin{subequations}
        \label{eq:mainproblemsuper}
        \begin{align}
            \minimize_{\supervarvec}~
                 & \superweightvec\transpose \supervarvec              \\
            \st~ & \supermat\supervarvec=\supervec                     \\
                 & \supervarvec [ \gridptnum+1:\text{end}] \geq \bm 0,
        \end{align}
    \end{subequations}
    where
    \begin{bullets}
        \blt[]$\supervarvec\define[\slackvec\transpose,
                \ratevec\transpose,\dslackovec\transpose,\vect\transpose(\dslacktmat),\vect\transpose(\dslackthmat)]\transpose$,
        \blt $\ratevec=\vect(\ratemat)$,
        \blt[]$\superweightvec\define[\sparsweight_1,\ldots,\sparsweight_\gridptnum,0,\ldots,0]\transpose$,
        \blt[]
        $\supervarvec [ \gridptnum+1:\text{end}]\define
            [\ratevec\transpose,\dslackovec\transpose,\vect\transpose(\dslacktmat),\vect\transpose(\dslackthmat)]\transpose$,
        \blt[]$\supervec\define[(\maxratevec)\transpose,\minrate\bm
                1\transpose,\vect\transpose(\capmat),\bm 0\transpose]\transpose$, and
        \blt[]
        \begin{align}
            \supermat \define \left[
                \begin{array}{c c c c c}
                    \bm 0                            & \bm I_\gridptnum \otimes \bm 1\transpose & \bm I_\gridptnum & \bm 0                      & \bm 0                      \\
                    \bm 0                            & \bm 1\transpose \otimes \bm I_\usernum   & \bm 0            & \bm 0                      & \bm 0                      \\
                    \bm 0                            & \bm I_{\gridptnum\usernum}               & \bm 0            & \bm I_{\gridptnum\usernum} & \bm 0                      \\
                    - \bm I_\gridptnum \otimes \bm 1 & \bm I_{\gridptnum\usernum}               & \bm 0            & \bm 0                      & \bm I_{\gridptnum\usernum}
                \end{array}
                \right].
        \end{align}

        \blt[complexity]Problem \eqref{eq:mainproblemsuper} can be solved by
        means of a standard interior-point algorithm. To derive a lower
        bound for its computational complexity, note that each inner
        iteration of the algorithm will involve solving a system of
        equations where the number of unknowns equals the number of
        variables of the optimization problem plus the number of linear
        constraints~\cite[Ch.~10 and 11]{boyd}. For
        \eqref{eq:mainproblemsuper}, the former equals
        $2\gridptnum+3\gridptnum\usernum $ whereas the latter is given by
        $\gridptnum + \usernum + 2\gridptnum\usernum$. Solving this system
        of equations without any tailor-made approach that exploits the
        specific structure of $\supermat$ in this problem therefore involves
        $\mathcal{O}(\gridptnum^3\usernum^3)$ arithmetic operations.

        It is worth remarking that this complexity is  prohibitive in practice: if,
        for example, $\gridptnum=\usernum=100$, then $10^{12}$ operations
        would be required per inner iteration.
    \end{bullets}
\end{bullets}

\section{Stopping Criterion}
\label{sec:stopping}
\begin{bullets}%
    \blt[admm] The stopping criterion of Algorithm
    \ref{algo:placement} follows the framework in
    \cite{boyd2011distributed}. Particularly, given the absolute and
    relative tolerance parameters $\admmstop_{\rm{abs}}$ and
    $\admmstop_{\rm{rel}}$, let $\admmstop_{\rm{pri}}\itnot{\itind+1}$ and
    $\admmstop_{\rm{dual}}\itnot{\itind+1}$ be
    \begin{subequations}
        \begin{align}
            \admmstop_{\rm{pri}}\itnot{\itind+1}  & \define \sqrt{\usernum\gridptnum }\admmstop_{\rm{abs}} + \admmstop_{\rm{rel}}\max\{ \| \Amat_1\Xmat\itnot{\itind+1}\Amat_2 \|_\frob, \| \Bmat_1\Zmat\itnot{\itind+1}\Bmat_2 \|_\frob  \}, \\
            \admmstop_{\rm{dual}}\itnot{\itind+1} & \define \sqrt{\usernum\gridptnum }\admmstop_{\rm{abs}} + \admmstop_{\rm{rel}}\| \admmstep\Amat_1\transpose\Umat\itnot{\itind+1}\Amat_2\transpose \|_\frob.
        \end{align}
    \end{subequations}
    Algorithm \ref{algo:placement} stops when both conditions
    \begin{subequations}
        \begin{align}
            \| \stopMatPrimal\itnot{\itind+1} \|_\frob^2  \leqslant \admmstop_{\rm{pri}}\itnot{\itind+1}\text{ and }
            \| \stopMatDual\itnot{\itind+1} \|_\frob^2  \leqslant \admmstop_{\rm{dual}}\itnot{\itind+1}
        \end{align}
    \end{subequations}
    are satisfied,
    where
    \begin{subequations}
        \label{eq:stopDefine}
        \begin{align}
            \stopMatPrimal\itnot{\itind+1} & \define \Amat_1\Xmat\itnot{\itind+1}\Amat_2 + \Bmat_1\Zmat\itnot{\itind+1}\Bmat_2                                \\
            \stopMatDual\itnot{\itind+1}   & \define \admmstep\Amat_1\transpose\Bmat_1 (\Zmat\itnot{\itind+1} - \Zmat\itnot{\itind}) \Bmat_2\Amat_2\transpose
        \end{align}
    \end{subequations}
    are the so-called primal and dual residuals.

\end{bullets}

\section{Extended Proof of  Theorem \ref{prop:nphard}}
\label{proof:nphard}

The idea is to establish that a special case of
\eqref{eq:problemfalphamat} is a multidimensional
knapsack problem. To this end, let the $\gridptind$-th
entry of $\maxratevec$ be at least as large as
$\bm 1\transpose \capvec_\gridptind$ and
note that, due to \eqref{eq:problemfalphamatrange},
constraint \eqref{eq:problemfalphamatmaxrate} holds
regardless of the choice of
$\{\acti_\gridptind\}_{\gridptind=1}^\gridptnum$ and
$\ratemat$, meaning that
\eqref{eq:problemfalphamatmaxrate} can be removed.

Next, note that if
$\{\acti_\gridptind\}_{\gridptind=1}^\gridptnum$ and
$\ratemat$ are feasible, then replacing any
$\ratevec_\gridptind$ with
$\acti_\gridptind\capvec_\gridptind$ yields another
feasible point that attains the same cost. This is because
none of the entries of the left-hand side of
\eqref{eq:problemfalphamatminrate} decreases after
modifying $\ratevec_\gridptind$  in this way.

The left-hand side of \eqref{eq:problemfalphamatminrate}
can then be written as $\ratemat \bm 1=\sum_\gridptind
    \ratevec_\gridptind = \sum_\gridptind
    \acti_\gridptind\capvec_\gridptind$, which yields the
following problem
\begin{subequations}
    \label{eq:problemfalphamatsimp}
    \begin{align}
        \minimize_{
        \{\acti_\gridptind\}_{\gridptind=1}^\gridptnum}~ & \sum_{\gridptind}^{}\acti_\gridptind \\
        \st~
                                                         & \sum_\gridptind
        \acti_\gridptind\capvec_\gridptind \geq \minrate \bm 1\label{eq:problemfalphamatminsimprate}
        \\
                                                         & \acti_\gridptind \in \{0,1\}.
    \end{align}
\end{subequations}

Finally, applying the change of variables
$\invacti_\gridptind\leftarrow 1-\acti_\gridptind$, the
objective becomes $\gridptnum-\sum_\gridptind
    \invacti_\gridptind$ and the left-hand side of
\eqref{eq:problemfalphamatminsimprate} becomes $\sum_\gridptind
    (1-\invacti_\gridptind)\capvec_\gridptind
    =\sum_\gridptind \capvec_\gridptind - \sum_\gridptind
    \invacti_\gridptind\capvec_\gridptind $, which implies
that \eqref{eq:problemfalphamatsimp} reads as
\begin{subequations}
    \begin{align}
        \maximize_{
        \{\invacti_\gridptind\}_{\gridptind=1}^\gridptnum}~ & \sum_{\gridptind}^{}\invacti_\gridptind \\
        \st~
                                                            &
        \sum_\gridptind
        \invacti_\gridptind\capvec_\gridptind \leq \sum_\gridptind \capvec_\gridptind - \minrate \bm 1
        \\
                                                            & \invacti_\gridptind \in \{0,1\}.
    \end{align}
\end{subequations}
This problem is an instance of the so-called
multidimensional knapsack problem, which has been shown
to be NP-hard unless P$=$NP~\cite{gens1980complexity}.

\section{Extended Proof of Proposition~\ref{prop:eqxstep1}}
\label{proof:eqxstep1}

Since Problem \eqref{eq:ratevecindividual} is convex differentiable
and Slater's conditions are satisfied, it follows that the
Karush-Kuhn-Tucker (KKT) conditions are sufficient and
necessary~\cite[Sec.~5.5.3]{boyd}. To obtain these
conditions, observe that the Lagrangian of
\eqref{eq:ratevecindividual} is given by
\begin{align}
    \lagrangian(\ratevec_\gridptind, \slack_\gridptind;\nuvec)
    = \sparsweight_\gridptind \slack_\gridptind
    + \frac{\admmstep}{2}\|\ratevec_\gridptind - \zvec_\gridptind\itnot{\itind
    }+\uvec_\gridptind\itnot{\itind}\|_2^2 + \nuvec\transpose(   \ratevec_\gridptind- \slack_\gridptind \bm 1)
\end{align}
and note that the KKT conditions can be stated as
\begin{subequations}
    \label{eq:ubkkt}
    \begin{align}
        \label{eq:ubkktgrate}
        \nabla_{\ratevec_\gridptind}\lagrangian(\ratevec_\gridptind, \slack_\gridptind;\nuvec) & =
        \admmstep(\ratevec_\gridptind - \zvec_\gridptind\itnot{\itind
        }+\uvec_\gridptind\itnot{\itind}) + \nuvec = \bm 0                                                                                                                                                       \\
        \label{eq:ubkktslack}
        \nabla_{\slack_\gridptind}\lagrangian(\ratevec_\gridptind, \slack_\gridptind;\nuvec)   & =
        \sparsweight_\gridptind  - \bm 1\transpose \nuvec=0                                                                                                                                                      \\
        \label{eq:ubkktprimal}
        \ratevec_\gridptind                                                                    & \leq \slack_\gridptind \bm 1                                                                                    \\
        \label{eq:ubkktnu}
        \nuvec                                                                                 & \geq \bm 0,~~\nus\entnot{\userind}(\rate_\gridptind\entnot{\userind} - \slack_\gridptind) = 0~\forall \userind.
    \end{align}
\end{subequations}

From \eqref{eq:ubkktgrate} and the inequality in \eqref{eq:ubkktnu}, it follows that
\begin{align}
    \label{eq:ubkktvnucomb}
    \nuvec = -\admmstep(\ratevec_\gridptind - \zvec_\gridptind\itnot{\itind
    }+\uvec_\gridptind\itnot{\itind})\geq \bm 0.
\end{align}
This implies that
$ \ratevec_\gridptind\leq \zvec_\gridptind\itnot{\itind } -
    \uvec_\gridptind\itnot{\itind}$.  Combining this inequality with
\eqref{eq:ubkktprimal} yields
\begin{align}
    \label{eq:ubkktratebound}
    \ratevec_\gridptind\leq \min(\zvec_\gridptind\itnot{\itind } -
    \uvec_\gridptind\itnot{\itind}, \slack_\gridptind \bm 1).
\end{align}
On the other hand, from the
equality in \eqref{eq:ubkktvnucomb} and the equality in
\eqref{eq:ubkktnu}, one finds that
\begin{align}
    -\admmstep(\rate_\gridptind\entnot{\userind} - \zs_\gridptind\itnot{\itind
    }\entnot{\userind}+\us_\gridptind\itnot{\itind}\entnot{\userind})(\rate_\gridptind\entnot{\userind} - \slack_\gridptind) = 0~\forall \userind.
\end{align}
This holds if and only if either $\rate_\gridptind\entnot{\userind} = \zs_\gridptind\itnot{\itind
    }\entnot{\userind}-\us_\gridptind\itnot{\itind}\entnot{\userind}$ or
$\rate_\gridptind\entnot{\userind} = \slack_\gridptind$. Therefore, it follows from \eqref{eq:ubkktratebound} that
\begin{align}
    \label{eq:ubkktrateeq}
    \ratevec_\gridptind= \min(\zvec_\gridptind\itnot{\itind } -
    \uvec_\gridptind\itnot{\itind}, \slack_\gridptind \bm 1),
\end{align}
which establishes \eqref{eq:xkktr}. Finally, combine this expression with
\eqref{eq:ubkktslack} and \eqref{eq:ubkktvnucomb} to arrive at
\begin{subequations}
    \begin{align}
        \sparsweight_\gridptind & =  -\admmstep \bm 1\transpose(\ratevec_\gridptind - \zvec_\gridptind\itnot{\itind
        }+\uvec_\gridptind\itnot{\itind})                                                                                          \\
                                & =  -\admmstep \bm 1\transpose(\min(\zvec_\gridptind\itnot{\itind } -
            \uvec_\gridptind\itnot{\itind}, \slack_\gridptind \bm 1) - \zvec_\gridptind\itnot{\itind
        }+\uvec_\gridptind\itnot{\itind})                                                                                          \\
                                & =  -\admmstep \bm 1\transpose\min(\bm 0, \slack_\gridptind \bm 1 - \zvec_\gridptind\itnot{\itind
        }+\uvec_\gridptind\itnot{\itind})                                                                                          \\
        \label{eq:ubslackeqp}
                                & =  \admmstep \bm 1\transpose\max(\bm 0,  \zvec_\gridptind\itnot{\itind
        }-\uvec_\gridptind\itnot{\itind}-\slack_\gridptind \bm 1 ),
    \end{align}
\end{subequations}
thereby recovering \eqref{eq:xkks}. The proof is complete by noting
that \eqref{eq:ubkkt} holds if and only if \eqref{eq:ubkktrateeq} and
\eqref{eq:ubslackeqp} hold.

\section{Extended Proof of Proposition~\ref{prop:rootssg}}
\label{proof:rootssg}

Consider the function
$\Ffun(\slack)\define \bm 1\transpose \max(\zvec_\gridptind\itnot{\itind
    }-\uvec_\gridptind\itnot{\itind} -\slack\bm 1, \bm 0) = \sum_\userind
    \max(\zs_\gridptind\itnot{\itind
    }\entnot{\userind}-\us_\gridptind\itnot{\itind}\entnot{\userind}
    -\slack, 0)$. Since $\Ffun$ is the sum of non-increasing piecewise
linear functions, so is $\Ffun$. Since $\Ffun(\slack)\rightarrow \infty$
as $\slack\rightarrow -\infty$ and $\Ffun(\slack)=0$ for a sufficiently
large $\slack$, it follows that \eqref{eq:xkks} has at least one
root. Uniqueness of the root follows readily by noting that
$\Ffun$ is strictly decreasing  whenever $\Ffun(\slack)>0$.

It remains to be shown that
$\Ffun(\slacklow\gii)\geq{\sparsweight_\gridptind}/
            \admmstep$ whereas
$\Ffun(\slackhigh\gii)\leq {\sparsweight_\gridptind}/
            \admmstep$.
\begin{bullets}%
    \blt[lower]
    For the first of these inequalities, observe that
    $\slacklow\gii \leq \zs\gii\entnot{\userind} -
        \us\gii\entnot{\userind} - {\sparsweight_\gridptind}/({\usernum
                \admmstep})$ for all $\userind$, which in turn implies that
    $\zs\gii\entnot{\userind} -
        \us\gii\entnot{\userind} - \slacklow\gii \geq {\sparsweight_\gridptind}/({\usernum
                \admmstep})$. Thus, $\max(\zs\gii\entnot{\userind} -
        \us\gii\entnot{\userind} - \slacklow\gii,0)=\zs\gii\entnot{\userind} -
        \us\gii\entnot{\userind} - \slacklow\gii \geq {\sparsweight_\gridptind}/({\usernum
                \admmstep})$, which yields $\Ffun(\slacklow\gii)\geq \sum_\userind{\sparsweight_\gridptind}/({\usernum
                \admmstep}) = {\sparsweight_\gridptind}/{
        \admmstep}$.
    \blt[upper]For the second inequality, note similarly that   $\zs\gii\entnot{\userind} -
        \us\gii\entnot{\userind} - \slackhigh\gii \leq {\sparsweight_\gridptind}/({\usernum
                \admmstep})$ for all $\userind$. This means that $\Ffun(\slackhigh\gii)\leq \sum_\userind\max(
        {\sparsweight_\gridptind}/({\usernum
                \admmstep}),0) =  {\sparsweight_\gridptind}/{
        \admmstep} $.

\end{bullets}%

\section{Extended Proof of Proposition~\ref{prop:eqxstep2}}
\label{proof:eqxstep2}

Again, the KKT conditions are sufficient and
necessary. Since the Lagrangian is
\begin{align}
    \lagrangian(\ratevec_\gridptind, \slack_\gridptind;\nuvec)\nonumber
    = & ~ \sparsweight_\gridptind \slack_\gridptind
    + \frac{\admmstep}{2}\|\ratevec_\gridptind - \zvec_\gridptind\itnot{\itind
    }+\uvec_\gridptind\itnot{\itind}\|_2^2                               \\
      & + \nuvec\transpose(   \ratevec_\gridptind- \slack_\gridptind \bm
    1) + \mumult (\bm 1\transpose\ratevec_\gridptind - \maxrate{\gridptind}),
\end{align}
the KKT conditions read as
\begin{subequations}
    \label{eq:equbkkt}
    \begin{align}
        \label{eq:equbkktgrate}
        \nabla_{\ratevec_\gridptind}\lagrangian(\ratevec_\gridptind, \slack_\gridptind;\nuvec) & =
        \admmstep(\ratevec_\gridptind - \zvec_\gridptind\itnot{\itind
        }+\uvec_\gridptind\itnot{\itind})
        +
        \nuvec
        +
        \mumult
        \bm
        1=
        \bm 0                                                                                                                 \\
        \label{eq:equbkktslack}
        \nabla_{\slack_\gridptind}\lagrangian(\ratevec_\gridptind, \slack_\gridptind;\nuvec)   & =
        \sparsweight_\gridptind  - \bm 1\transpose \nuvec=0                                                                   \\
        \label{eq:equbkktprimal}
        \ratevec_\gridptind                                                                    & \leq \slack_\gridptind \bm 1 \\
        \label{eq:equbkktnu}
        \nuvec                                                                                 & \geq \bm
        0,~~\nus\entnot{\userind}(\rate_\gridptind\entnot{\userind}
        - \slack_\gridptind) = 0~\forall \userind                                                                             \\
        \label{eq:equbkktmaxrate}
        \bm 1\transpose\ratevec_\gridptind                                                     & = \maxrate{\gridptind}.
    \end{align}
\end{subequations}

From \eqref{eq:equbkktgrate} and the inequality in \eqref{eq:equbkktnu}, it follows that
\begin{align}
    \label{eq:equbkktvnucomb}
    \nuvec = -\admmstep(\ratevec_\gridptind - \zvec_\gridptind\itnot{\itind
    }+\uvec_\gridptind\itnot{\itind})-\mumult \bm 1\geq \bm 0.
\end{align}
This implies that
$ \ratevec_\gridptind\leq \zvec_\gridptind\itnot{\itind } -
    \uvec_\gridptind\itnot{\itind} - (\mumult/\admmstep)\bm 1$.  Combining this inequality with
\eqref{eq:equbkktprimal} yields
\begin{align}
    \label{eq:equbkktratebound}
    \ratevec_\gridptind\leq \min(\zvec_\gridptind\itnot{\itind } -
    \uvec_\gridptind\itnot{\itind} - (\mumult/\admmstep)\bm 1, \slack_\gridptind \bm 1).
\end{align}
On the other hand, from the
equality in \eqref{eq:equbkktvnucomb} and the equality in
\eqref{eq:equbkktnu}, one finds that
\begin{align}
    [-\admmstep(\rate_\gridptind\entnot{\userind} - \zs_\gridptind\itnot{\itind
        }\entnot{\userind}+\us_\gridptind\itnot{\itind}\entnot{\userind})-\mumult](\rate_\gridptind\entnot{\userind} - \slack_\gridptind) = 0~\forall \userind.
\end{align}
This holds if and only if either $\rate_\gridptind\entnot{\userind} = \zs_\gridptind\itnot{\itind
    }\entnot{\userind}-\us_\gridptind\itnot{\itind}\entnot{\userind} - \mumult/\admmstep$ or
$\rate_\gridptind\entnot{\userind} = \slack_\gridptind$. Therefore,
it follows that
\begin{align}
    \label{eq:eqratevecxj}
    \ratevec_\gridptind = \min(\zvec_\gridptind\itnot{\itind
    }-\uvec_\gridptind\itnot{\itind} - (\mumult/\admmstep)\bm 1, \slack_\gridptind\bm 1),
\end{align}
which establishes \eqref{eq:eqxkktr}.

To find $\mumult$, substitute the equality in
\eqref{eq:equbkktvnucomb} into \eqref{eq:equbkktslack} to obtain
\begin{align}
    \label{eq:eqcondxyz}
    \bm 1\transpose[-\admmstep(\ratevec_\gridptind - \zvec_\gridptind\itnot{\itind
        }+\uvec_\gridptind\itnot{\itind})-\mumult \bm 1] =\sparsweight_\gridptind.
\end{align}
Solving for $\mumult$ yields
\begin{align}
    \mumult 
     & =\frac{ -\admmstep\bm 1\transpose\ratevec_\gridptind
        +\admmstep\bm 1\transpose ( \zvec_\gridptind\itnot{\itind
        }-\uvec_\gridptind\itnot{\itind})-\sparsweight_\gridptind}{\usernum}
\end{align}
and using \eqref{eq:equbkktmaxrate} results in \eqref{eq:eqmumult}.

Finally, substitute \eqref{eq:eqratevecxj}  into \eqref{eq:eqcondxyz} to
arrive at
\begin{align}
    \sparsweight_\gridptind 
     & =
    \bm 1\transpose\max(  \mumult\bm 1, \admmstep(\zvec_\gridptind\itnot{\itind
        }-\uvec_\gridptind\itnot{\itind}-\slack_\gridptind\bm 1)) -\mumult\usernum,
\end{align}
thereby recovering \eqref{eq:eqxkks}. The proof is complete by noting
that \eqref{eq:equbkkt} holds if and only if \eqref{eq:eqxkktr} and
\eqref{eq:eqxkks} hold.

\section{Extended Proof of Proposition~\ref{prop:eqrootssg}}
\label{proof:eqrootssg}

Consider the function
$\Ffun(\slack)\define \bm 1\transpose\max( \mumult\bm 1,
    \admmstep(\zvec_\gridptind\itnot{\itind
        }-\uvec_\gridptind\itnot{\itind}-\slack\bm 1)) = \sum_\userind
    \max(\mumult, \admmstep(\zs_\gridptind\itnot{\itind
        }\entnot{\userind}-\us_\gridptind\itnot{\itind}\entnot{\userind}
        -\slack))$. Due to the same argument as in the proof of
\thref{prop:rootssg}, this function has a unique root.

\begin{bullets}%

    \blt[lower]To show that $\Ffun(\slacklow\gii)\geq  \sparsweight_\gridptind+{\mumult}\usernum$, observe that   \begin{subequations}
        \begin{align}
            \Ffun(\slacklow\gii) & \geq \sum_\userind \admmstep(\zs_\gridptind\itnot{\itind
            }\entnot{\userind}-\us_\gridptind\itnot{\itind}\entnot{\userind}
            -\slacklow\gii)                                                                                                                      \\
                                 & \geq \usernum \min_\userind[ \admmstep(\zs_\gridptind\itnot{\itind
                }\entnot{\userind}-\us_\gridptind\itnot{\itind}\entnot{\userind}
            -\slacklow\gii)]                                                                                                                     \\
                                 & =\usernum \admmstep\left( \frac{\sparsweight_\gridptind}{\usernum \admmstep}+\frac{\mumult}{\admmstep}\right) \\
                                 & ={\sparsweight_\gridptind}+{\mumult}\usernum.
        \end{align}
    \end{subequations}

    \blt[upper]To show that $\Ffun(\slackhigh\gii)\leq  \sparsweight_\gridptind+{\mumult}\usernum$, observe that   \begin{subequations}
        \begin{align}
            \Ffun(\slackhigh\gii) & \leq \usernum \max_\userind[ \max(\mumult, \admmstep(\zs_\gridptind\itnot{\itind
                }\entnot{\userind}-\us_\gridptind\itnot{\itind}\entnot{\userind}
            -\slackhigh\gii))]                                                                                                                                               \\
                                  & = \usernum  \max(\mumult, \admmstep(\max_\userind[\zs_\gridptind\itnot{\itind
                    }\entnot{\userind}-\us_\gridptind\itnot{\itind}\entnot{\userind}]
            -\slackhigh\gii))                                                                                                                                                \\
                                  & = \usernum  \max\left(\mumult, \admmstep\left(\frac{\sparsweight_\gridptind}{\usernum \admmstep}+\frac{\mumult}{\admmstep}\right)\right) \\
                                  & \leq  \sparsweight_\gridptind+{\mumult}\usernum.
        \end{align}
    \end{subequations}

\end{bullets}%

\section{Extended Proof of Proposition~\ref{prop:zsol}}
\label{proof:zsol}

The fact that $\bm 1\transpose \caprvec_\userind< \minrate$ implies
that \eqref{eq:ubzrvecs} is infeasible is trivial and, therefore,
the rest of the proof focuses on the case where
$\bm 1\transpose \caprvec_\userind\geq \minrate$.

As before, the KKT conditions are sufficient and necessary in this
case. Noting that  the Lagrangian is given by
\begin{align}
     & \lagrangian(\zrvec_\userind;\lambdas, \nuvec,\muvec ) =~
    \frac{1}{2}\| \ratervec_\userind\itnot{\itind+1} -\zrvec_\userind + \urvec_\userind\itnot{\itind}\|_\frob^2 \nonumber \\&\quad\quad+ \lambdas( \bm 1\transpose \zrvec_\userind - \minrate) - \nuvec\transpose \zrvec_\userind + \muvec\transpose( \zrvec_\userind - \caprvec_\userind)
\end{align}
yields the KKT conditions
\begin{subequations}
    \begin{align}
        \nonumber
         & \nabla_{\zrvec_\userind}    \lagrangian(\zrvec_\userind;\lambdas, \nuvec,\muvec ) =                                                                                         \\
        \label{eq:ubkktz}
         & \quad
        -( \ratervec_\userind\itnot{\itind+1} -\zrvec_\userind + \urvec_\userind\itnot{\itind}) + \lambdas \bm 1 - \nuvec + \muvec = \bm 0,                                            \\
        \label{eq:ubkktzlambda}
         & \bm 1\transpose \zrvec_\userind = \minrate,                                                                                                                                 \\
        \label{eq:ubkktznu}
         & \zrvec_\userind \geq \bm 0,~\nuvec \geq 0,~\nus\entnot{\gridptind}\zrs_\userind\entnot{\gridptind}=0~\forall \gridptind,                                                    \\
        \label{eq:ubkktmu}
         & \zrvec_\userind \leq \caprvec_\userind,~\muvec \geq 0,~\mus\entnot{\gridptind}(\zrs_\userind\entnot{\gridptind} -\caprs_\userind\entnot{\gridptind} )=0~\forall \gridptind.
    \end{align}
\end{subequations}
From \eqref{eq:ubkktz} and the second inequality in \eqref{eq:ubkktmu}, it follows that
\begin{align}
    \label{eq:ubmuvecin}
    \muvec =  \ratervec_\userind\itnot{\itind+1} -\zrvec_\userind + \urvec_\userind\itnot{\itind}
    - \lambdas \bm 1 + \nuvec \geq \bm 0,
\end{align}
which in turn implies that
\begin{align}
    \zrvec_\userind \leq  \ratervec_\userind\itnot{\itind+1} + \urvec_\userind\itnot{\itind}
    - \lambdas \bm 1 + \nuvec .
\end{align}
Combining this expression with the first inequality in \eqref{eq:ubkktmu} yields
\begin{align}
    \zrvec_\userind \leq \min(\caprvec_\userind,  \ratervec_\userind\itnot{\itind+1} + \urvec_\userind\itnot{\itind}
    - \lambdas \bm 1 + \nuvec ).
\end{align}
To show that this expression holds with equality, substitute \eqref{eq:ubmuvecin} into the equality of \eqref{eq:ubkktmu} to obtain
\begin{align}
    (\raters_\userind\itnot{\itind+1}\entnot{\gridptind} -\zrs_\userind\entnot{\gridptind} + \urs_\userind\itnot{\itind}\entnot{\gridptind}
    - \lambdas  + \nus\entnot{\gridptind}
    )(\zrs_\userind\entnot{\gridptind} -\caprs_\userind\entnot{\gridptind} )=0,
\end{align}
which implies that either $ \zrs_\userind\entnot{\gridptind} = \raters_\userind\itnot{\itind+1}\entnot{\gridptind}  + \urs_\userind\itnot{\itind}\entnot{\gridptind}
    - \lambdas  + \nus\entnot{\gridptind}$ or $\zrs_\userind\entnot{\gridptind} =\caprs_\userind\entnot{\gridptind}$. Therefore,
\begin{align}
    \label{eq:ubzrvec}
    \zrvec_\userind = \min(\caprvec_\userind,  \ratervec_\userind\itnot{\itind+1} + \urvec_\userind\itnot{\itind}
    - \lambdas \bm 1 + \nuvec ).
\end{align}
To obtain an expression for $\zrvec_\userind$ that does not depend on
$\nuvec$, one may consider three cases for each $\gridptind$:
\begin{itemize}
    \item C1:
          $\raters_\userind\itnot{\itind+1}\entnot{\gridptind} +
              \urs_\userind\itnot{\itind}\entnot{\gridptind} - \lambdas< 0$. In
          this case, if $\nus\entnot{\gridptind}=0$, expression
          \eqref{eq:ubzrvec} would imply that
          $\zrs_\userind\entnot{\gridptind}<0$, which would violate the
          first inequality in \eqref{eq:ubkktznu}. Therefore,
          $\nus\entnot{\gridptind}>0$ and, due to the equality in
          \eqref{eq:ubkktznu}, $\zrs_\userind\entnot{\gridptind}=0$. If
          $\caprs_\userind\entnot{\gridptind}>0$, it is then clear from
          \eqref{eq:ubzrvec} that
          $ \nus\entnot{\gridptind} =
              -(\raters_\userind\itnot{\itind+1}\entnot{\gridptind} +
              \urs_\userind\itnot{\itind}\entnot{\gridptind} - \lambdas)$.  If
          $\caprs_\userind\entnot{\gridptind}=0$, then greater values of
          $\nus\entnot{\gridptind}$ will also satisfy the KKT conditions but
          this is not relevant since  the only feasible
          $\zrs_\userind\entnot{\gridptind}$ in case C1 is
          $\zrs_\userind\entnot{\gridptind}=0$.

    \item C2:
          $\raters_\userind\itnot{\itind+1}\entnot{\gridptind} +
              \urs_\userind\itnot{\itind}\entnot{\gridptind} - \lambdas = 0$. In
          this case, \eqref{eq:ubzrvec} becomes
          $\zrs_\userind\entnot{\gridptind} =
              \min(\caprs_\userind\entnot{\gridptind},
              \nus\entnot{\gridptind})$. Due to the equality in
          \eqref{eq:ubkktznu}, it then follows that either
          $\caprs_\userind\entnot{\gridptind}=0$ and
          $\nus\entnot{\gridptind}\geq 0$, or
          $\zrs_\userind\entnot{\gridptind} = \nus\entnot{\gridptind}=0$.

    \item C3:
          $\raters_\userind\itnot{\itind+1}\entnot{\gridptind} +
              \urs_\userind\itnot{\itind}\entnot{\gridptind} - \lambdas>0$. If
          $\caprs_\userind\entnot{\gridptind}=0$, then necessarily
          $\zrs_\userind\entnot{\gridptind} = 0$ and any
          $\nus\entnot{\gridptind}\geq 0$ satisfies the KKT conditions. On
          the other hand, if $\caprs_\userind\entnot{\gridptind}>0$, then it
          is clear that $\zrs_\userind\entnot{\gridptind} > 0$ and, due to
          the equality in \eqref{eq:ubkktznu}, one has that
          $\nus\entnot{\gridptind} =0$, which in turn implies that
          $\zrs_\userind\entnot{\gridptind} =
              \min(\caprs_\userind\entnot{\gridptind}, \raters_\userind\itnot{\itind+1}\entnot{\gridptind} +
              \urs_\userind\itnot{\itind}\entnot{\gridptind} - \lambdas)$.
\end{itemize}
Combining C1-C3 yields
\begin{align}
    \label{eq:ubcombzrs}
    \zrs_\userind\entnot{\gridptind} =
    \max(0,    \min(\caprs_\userind\entnot{\gridptind}, \raters_\userind\itnot{\itind+1}\entnot{\gridptind} +
        \urs_\userind\itnot{\itind}\entnot{\gridptind} - \lambdas)),
\end{align}
which is just the scalar version of
\eqref{eq:ubzrvecsol}. Finally, to obtain $\lambdas$, one may
substitute \eqref{eq:ubcombzrs} into \eqref{eq:ubkktzlambda},
which produces \eqref{eq:ubzrveclambda}.

\section{Extended Proof of Proposition~\ref{prop:ubrootlambda}}
\label{proof:ubrootlambda}

Denote by $\Gfun(\lambdas)$ the left-hand side of
\eqref{eq:ubzrveclambda}, i.e.,
\begin{align}
    \Gfun(\lambdas) \define \sum_\gridptind
    \max(0,    \min(\caprs_\userind\entnot{\gridptind}, \raters_\userind\itnot{\itind+1}\entnot{\gridptind} +
        \urs_\userind\itnot{\itind}\entnot{\gridptind} - \lambdas)).
\end{align}
This is a sum of non-increasing piecewise continuous functions and
therefore $\Gfun$ is also non-increasing piecewise continuous. The
maximum value is attained for sufficiently small $\lambdas$ and
equals
$ \sum_\gridptind\caprs_\userind\entnot{\gridptind} = \bm
    1\transpose \caprvec_\userind$. If
$ \bm 1\transpose \caprvec_\userind < \minrate$, then
$\Gfun(\lambdas)<\minrate~\forall \lambda$ and
\eqref{eq:ubzrveclambda} admits no solution. Conversely, if
$ \bm 1\transpose \caprvec_\userind > \minrate$, then a solution
can be found since $\Gfun(\lambdas)>\minrate$ for sufficiently
small $\lambdas$ and $\Gfun(\lambdas)=0$ for sufficiently large
$\lambdas$. Uniqueness follows from the fact that $\Gfun$ is
strictly decreasing except when $\Gfun(\lambdas)=0$ or
$\Gfun(\lambdas)=\bm 1\transpose \caprvec_\userind$.

\begin{bullets}%
    \blt[lower]    To show that
    $\Gfun(\lambdaslow_\userind\itnot{\itind})\geq \minrate$ just note
    from \eqref{eq:ublambdaslow} that
    $\lambdaslow_\userind\itnot{\itind} \leq
        \raters_\userind\itnot{\itind+1}\entnot{\gridptind} +
        \urs_\userind\itnot{\itind}\entnot{\gridptind}
        -\caprs_\userind\entnot{\gridptind}$ or, equivalently,
    $ \caprs_\userind\entnot{\gridptind} \leq
        \raters_\userind\itnot{\itind+1}\entnot{\gridptind} +
        \urs_\userind\itnot{\itind}\entnot{\gridptind}
        -\lambdaslow_\userind\itnot{\itind} $, for all $\gridptind$. This clearly yields
    $ \Gfun(\lambdaslow_\userind\itnot{\itind}) = \sum_\gridptind
        \max(0, \caprs_\userind\entnot{\gridptind}) = \sum_\gridptind
        \caprs_\userind\entnot{\gridptind} $, which is greater than or
    equal to $\minrate$ by assumption.

    \blt[upper]To show that
    $\Gfun(\lambdashigh_\userind\itnot{\itind})\leq \minrate$, note
    from \eqref{eq:ublambdashigh} that
    $ \lambdashigh_\userind\itnot{\itind} \geq
        \raters_\userind\itnot{\itind+1}\entnot{\gridptind} +
        \urs_\userind\itnot{\itind}\entnot{\gridptind}-{\minrate}/{\gridptnum}$ for all
    $\gridptind$ such that
    $
        \caprs_\userind\entnot{\gridptind}>{\minrate}/{\gridptnum}$. This
    clearly implies that
    $     \raters_\userind\itnot{\itind+1}\entnot{\gridptind} +
        \urs_\userind\itnot{\itind}\entnot{\gridptind} - \lambdashigh_\userind\itnot{\itind}\leq {\minrate}/{\gridptnum}$ for all $\gridptind$ such that
    $
        \caprs_\userind\entnot{\gridptind}>{\minrate}/{\gridptnum}$ and, as a consequence,
    $ \min(\caprs_\userind\entnot{\gridptind},
        \raters_\userind\itnot{\itind+1}\entnot{\gridptind} +
        \urs_\userind\itnot{\itind}\entnot{\gridptind} - \lambdashigh_\userind\itnot{\itind}))\leq
        {\minrate}/{\gridptnum}$ and the inequality $\Gfun(\lambdashigh_\userind\itnot{\itind})\leq \minrate$ follows.

\end{bullets}%

\section{Additional Experiments with the Tomographic Model}
\label{sec:additionaltomo}

\begin{bullets}

    \begin{figure}[t]
        \centering
        \captionsetup{width=.9\linewidth}
        \includegraphics[width=1.05\linewidth]{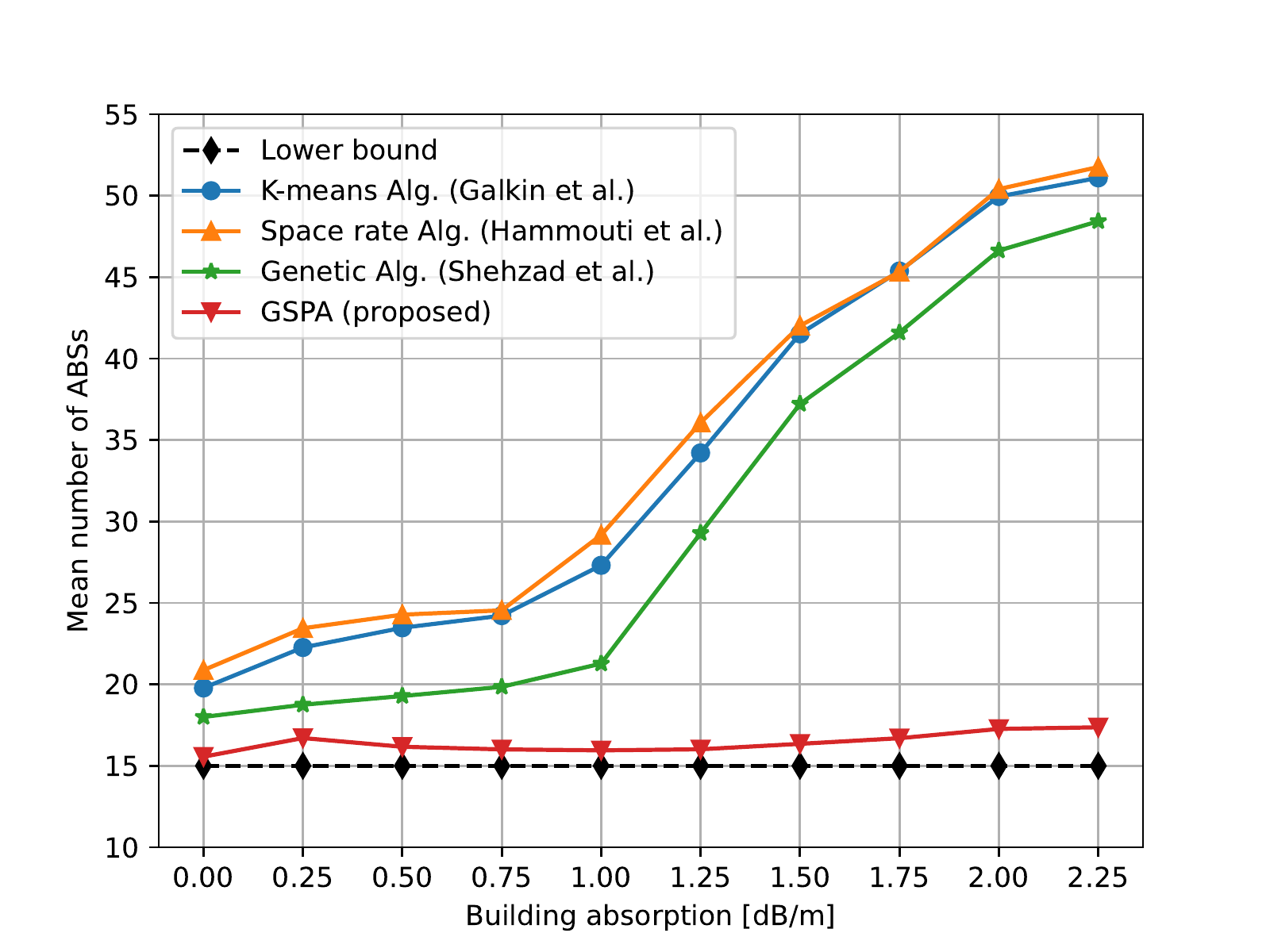}
        \captionof{figure}{Mean number of ABSs vs. building
            absorption  ($\minrate= 17$ Mbps, $\maxrate{}=
                84$ Mbps).}
        \label{fig:vs_absorption}
    \end{figure}
    \blt[Exp. - vs buildingAbsorption]To
    study the influence of the channel,
    \begin{bullets}%
        \blt[Figure]Fig.~\ref{fig:vs_absorption} shows
        the mean number of ABSs vs. the absorption undergone by
        the communication signals when propagating through the
        buildings.
        %
        \blt[Interpretation]
        \begin{bullets}
            \blt[low absorption]
            \begin{bullets}
                \blt[free space]When the absorption is zero, the
                propagation conditions are those of free
                space. In this case, the proposed algorithm
                outperforms the rest only because of a better
                ability to perform the rate allocation.
                \blt[not free space]As the absorption increases,
                the benchmark algorithms are dramatically
                affected, which suggests that these algorithms
                are not well suited to scenarios without
                line-of-sight. In contrast, the proposed
                algorithm remains unaffected since there are
                always sufficiently good flight grid points
                regardless of the building absorption considered
                in the figure. Informally speaking,
                matrix $\ratemat$ in \eqref{eq:mainproblem} is
                upper bounded by $\maxratevec$ and
                $\capmat$. When the former constraint is tighter
                than the latter, the resulting number of ABSs
                will not depend on $\capmat$.

            \end{bullets}
            %
        \end{bullets}
    \end{bullets}

    \blt[Exp. - vs min-user-rate]
    This phenomenon
    is investigated further in
    Fig.~\ref{fig:vs_minRate_GroupSparsePlacer_ray_tracing}, which is the
    counterpart of
    Fig.~\ref{fig:vs_minRate_GroupSparsePlacer_ray_tracing} for
    tomographic channels. The purpose of this simulation is to confirm
    that the
    approximate linearity and proximity to the bound of GSPA
    observed in  Fig.~\ref{fig:vs_minRate} take
    place for a wide range of parameters.

    \begin{figure}[t]
        \centering
        \includegraphics[width=.5\textwidth]{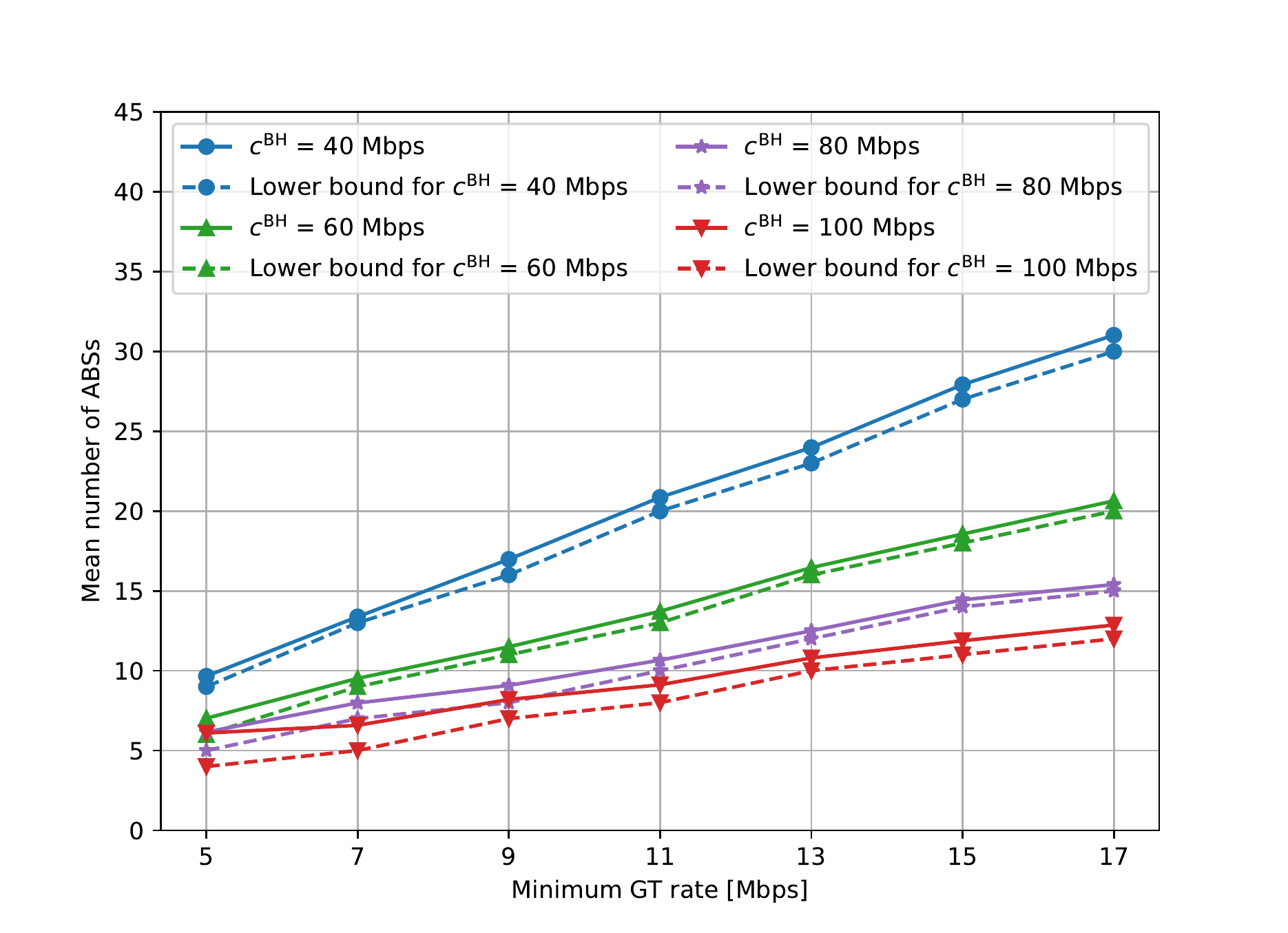}
        \caption{Mean number of ABSs vs. $\minrate$ of the
            proposed GSPA algorithm.}
        \label{fig:vs_minRate_GroupSparsePlacer}
    \end{figure}

    \blt[Exp. - vs minflyHeight]Fig. \ref{fig:vs_minFlyHeight}
    investigates the influence of the  minimum flight height
    on the performance of the considered algorithms.
    \begin{bullets}
        \blt[Figure]
        \blt[vs min fly height]
        \begin{figure}[t]
            \centering
            \includegraphics[width=.5\textwidth]{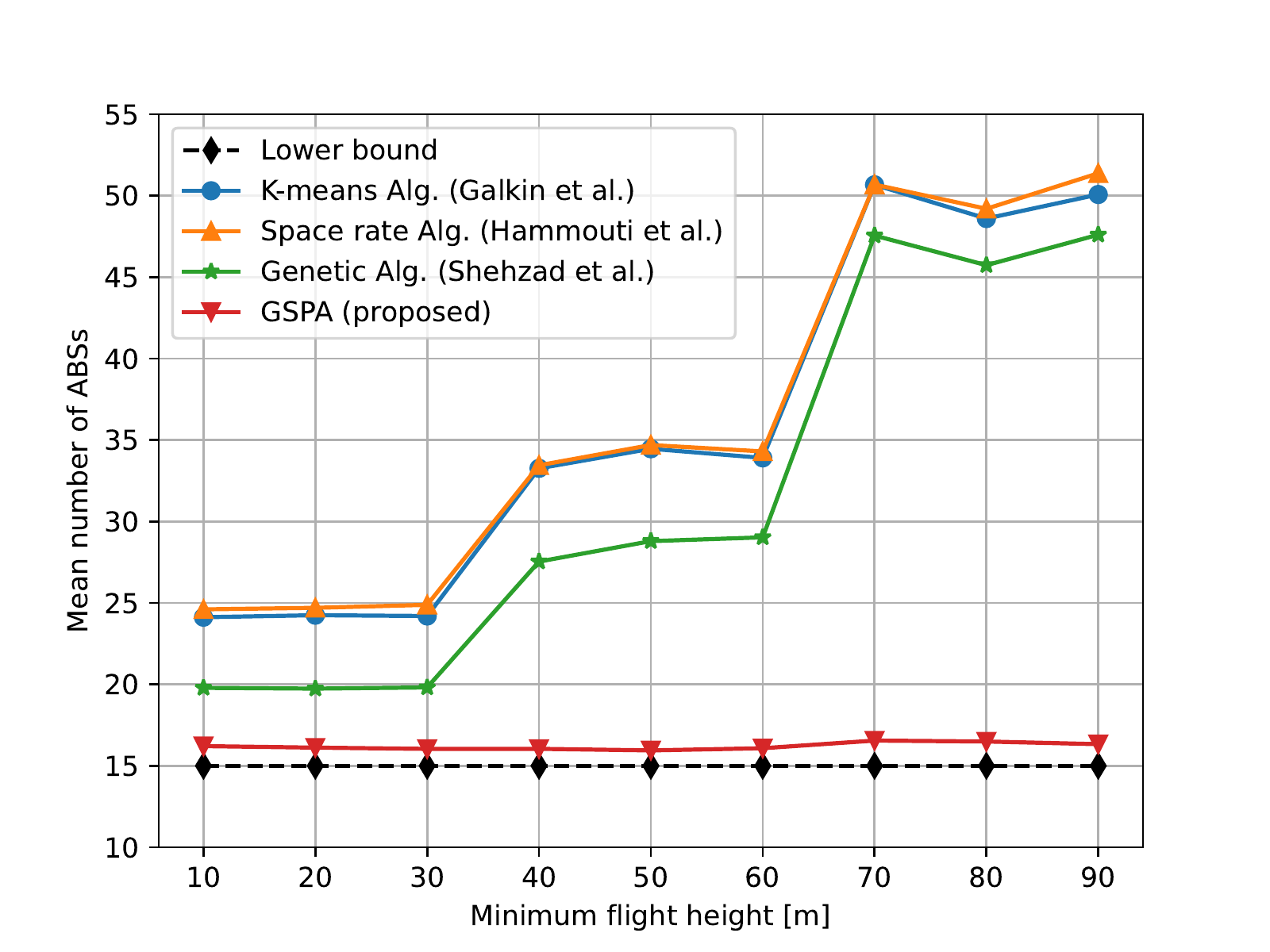}
            \caption{Mean number of ABSs vs. minimum flight
                height ($\minrate= 17$ Mbps, $\maxrate{}= 84$
                Mbps).}
            \label{fig:vs_minFlyHeight}
        \end{figure}
        \blt[Interpretation]
        \begin{bullets}%
            \blt[stairs]The staircase behavior of the benchmarks
            can be explained by noting that the flight grid
            initially comprises points with heights 0, 30 m, 60
            m, 90 m, etc. Then, the points inside buildings and
            the points below the minimum flight height are
            removed. Thus, the allowed flight points are the
            same e.g. when the minimum flight height is 10 m as
            when it is 20 m.

            \blt[constant] As already observed in
            Sec.~\ref{sec:experiments}, the performance of GSPA
            is not degraded for increasing flight height because
            the backhaul capacity poses a more stringent
            constraint than the one imposed by $\capmat$ even
            for the maximum flight height considered in the
            figure.
        \end{bullets}
    \end{bullets}

    \blt[Exp. - vs building height]
    \begin{bullets}
        \blt[How to run]Fig. \ref{fig:vs_buildingHeights}
        studies the impact of the height of the buildings. To
        reduce the spatial quantization effect of the
        tomographic integral approximation, the numbers of SLF
        grid points in the x, y, and z axes were set to 50, 40,
        and 150, respectively.
        \blt[Figure]
        \begin{figure}[t!]
            \centering
            \includegraphics[width=.5\textwidth]{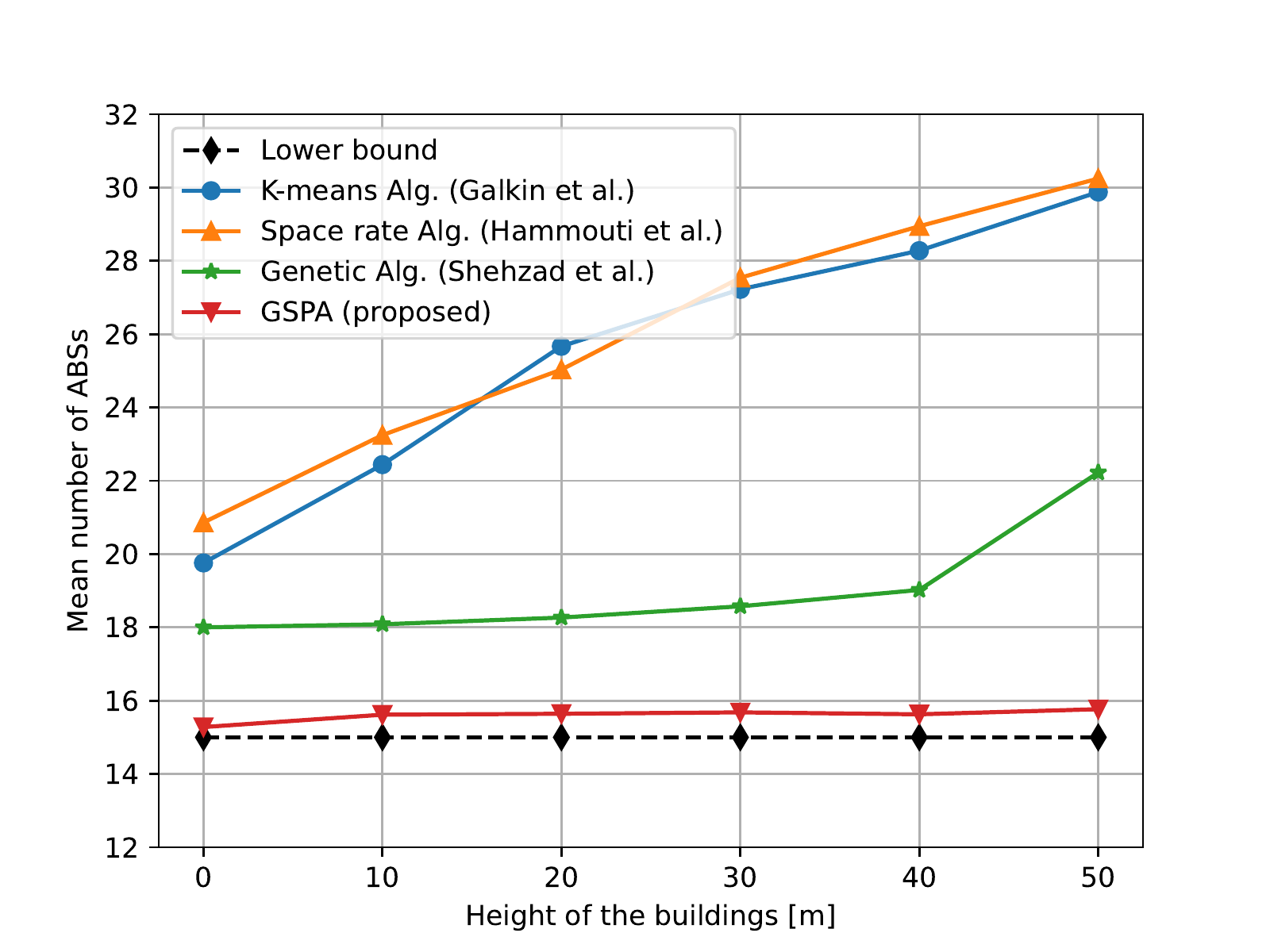}
            \caption{Mean number of ABSs vs. height of the
                buildings ($\minrate= 17$ Mbps, $\maxrate{}= 84$
                Mbps).
            }
            \label{fig:vs_buildingHeights}
        \end{figure}
        \blt[Interpretation]
        \begin{bullets}
            \blt[general]It is observed that the K-means and
            space rate K-means algorithms are negatively
            affected by the increased attenuation introduced by
            higher buildings.
            \blt[constant]In  contrast,  the  genetic  placement
            algorithm and  GSPA exhibit  a milder  dependence on
            the  building  height.  The reason  is  the  greater
            capacity of these algorithms to select ABS locations
            with favorable propagation  conditions. The backhaul
            capacity  is again  the  limiting constraint,  which
            explains why  the performance of GSPA  is unaffected
            by a greater building height.

            %
        \end{bullets}

    \end{bullets}

\end{bullets}


\begin{thebibliography}{10}

\bibitem{romero2022aerial}
D.~Romero, P.~Q. Viet, and G.~Leus,
\newblock ``Aerial base station placement leveraging radio tomographic maps,''
\newblock in {\em IEEE Int. Conf. Acoustics Speech Signal Process.}, Singapore,
  2022, IEEE, pp. 5358--5362.

\bibitem{zeng2019accessing}
Y.~Zeng, Q.~Wu, and R.~Zhang,
\newblock ``Accessing from the sky: A tutorial on uav communications for 5g and
  beyond,''
\newblock {\em Proceedings IEEE}, vol. 107, no. 12, pp. 2327--2375, 2019.

\bibitem{viet2022introduction}
P.~Q. Viet and D.~Romero,
\newblock ``Aerial base station placement: A tutorial introduction,''
\newblock {\em IEEE Commun. Mag.}, vol. 60, no. 5, pp. 44--49, 2022.

\bibitem{han2009manet}
Z.~Han, A.~L. Swindlehurst, and K.~J.~R. Liu,
\newblock ``Optimization of manet connectivity via smart deployment/movement of
  unmanned air vehicles,''
\newblock {\em IEEE Trans. Veh. Technol.}, vol. 58, no. 7, pp. 3533--3546,
  2009.

\bibitem{lee2011climbing}
D.-J. Lee,
\newblock ``Autonomous unmanned flying robot control for reconfigurable
  airborne wireless sensor networks using adaptive gradient climbing
  algorithm,''
\newblock {\em J. Korea Robotics Society}, vol. 6, no. 2, pp. 97–107, May
  2011.

\bibitem{boryaliniz2016placement}
I.~Bor-Yaliniz, A.~El-Keyi, and H.~Yanikomeroglu,
\newblock ``Efficient 3-d placement of an aerial base station in next
  generation cellular networks,''
\newblock in {\em Proc. IEEE Int. Conf. Commun.} IEEE, 2016, pp. 1--5.

\bibitem{chen2017map}
J.~Chen and D.~Gesbert,
\newblock ``Optimal positioning of flying relays for wireless networks: A {LOS}
  map approach,''
\newblock in {\em Proc. IEEE Int. Conf. Commun.}, Paris, France, May 2017, pp.
  1--6.

\bibitem{wang2018adaptive}
Z.~Wang, L.~Duan, and R.~Zhang,
\newblock ``Adaptive deployment for {UAV}-aided communication networks,''
\newblock {\em IEEE Trans. Wireless Commun.}, vol. 18, no. 9, pp. 4531--4543,
  2019.

\bibitem{lee2010decentralized}
D.-J. Lee and R.~Mark,
\newblock ``Decentralized control of unmanned aerial robots for wireless
  airborne communication networks,''
\newblock {\em Int. J. Advanced Robotic Syst.}, vol. 7, no. 3, pp. 22, 2010.

\bibitem{andryeyev2016selforganized}
O.~Andryeyev and A.~Mitschele-Thiel,
\newblock ``Increasing the cellular network capacity using self-organized
  aerial base stations,''
\newblock in {\em Proc. Workshop Micro Aerial Veh. Netw., Syst., Appl.} ACM,
  2017, pp. 37--42.

\bibitem{galkin2016deployment}
B.~Galkin, J.~Kibilda, and L.A. DaSilva,
\newblock ``Deployment of {UAV}-mounted access points according to spatial user
  locations in two-tier cellular networks,''
\newblock in {\em Wireless Days}. IEEE, 2016, pp. 1--6.

\bibitem{lyu2017mounted}
J.~Lyu, Y.~Zeng, R.~Zhang, and T.J. Lim,
\newblock ``Placement optimization of {UAV}-mounted mobile base stations,''
\newblock {\em IEEE Commun. Letters}, vol. 21, no. 3, pp. 604--607, 2017.

\bibitem{romero2019noncooperative}
D.~Romero and G.~Leus,
\newblock ``Non-cooperative aerial base station placement via stochastic
  optimization,''
\newblock in {\em Proc. IEEE Mobile Ad-hoc Sensor Netw.}, Shenzhen, China, Dec.
  2019, pp. 131--136.

\bibitem{huang2020sparse}
M.~Huang, L.~Huang, S.~Zhong, and P.~Zhang,
\newblock ``{UAV}-mounted mobile base station placement via sparse recovery,''
\newblock {\em IEEE Access}, vol. 8, pp. 71775--71781, 2020.

\bibitem{mach2021realistic}
P.~Mach, Z.~Becvar, and M.~Najla,
\newblock ``Power allocation, channel reuse, and positioning of flying base
  stations with realistic backhaul,''
\newblock {\em IEEE Internet Things J.}, pp. 1--1, 2021.

\bibitem{yin2021multiagent}
S.~Yin and F.~R. Yu,
\newblock ``Resource allocation and trajectory design in uav-aided cellular
  networks based on multi-agent reinforcement learning,''
\newblock {\em IEEE Internet Things J.}, pp. 1--1, 2021.

\bibitem{park2018formation}
S.~Park, K.~Kim, H.~Kim, and H.~Kim,
\newblock ``Formation control algorithm of multi-uav-based network
  infrastructure,''
\newblock {\em Applied Sciences}, vol. 8, no. 10, pp. 1740, 2018.

\bibitem{kim2018topology}
D.-Y. Kim and J.-W. Lee,
\newblock ``Integrated topology management in flying ad hoc networks: Topology
  construction and adjustment,''
\newblock {\em IEEE Access}, vol. 6, pp. 61196--61211, 2018.

\bibitem{alhourani2014urban}
A.~Al-Hourani, S.~Kandeepan, and A.~Jamalipour,
\newblock ``Modeling air-to-ground path loss for low altitude platforms in
  urban environments,''
\newblock in {\em IEEE Global Commun. Conf.}, 2014, pp. 2898--2904.

\bibitem{alhourani2014lap}
A.~Al-Hourani, S.~Kandeepan, and S.~Lardner,
\newblock ``Optimal {LAP} altitude for maximum coverage,''
\newblock {\em IEEE Wireless Commun. Lett.}, vol. 3, no. 6, pp. 569--572, 2014.

\bibitem{kalantari2016number}
E.~Kalantari, H.~Yanikomeroglu, and A.~Yongacoglu,
\newblock ``On the number and {3D} placement of drone base stations in wireless
  cellular networks,''
\newblock in {\em IEEE Vehicular Tech. Conf.}, 2016, pp. 1--6.

\bibitem{hammouti2019mechanism}
H.~El~Hammouti, M.~Benjillali, B.~Shihada, and M.-S. Alouini,
\newblock ``A distributed mechanism for joint {3D} placement and user
  association in {UAV}-assisted networks,''
\newblock in {\em IEEE Wireless Commun. Netw. Conf.}, Marrakech, Morocco, Apr.
  2019.

\bibitem{perabathini2019qos}
B.~Perabathini, K.~Tummuri, A.~Agrawal, and V.S. Varma,
\newblock ``Efficient {3D} placement of {UAVs with QoS Assurance in Ad Hoc
  Wireless Networks},''
\newblock in {\em Int. Conf. Comput. Commun. Netw.}, 2019, pp. 1--6.

\bibitem{liu2019deployment}
X.~Liu, Y.~Liu, and Y.~Chen,
\newblock ``Reinforcement learning in multiple-{UAV} networks: Deployment and
  movement design,''
\newblock {\em IEEE Trans. Veh. Tech.}, vol. 68, no. 8, pp. 8036--8049, 2019.

\bibitem{shehzad2021backhaul}
M.K. Shehzad, A.~Ahmad, S.A. Hassan, and H.~Jung,
\newblock ``Backhaul-aware intelligent positioning of {UAVs} and association of
  terrestrial base stations for fronthaul connectivity,''
\newblock {\em IEEE Trans. Netw. Sci. Eng.}, pp. 1--1, 2021.

\bibitem{qiu2020reinforcement}
J.~Qiu, J.~Lyu, and L.~Fu,
\newblock ``Placement optimization of aerial base stations with deep
  reinforcement learning,''
\newblock in {\em IEEE Int. Conf. Commun.}, 2020, pp. 1--6.

\bibitem{sabzehali2021orientation}
J.~Sabzehali, V.K. Shah, H.S. Dhillon, and J.H. Reed,
\newblock ``{3D} placement and orientation of {mmWave-based UAVs for Guaranteed
  LoS Coverage},''
\newblock {\em IEEE Wireless Commun. Letters}, pp. 1--1, 2021.

\bibitem{romero2022cartography}
D.~Romero and S.-J. Kim,
\newblock ``Radio map estimation: A data-driven approach to spectrum
  cartography,''
\newblock {\em IEEE Signal Process. Mag.}, vol. 39, no. 6, pp. 53--72, 2022.

\bibitem{alayafeki2008cartography}
A.~Alaya-Feki, S.~B. Jemaa, B.~Sayrac, P.~Houze, and E.~Moulines,
\newblock ``Informed spectrum usage in cognitive radio networks: Interference
  cartography,''
\newblock in {\em Proc. IEEE Int. Symp. Personal, Indoor Mobile Radio Commun.},
  Cannes, France, Sep. 2008, pp. 1--5.

\bibitem{shrestha2022surveying}
R.~Shrestha, D.~Romero, and S.~P. Chepuri,
\newblock ``Spectrum surveying: Active radio map estimation with autonomous
  {UAVs},''
\newblock {\em IEEE Trans. Wireless Commun.}, vol. 22, no. 1, pp. 627--641,
  2022.

\bibitem{teganya2019locationfree}
Y.~Teganya, D.~Romero, L.~M. Lopez-Ramos, and B.~Beferull-Lozano,
\newblock ``Location-free spectrum cartography,''
\newblock {\em IEEE Trans. Signal Process.}, vol. 67, no. 15, pp. 4013--4026,
  Aug. 2019.

\bibitem{teganya2020rme}
Y.~Teganya and D.~Romero,
\newblock ``Deep completion autoencoders for radio map estimation,''
\newblock {\em IEEE Trans. Wireless Commun.}, vol. 21, no. 3, pp. 1710--1724,
  2021.

\bibitem{zeng2021toward}
Y.~Zeng and X.~Xu,
\newblock ``Toward environment-aware 6g communications via channel knowledge
  map,''
\newblock {\em IEEE Wireless Communications}, vol. 28, no. 3, pp. 84--91, Jun.
  2021.

\bibitem{patwari2008nesh}
N.~Patwari and P.~Agrawal,
\newblock ``{NeSh}: a joint shadowing model for links in a multi-hop network,''
\newblock in {\em Proc. IEEE Int. Conf. Acoust., Speech, Signal Process.}, Las
  Vegas, NV, Mar. 2008, pp. 2873--2876.

\bibitem{patwari2008correlated}
N.~Patwari and P.~Agrawal,
\newblock ``Effects of correlated shadowing: Connectivity, localization, and
  {RF} tomography,''
\newblock in {\em Proc. Int. Conf. Info. Process. Sensor Networks}, St. Louis,
  MO, Apr. 2008, pp. 82--93.

\bibitem{boyd2011distributed}
S.~Boyd, N.~Parikh, E.~Chu, B.~Peleato, and J.~Eckstein,
\newblock ``Distributed optimization and statistical learning via the
  alternating direction method of multipliers,''
\newblock {\em Found. Trends Mach. Learn.}, vol. 3, no. 1, pp. 1--122, Jan.
  2011.

\bibitem{romero2023placementarxiv}
D.~Romero, P.~Q. Viet, and R.~Shrestha,
\newblock ``Aerial base station placement via propagation radio maps,''
\newblock {\em arXiv preprint arXiv:2301.04966}, 2023.

\bibitem{chen2021relay}
J.~Chen, U.~Mitra, and D.~Gesbert,
\newblock ``{3D} urban {UAV} relay placement: Linear complexity algorithm and
  analysis,''
\newblock {\em IEEE Trans. Wireless Commun.}, pp. 1--1, 2021.

\bibitem{mozaffari2016coverage}
M.~Mozaffari, W.~Saad, M.~Bennis, and M.~Debbah,
\newblock ``Efficient deployment of multiple unmanned aerial vehicles for
  optimal wireless coverage,''
\newblock {\em IEEE Commun. Lett.}, vol. 20, no. 8, pp. 1647--1650, 2016.

\bibitem{agrawal2009correlated}
P.~Agrawal and N.~Patwari,
\newblock ``Correlated link shadow fading in multi-hop wireless networks,''
\newblock {\em IEEE Trans. Wireless Commun.}, vol. 8, no. 9, pp. 4024--4036,
  Aug. 2009.

\bibitem{romero2018blind}
D.~Romero, D.~Lee, and G.~B. Giannakis,
\newblock ``Blind radio tomography,''
\newblock {\em IEEE Trans. Signal Process.}, vol. 66, no. 8, pp. 2055--2069,
  Jan. 2018.

\bibitem{kim2011cooperative}
S.-J. Kim, E.~Dall'Anese, and G.~B. Giannakis,
\newblock ``Cooperative spectrum sensing for cognitive radios using {K}riged
  {K}alman filtering,''
\newblock {\em IEEE J. Sel. Topics Signal Process.}, vol. 5, no. 1, pp. 24--36,
  Jun. 2010.

\bibitem{wilson2009regularization}
J.~Wilson, N.~Patwari, and O.~G. Vasquez,
\newblock ``Regularization methods for radio tomographic imaging,''
\newblock in {\em Virginia Tech Symp. Wireless Personal Commun.}, Blacksburg,
  VA, Jun. 2009.

\bibitem{hamilton2014modeling}
B.~R. Hamilton, X.~Ma, R.~J. Baxley, and S.~M. Matechik,
\newblock ``Propagation modeling for radio frequency tomography in wireless
  networks,''
\newblock {\em IEEE J. Sel. Topics Signal Process.}, vol. 8, no. 1, pp. 55--65,
  Feb. 2014.

\bibitem{imai2017survey}
T.~Imai,
\newblock ``A survey of efficient ray-tracing techniques for mobile radio
  propagation analysis,''
\newblock {\em IEICE Trans. Commun.}, vol. 100, no. 5, pp. 666--679, 2017.

\bibitem{kanso2009compressed}
M.~A. Kanso and M.~G. Rabbat,
\newblock ``Compressed {RF} tomography for wireless sensor networks:
  Centralized and decentralized approaches,''
\newblock in {\em Int. Conf. Distributed Comput. Sensor Syst.}, Marina del Rey,
  CA, 2009, Springer, pp. 173--186.

\bibitem{gutierrezestevez2021hybrid}
M.~A. Gutierrez-Estevez, M.~Kasparick, R.~L.~G. Cavalvante, and
  S.~Sta{\'n}czak,
\newblock ``Hybrid model and data driven algorithm for online learning of
  any-to-any path loss maps,''
\newblock {\em arXiv preprint arXiv:2107.06677}, 2021.

\bibitem{mitchell1990comparison}
J.R. Mitchell, P.~Dickof, and A.G. Law,
\newblock ``A comparison of line integral algorithms,''
\newblock {\em Comput. Physics}, vol. 4, no. 2, pp. 166--172, 1990.

\bibitem{thrun2008simultaneous}
S.~Thrun,
\newblock ``Simultaneous localization and mapping,''
\newblock in {\em Robotics and cognitive approaches to spatial mapping}, pp.
  13--41. Springer, 2008.

\bibitem{liu2022integrated}
F.~Liu, Y.~Cui, C.~Masouros, J.~Xu, T.~X. Han, Y.~C. Eldar, and S.~Buzzi,
\newblock ``Integrated sensing and communications: Toward dual-functional
  wireless networks for 6g and beyond,''
\newblock {\em IEEE J. Selected Areas in Commun.}, vol. 40, no. 6, pp.
  1728--1767, 2022.

\bibitem{gens1980complexity}
G.~Gens and E.~Levner,
\newblock ``Complexity of approximation algorithms for combinatorial problems:
  a survey,''
\newblock {\em ACM SIGACT News}, vol. 12, no. 3, pp. 52--65, 1980.

\bibitem{candes2008reweighted}
E.J. Candes, M.B. Wakin, and S.P. Boyd,
\newblock ``Enhancing sparsity by reweighted $\ell_1$ minimization,''
\newblock {\em J. Fourier Analysis App.}, vol. 14, no. 5, pp. 877--905, 2008.

\bibitem{lin2021admm}
T.~Lin, S.~Ma, Y.~Ye, and S.~Zhang,
\newblock ``An {ADMM}-based interior-point method for large-scale linear
  programming,''
\newblock {\em Optim. Methods Software}, vol. 36, no. 2-3, pp. 389--424, 2021.

\bibitem{boyd}
S.~Boyd and L.~Vandenberghe,
\newblock {\em Convex Optimization},
\newblock Cambridge University Press, Cambridge, UK, 2004.

\bibitem{yuan2006grouplasso}
M.~Yuan and Y.~Lin,
\newblock ``Model selection and estimation in regression with grouped
  variables,''
\newblock {\em J. Royal Statist. Soc.: Series B (Statist. Method.)}, vol. 68,
  no. 1, pp. 49--67, 2006.

\bibitem{bertsekas1999}
Dimitri~P Bertsekas,
\newblock {\em Nonlinear programming},
\newblock Athena scientific Belmont, 1999.

\bibitem{tardella2011linearprogramming}
F.~Tardella,
\newblock ``The fundamental theorem of linear programming: extensions and
  applications,''
\newblock {\em Optimization}, vol. 60, no. 1-2, pp. 283--301, 2011.

\bibitem{romero2015onlinesemiparametric}
D.~Romero, S.-J. Kim, and G.~B. Giannakis,
\newblock ``Stochastic semiparametric regression for spectrum cartography,''
\newblock in {\em Proc. IEEE Int. Workshop Comput. Advan. Multi-Sensor Adapt.
  Process.}, Cancun, Mexico, Dec. 2015, pp. 513--516.

\bibitem{romero2017spectrummaps}
D.~Romero, S-J. Kim, G.~B. Giannakis, and R.~López-Valcarce,
\newblock ``Learning power spectrum maps from quantized power measurements,''
\newblock {\em IEEE Trans. Signal Process.}, vol. 65, no. 10, pp. 2547--2560,
  May 2017.

\end{thebibliography}
\end{document}